\RequirePackage{fix-cm}
\documentclass[smallextended]{svjour3}       
\smartqed  

\usepackage[numbers]{natbib}

\usepackage[utf8]{inputenc}
\usepackage[dvipsnames]{xcolor}
\usepackage{graphicx}
\usepackage{times}
\usepackage[utf8]{inputenc}
\usepackage{mathtools}
\usepackage{amsfonts}
\usepackage[T1]{fontenc}
\usepackage[pdfpagelabels]{hyperref}
\usepackage{bm}
\usepackage{verbatim}
\usepackage{algorithm}
\usepackage{framed}
\usepackage{booktabs,array}
\usepackage{xcolor}
\usepackage{comment}
\usepackage{marginnote}
\usepackage{mathtools}
\usepackage{amssymb}
\usepackage{yhmath}
\usepackage{subcaption}
\usepackage{ulem}
\usepackage{cancel}
\usepackage{soul,xcolor}

\newcommand{\MD}[2][black]{{\textcolor{#1}{#2}}}

\newcommand{\MG}[2][black]{{\textcolor{#1}{#2}}}

\newcommand{\mbRd}{{\mathbb{R}^d}}

\newcommand{\thetab}{{\boldsymbol \theta}}

\newcommand{\xb}{\mathbf{x}}

\newcommand{\yb}{\mathbf{y}}

\renewcommand\cite{\citep}

\begin{document}

\title{\center Fractional Modeling in Action: \\A Survey of Nonlocal Models for Subsurface Transport, Turbulent Flows, and Anomalous Materials}
%
\titlerunning{
Fractional Modeling in Action
}      

\author{\center 
Jorge L. Suzuki \and Mamikon Gulian \and Mohsen Zayernouri \and Marta D'Elia$^*$}

\authorrunning{Suzuki, Gulian, Zayernouri, D'Elia}

\institute{J.L. Suzuki \at Department of Mechanical Engineering \& Computational Mathematics, Science and Engineering, Michigan State University, East Lansing, MI, USA\\
\email{suzukijo@msu.edu}
\and M. Gulian \at Center for Computing Research, Sandia National Laboratories, Albuquerque, NM, USA\\
\email{mgulian@sandia.gov}
\and M. Zayernouri \at Department of Mechanical Engineering \& Statistics and Probability, Michigan State University, East Lansing, MI, USA\\
\email{zayern@msu.edu}
\and M. D'Elia, corresponding author \at Computational Science and Analysis, Sandia National Laboratories, Livermore, CA, USA\\
\email{mdelia@sandia.gov}}

\date{Received: date / Accepted: date}

\maketitle

\begin{abstract}
Modeling of phenomena such as anomalous transport via fractional-order differential equations has been established as an effective alternative to partial differential equations, due to the inherent ability to describe large-scale behavior with greater efficiency than fully-resolved classical models. In this review article, we first provide a broad overview of fractional-order derivatives with a clear emphasis on the stochastic processes that underlie their use. We then survey three exemplary application areas -- subsurface transport, turbulence, and anomalous materials -- in which fractional-order differential equations provide accurate and predictive models. For each area, we report on the evidence of anomalous behavior that justifies the use of fractional-order models, and survey both foundational models as well as more expressive state-of-the-art models. We also propose avenues for future research, including more advanced and physically sound models, as well as tools for calibration and discovery of fractional-order models.
\end{abstract}


\section{Introduction}\label{sec:intro}

Understanding and applying the theory of \textit{anomalous transport} opens up rich fields of study in science and engineering, transforming our perspective and facilitating extraordinary discoveries that would not be possible otherwise. 
This class of phenomena refers to fascinating and widespread\footnote{\MG{The adjective ``anomalous'' means ``not normal'' -- an apt description of such processes in a statistical context, as discussed in Section \ref{sec:classification}. However, it belies the very widespread nature of such processes, an irony that \citet{klafter2005anomalous} and \citet{sancho2004diffusion} point out with the statement ``anomalous is normal!''}} processes that, viewed at appropriate scales, exhibit non-Markovian long-term memory effects, non-Fickian long-range interactions, nonergodic statistics, and non-equilibrium dynamics \cite{klages2008anomalous}. Anomalous transport is observed in a wide variety of complex, multi-scale, and multi-physics systems such as subdiffusion and superdiffusion in porous media, kinetic plasma turbulence, aging of polymers, glassy materials, amorphous semiconductors, biological cells, heterogeneous tissues, and disordered media \cite{klafter2005anomalous, west2016complexity,Wong2004,Sollich1998}. The crucial point that prompts this work is that conventional mathematical models cannot describe such processes in a succinct, compact way that directly expresses their anomalous and nonlocal character. 

This work is founded on the use of \textit{fractional-order partial differential equations} (FPDEs), which seamlessly generalize standard PDEs of integer order to real-valued order. In practice, FPDEs appear within tractable mathematical models for anomalous transport, ranging from complex fluids to non-Newtonian rheology and the design of aging materials \cite{West2003, klages2008anomalous, meerschaert2019stochastic, Jaishankar2013, Jaishankar2014}, but also in modeling transport phenomena when rates of change in the quantity of interest depend on space or time. In this context, FPDEs with ``variable orders'' can be exploited in diverse physical and biological applications \cite{patnaik2020applications,patnaik2021variable,zhao2015second} to capture transitions between different transport regimes. Moreover, even classical long-standing issues such as monotonicity, anisotropy, and multi-fractal scaling laws in turbulence can be reformulated and reinterpreted in the context of fractional calculus and probability theory. FPDEs therefore emerge as an \textit{expressive} approach to modeling such physics, transforming the current practice in mathematical modeling and giving rise to a new generation of flexible, high-fidelity, and direct approaches \cite{atanackovic2014, fallahgoul2016fractional, west2017Nature_Patterns,west2016complexity}.

In this review article, we focus on three important applications of FPDEs, reporting the scientific evidence of how and why fractional modeling naturally emerges in each case, along with a review of selected nonlocal mathematical models that have been proposed.
For brevity, throughout this article we use the term ``fractional'' to mean ``fractional-order''. Despite conflicting with the most common usage of the adjective ``fractional'' in the English language, this is standard in the literature; thus, fractional-order derivatives are referred to as ``fractional derivatives'' and fractional-order models as ``fractional models''. 

\vspace{1ex}
\noindent \textbf{Anomalous Subsurface Transport (Section \ref{sec:subsurface})} The accurate prediction at large scales of contaminant transport in both surface and subsurface water is fundamental for the management of water resources and critical for environmental safety. However, the explicit description of the systems where transport takes place is extremely challenging, especially at large scales, due to the complexity the medium. Such media feature heterogeneities that are either difficult or impossible to observe, and hence cannot be described with certainty at all relevant scales and locations. Moreover, even when the environment's microstructure can be captured, numerical simulations of appropriate PDE models such as systems of advection-diffusion equations may be infeasibly expensive if conducted at fully-resolved small scales \citep{Sun2020}. In fact, the same types of equations that are accurate at small scales \textit{do not} extrapolate and predict solutes' behavior at larger scales, due to the appearance of ``anomalous'', or ``non-Fickian'' behavior \cite{Benson2000,Benson2001,Levy2003,Neuman2009}. At large scales, FPDEs are called for. 

\vspace{1ex}
\noindent \textbf{Turbulent Flows (Section \ref{sec:turbulence})} Turbulence ``remembers'' and is fundamentally nonlocal. Coherent motions and ``turbulence spots'' structures inherently give rise to intermittent signals with sharp peaks, heavy-skirts, and skewed distributions of velocity increments \citep{batchelor1953theory, frisch1995turbulence}, manifesting the non-Markovian, non-Fickian nature of turbulence. This suggests that nonlocal interactions cannot be ruled out in the physics of turbulence \cite{davidson2015turbulence}. In addition to such an inherent nonlocality, filtering the Navier-Stokes and energy equations in the corresponding large eddy simulation (LES) of turbulent flows and scalar turbulence, in which large-scale motions are ``resolved'' and only the small scales are ``modeled'', would make the existing nonlocality in the corresponding subgrid stochastic processes (i.e., turbulent fluctuations) even more pronounced \cite{Scalar_FSGS,FSGS, Zayernouri_2021_FracLES}. This requires the development of new modeling paradigms in addition to new statistical measures that can meticulously highlight the nonlocal character of turbulence and their absence in the common turbulence modeling practice. 

\vspace{1ex}
\noindent \textbf{Anomalous Materials/Rheology (Section \ref{sec:material})}  Accurate modeling of the evolution of material response and failure across multiple time and length scales is essential for life-cycle prediction and design of new materials. While the mechanical behavior of a number of standard engineering materials (e.g., metals, polymers, rubbers) is quite well-understood, a significant modeling effort still needs to be conducted for complex materials, where microstructure heterogeneities, randomness and small scale physical mechanisms \cite{Wong2004,zapperi1997plasticity} (e.g. trapping effects and collective behavior) lead to non-standard and, at times, counter-intuitive responses. Two examples are bio-tissues and natural materials, e.g. biopolymers, which are multi-functional products of millions of years of evolution, locally optimized for their hosts and environment, and constrained by a limited set of building blocks and available resources \cite{Imbeni2005,Wegst20151}. These materials possess unprecedented properties at low densities, especially due to their hierarchical and multi-scale structure, leading to a wide spectrum of behaviors, such as power-law viscoelasticity , viscoplastic strains under hysteresis loading, damage and failure, fractal avalanche ruptures and self-healing mechanisms \cite{Stamenovic2007,bonfanti2020fractional,Bonadkar2016,Bonamy2008,Richeton2005,Wegst20151}.

\subsection{Outline of the article}\label{sec:outline}

Before describing each of the aforementioned applications, we review the foundations of fractional calculus. \MG{We} classify fractional models via their connection with the underlying stochastic processes that serve as the statistical backbone of fractional modeling. The organization of the rest of the review article is as follows: Sections \ref{sec:subsurface}, \ref{sec:turbulence}, and \ref{sec:material} are dedicated to subsurface transport, turbulent flows, and anomalous materials, respectively. Each section has the same structure: first, we motivate the need for fractional modeling and provide results or tools necessary for a full understanding of the section. Next, we provide evidence of fractional behavior, reporting state-of-the-art results that highlight the improved accuracy of FPDEs as opposed to classical PDEs. Then, the core of each section is a description of past and current models, with some insights on discretization techniques currently in use. At the end of each section, a paragraph on future directions gives our perspective on fruitful research directions in each area.

\section{An overview of fractional derivatives}

\subsection{Classification of fractional derivatives and models}\label{sec:classification}
We introduce and classify the most commonly used fractional-order differential operators in the context of diffusion models based on random walks. 
For simplicity, we restrict our discussion to one spatial dimension except for a few remarks in which the extension to higher dimensions is touched upon. 

{To avoid mathematical \MG{intricacies}, we discuss stochastic processes in terms of their discretizations, thinking of them as sequences of random variables $X_{n \Delta t}$ for time step $\Delta t > 0$ and integer $n$ \MG{that are} defined as cumulative sums of increments. Strictly speaking, FPDEs govern the statistical properties of \textit{continuous-time} random walks, which are appropriate scaling limits or long-time limits of the discrete random walks, limits in which $n$ becomes large relative to $\Delta t$ \citep{Meerschaert2012}. 
However, the rigorous definition of such stochastic processes requires significant excursions into probability theory; this is true even for the classical case of Brownian motion \citep{rogers1994diffusions, rogers2000diffusions}. Thus, while not entirely precise, in introducing fractional operators we characterize the related process in their discretized form, providing references where rigorous definitions of the process, as well as proofs of convergence of the discretization to the continuous-time process in appropriate limits, are given. }

\begin{figure}[t]
\includegraphics[width=\textwidth]{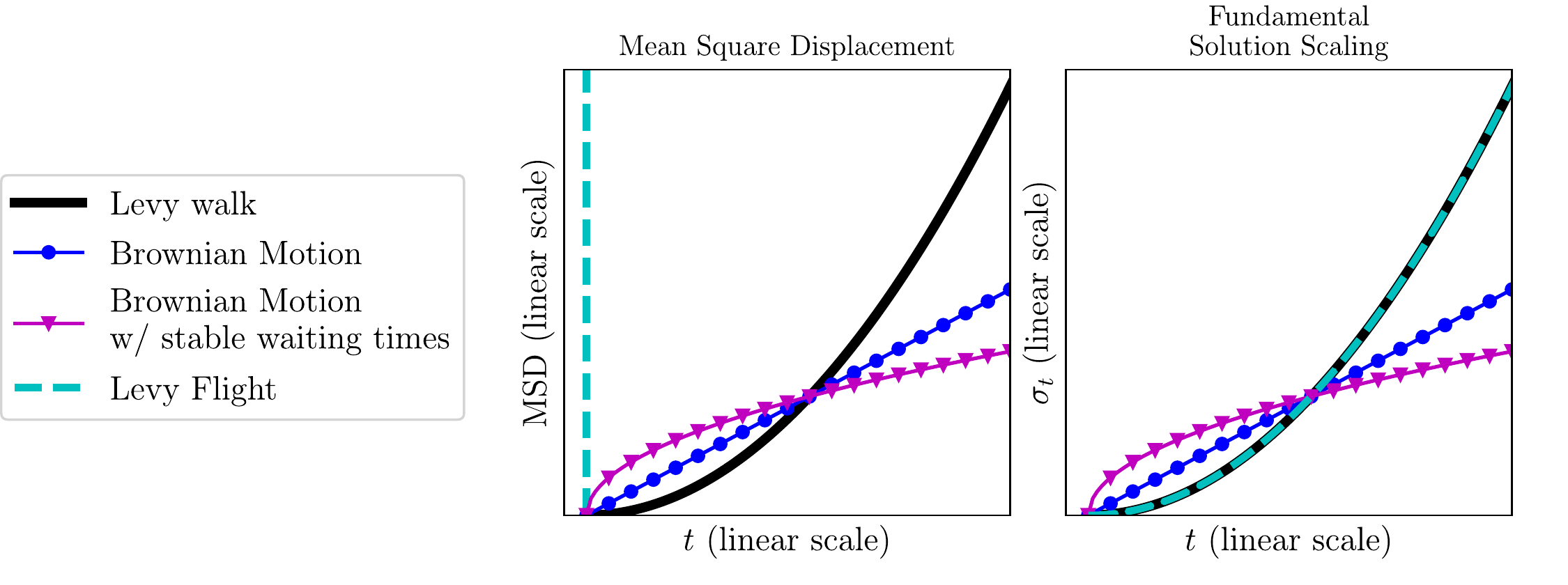}
\caption{Comparison of normal diffusion, superdiffusion, and subdiffusion via mean-square displacement (MSD) of the particle models and the scaling-in-time of the fundamental solution the diffusion equations governing the density functions. \MG{These are important characteristic properties that distinguish classes of diffusion.} Brownian motion exhibits both an MSD and scaling factor that are linear in time. Superdiffusive L\'evy flight exhibits infinite MSD, and a fundamental solution scaling factor $t^{\alpha}$ for $\alpha > 1$, while the superdiffusive L\'evy walk exhibits the same scaling of the fundamental solution as well as a finite MSD that scales as $t^{\alpha}$ for long times. The subdiffusive Brownian Motion with waiting times exhibits sublinear MSD and fundamental solution scaling, proportional to $t^{1/\alpha}$ for $\alpha < 1$. }
\label{fig:motion-comparison}
\end{figure}
 
\subsubsection{Normal, or Fickian, diffusion}
The connection between Brownian motion $B_t$ and the classical diffusion equation was studied in seminal works by \citet{bachelier1900theorie}, \citet{einstein1905motion}, and \citet{smoluchowski}. The diffusion equation is posed in an initial value problem,
\begin{equation}\label{eq:normal_diffusion_equation}
\left\{
\begin{split}
\frac{\partial u}{\partial t}(x,t) &= k^2 \Delta u(x,t), & \quad x \in \mathbb{R}, t > 0, \\
u(x, t = 0) &= u_0 (x), & \quad x \in \mathbb{R},
\end{split}
\right.
\end{equation}
in which $k^2 > 0$ is the diffusion coefficient and $\Delta u = \partial^2 u / \partial x^2$ denotes the Laplacian. 
Brownian motion $B_t$ is a continuous-time stochastic process defined for $t \ge 0$, which when discretized in time steps of size $\Delta t$ has the property that $B_{t = 0} = 0$ and
\begin{equation}\label{eq:brownian_discretization}
B_{t+\Delta t} = B_{t} + \Delta B; \quad  \Delta B \sim \mathcal{N}(\mu = 0, \sigma = k\sqrt{\Delta t}). 
\end{equation}
The above notation indicates that the increment $\Delta B$ at each time step of $\Delta t$ is drawn from a normal distribution $\mathcal{N}(0,k\sqrt{\Delta t})$ with mean $\mu = 0$ and standard deviation $\sigma = k\sqrt{\Delta t}$. This has probability density function
\begin{equation}\label{eq:normal_density}
p_{\mathcal{N}}(x; \mu,\sigma) = \frac{1}{\sigma\sqrt{2 \pi}}
e^{-\frac{1}{2} \left( \frac{x-\mu}{\sigma} \right)^2}
\end{equation}
The rule \eqref{eq:brownian_discretization} for sampling a path of $B_t$ at times $m \Delta t$, $m = 0, 1, 2, ...$ is an example of a discrete stochastic differential equation (SDE), and is referred to as the Euler-Maruyama 
discretization\footnote{Another frequently used discrete random walk that leads to Brownian motion simply involves steps of fixed length to the left or right with probability $1/2$ each; see \citet{lawler2010random}. In the long-time limit, all such discrete walks that draw increments from a finite-variance distribution lead to Brownian motion, due to the central limit theorem \cite{zaburdaev2015levy}} of Brownian motion \citep{kloeden1992stochastic}. 

\begin{figure}[t]
\includegraphics[width=\textwidth]{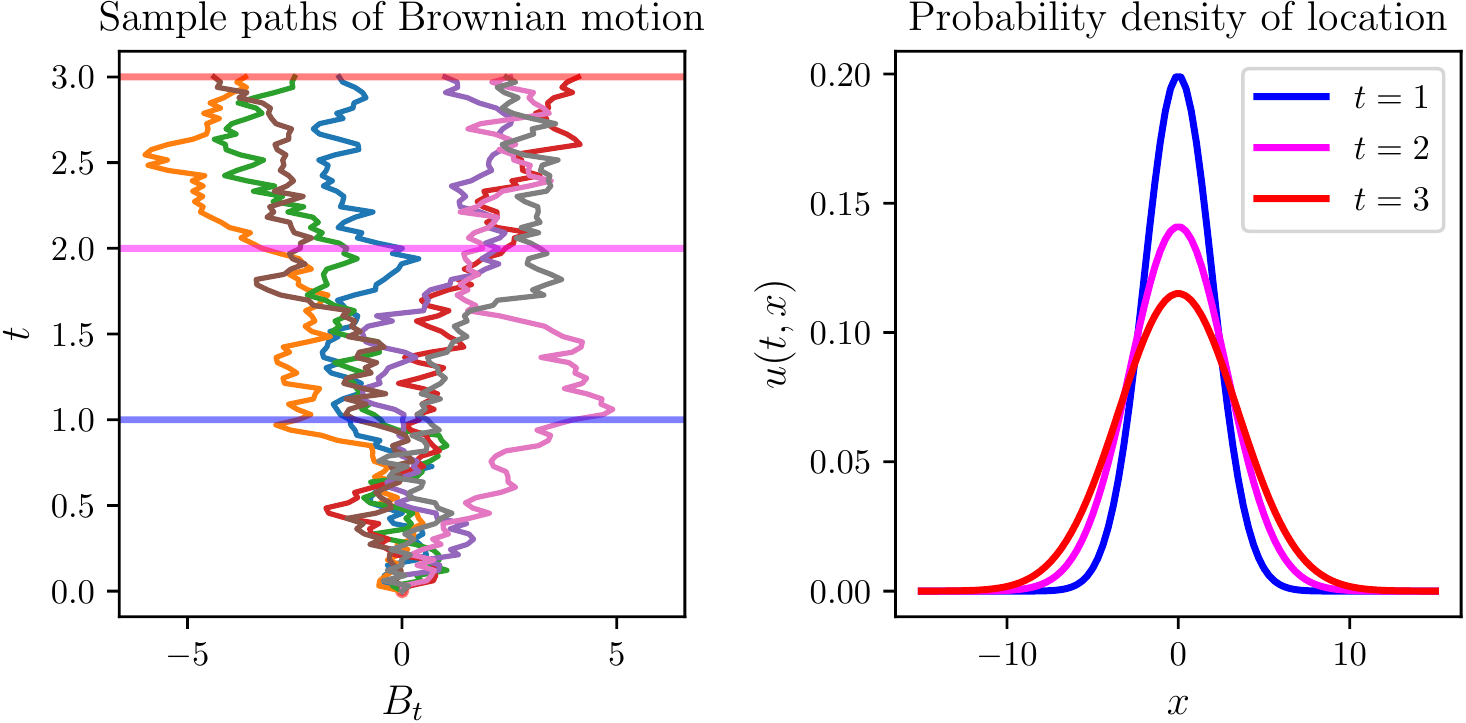}
\caption{\textit{(Left)} Eight independent sample paths of Brownian motion representing the path of a particle starting at the origin and stepping according to the rule \eqref{eq:brownian_discretization}. \textit{(Right)} For $t= 1,2,3$, the probability density of the location of the particle, i.e. the fundamental solution to the classical heat equation \eqref{eq:normal_diffusion_equation}.
}
\label{fig:bm_paths_density}
\end{figure}

The discrete process $B_{ m \Delta t}$ should be thought of as tracing a path in $\mathbb{R}$ of a particle undergoing ``jumps'' in a random direction at time intervals of size $\Delta t$. At each time $t$, the position $B_{t}$ of the particle is a random variable. 
It can be shown that the paths of the continuous-time process $B_t$ are almost surely continuous in time \cite{revuz2013continuous}. From \eqref{eq:brownian_discretization} and the central limit theorem, it follows that Brownian motion satisfies the scaling property
\begin{equation}\label{eq:brownian_scaling}
\langle B_t^2 \rangle = 2 k^2t,
\end{equation} 
where the left-hand side denotes the variance or second moment of the random variable $B_t$\MG{; see Figure \ref{fig:motion-comparison}.}
Given an initial distribution of particles $u_0(x)$ in $\mathbb{R}$ which then undergo Brownian motion, the distribution $u(x,t)$ of particles in $\mathbb{R}$ is governed by \eqref{eq:normal_diffusion_equation}. In other words, diffusing particles described at a microscopic scale by Brownian motion, i.e. by \eqref{eq:brownian_discretization} in discrete time, have their distribution in space -- a macroscopic property -- governed by the heat equation \citep[Section 1.1]{Meerschaert2012}. This is illustrated in Figure \ref{fig:bm_paths_density}. 

The consistency between this macroscopic description and the microscopic model is illustrated by scaling properties.
A necessary property of \MG{the} Brownian motion model is the second-moment condition \eqref{eq:brownian_scaling}, which states that on average, particles travel a distance $k\sqrt{t}$ from their initial position after time $t$. This is reflected in the fact that the solution of \eqref{eq:normal_diffusion_equation} with initial condition $u_0(x) = \delta_0(x)$ is
\begin{equation}\label{eq:classical_diffusion_fundamental_solution}
u(x,t) = \frac{1}{\sqrt{4\pi t}} e^{-{x^2}/{4k^2t}},
\end{equation}
which is the normal density \eqref{eq:normal_density} with standard deviation $k\sqrt{t}$. Note that this solution has the property that
\begin{equation}
u(x,t_2) = \left(\frac{t_2}{t_1}\right)^{-{1}/{2}} \, u\left(\frac{x}{\left({t_2}/{t_1}\right)^{-{1}/{2}}},t_1\right), \quad t_2 > t_1 > 0.
\end{equation}
Thus, the distribution of plume of particles in this diffusion model spreads out as $(t_2/t_1)^{1/2}$ as time elapses from $t_1$ to $t_2$, consistent with \eqref{eq:brownian_scaling}.
\MG{This scaling of the fundamental solution is illustrated in Figure \ref{fig:motion-comparison}}.

The model for normal diffusion reviewed here is also referred to as Fickian diffusion. The heat equation \eqref{eq:normal_diffusion_equation} can be derived from the mass conservation with flux term $J$,
\begin{equation}
\frac{\partial u}{\partial t} 
+
\frac{\partial J}{\partial x}
=0
\end{equation}
under Fick's law $J = \nabla u$. As discussed by \citet{Schumer2001}, the fractional diffusion equations we introduce below follow from mass conservation with non-Fickian fluxes. 

\subsubsection{$\alpha$-stable L\'evy flights and the fractional Laplacian}
\label{sec:frac_lapl}
Many important systems exhibit diffusive behavior, but do not satisfy the scaling property \eqref{eq:brownian_scaling} \citep{klafter2005anomalous}. This type of diffusion is referred to as \textit{anomalous diffusion}, as it cannot be described by \eqref{eq:brownian_discretization} with normally distributed increments. We desire a microscopic model that generalizes Brownian motion $B_t$, and a corresponding macroscopic model that generalizes the diffusion equation \eqref{eq:normal_diffusion_equation}. The first model we propose remains in the framework of a discrete SDE with independent identically distributed (i.i.d.) increments, 
\begin{equation}\label{eq:iid_rule}
X_{t+\Delta t} = X_{t} + \Delta X,  \quad
X_0 = 0,
\end{equation}
but the increments $\Delta X$ are no longer drawn from a normal distribution. It follows from the central limit theorem that the only way to obtain a microscopic model in this framework that is statistically distinct from $B_t$, i.e., not equivalent in distribution, is to draw step sizes from a probability density function with infinite variance \citep{zaburdaev2015levy, Meerschaert2012}. 

We introduce the isotropic $\alpha$-stable random variable $S_\alpha(\gamma,\sigma,\mu)$. This family of random variables is defined\footnote{Several parametrizations of the $\alpha$-stable characteristic function exist. The parametrization \eqref{eq:alpha_stable_characteristic} is due to \citet{samorodnitsky1994}. See \citet{nolan2020univariate,nolan1998parameterizations} for discussions of alternate forms.}  most simply by their characteristic function. For a general random variable $X$, the characteristic function $\varphi_X$ is the Fourier transform of the probability density function $p_X$ of $X$, i.e., 
\begin{equation}
\varphi_X(\xi) = \int e^{i\xi x} p_X(x) dx.
\end{equation}
Thus, the characteristic function of the normal random variable \MG{$\mathcal{N}(\mu,\sigma)$} is $e^{i\xi\mu - \sigma^2 \xi^2/2}$. \MG{Generalizing this, }the $\alpha$-stable random variable \MG{$S_\alpha(\gamma,\sigma,\mu)$} has characteristic function
\begin{equation}\label{eq:alpha_stable_characteristic}
\varphi_\alpha(\xi;\gamma,\sigma,\mu)=e^{i\xi\mu-|\sigma \xi|^\alpha (1 - i\gamma \text{sgn}(\xi) \Phi ) }, 
\end{equation}
where
\begin{equation}
\Phi = 
\begin{cases}
\tan\left(\frac{\pi \alpha}{2}\right) \text{ if $\alpha \neq 1$, } \\
-\frac{2}{\pi} \log(|\sigma\xi|) \text{ if $\alpha = 1$}.
\end{cases}
\end{equation}
The parameter $\alpha\in (0,2]$ is referred to as the stability parameter of the distribution, $\mu \in \mathbb{R}$ as the center, $\gamma\in [-1,1]$ as the skewness, and $\sigma \in (0,\infty)$ as the scale. 
The \textit{isotropic} or \textit{symmetric} $\alpha$-stable distribution $S_\alpha(\gamma = 0, \sigma, \mu = 0)$ therefore has characteristic function
\begin{equation}\label{eq:iso_alpha_stable_characteristic}
\varphi_\alpha(\xi;\gamma = 0,\sigma,\mu = 0)=e^{-\sigma|\xi|^\alpha }, 
\end{equation}
generalizing the characteristic function of the normal distribution with mean $\mu = 0$ and standard deviation $\sigma/\sqrt{2}$ and reducing to it when $\alpha = 2$. 
By definition, the probability density function of $S_\alpha(\gamma,\sigma,\mu)$ can be written
\begin{equation}\label{eq:iso_alpha_stable_density}
f_\alpha(x;\gamma,\sigma,\mu) = \int e^{i \xi x} \varphi_\alpha(\xi; \gamma,\sigma,\mu) d\xi. 
\end{equation}
In general, the $\alpha$-stable density does not admit a closed-form expression\footnote{Special cases are $\alpha = 2$ corresponding to the normal distribution, $\alpha = 1$ and $\gamma = 0$ corresponding to the Cauchy distribution, and $\alpha = 1/2$ and $\gamma = 1$ corresponding to the L\'evy distribution.}, but in the symmetric case where $\gamma = \mu = 0$, it has the property that
\begin{equation}\label{eq:pareto}
f_\alpha(x;\gamma = 0,\sigma,\mu = 0) \sim \frac{1}{|\sigma x|^{1+\alpha}} \text{ for large } x, \quad 0 < \alpha < 2, 
\end{equation}
as discussed in, e.g., \citet{nolan2020univariate} or 
\citet{tankov2003financial}.
In other words, the density exhibits \textit{Paretian} or power-law tails. This is in contrast to the rapidly decaying square-exponential tails of the normal distribution. In many settings, such tails are informally referred to as being examples of heavy or fat tails \citep{adler1998practical, haas2009financial}. 

\begin{figure}[t]
\includegraphics[width=\textwidth]{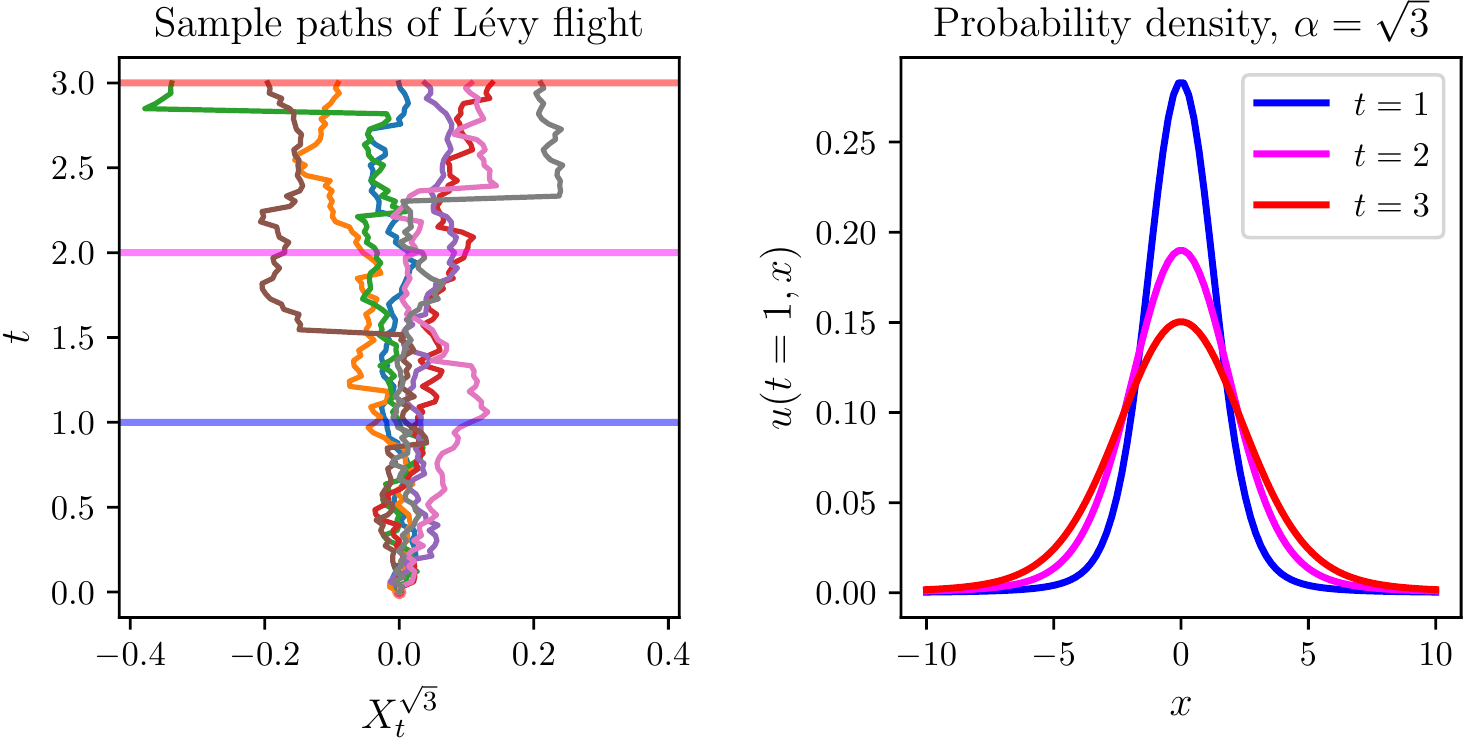}
\caption{\textit{(Left)} Eight independent sample paths of symmetric $\alpha$-stable L\'evy flight with $\alpha = \sqrt{3}$ representing the path of particle starting at the origin and stepping according to the rule \eqref{eq:isotropic_levy_flight_discretization}. \textit{(Right)} For $t= 1,2,3$, the probability density of the location of the particle given by \eqref{eq:iso_alpha_stable_density}, i.e. the fundamental solution to the fractional diffusion equation \eqref{eq:fractional_laplacian_equation}. Compared to Figure \ref{fig:bm_paths_density}, note that despite some qualitative similarity between the shapes of the density functions, the presence of long jumps signifies a striking difference between the paths of a particle undergoing L\'evy flight versus Brownian motion.}
\label{fig:symmeric_stable_paths_density}
\end{figure}

\begin{figure}[t]
\includegraphics[width=\textwidth]{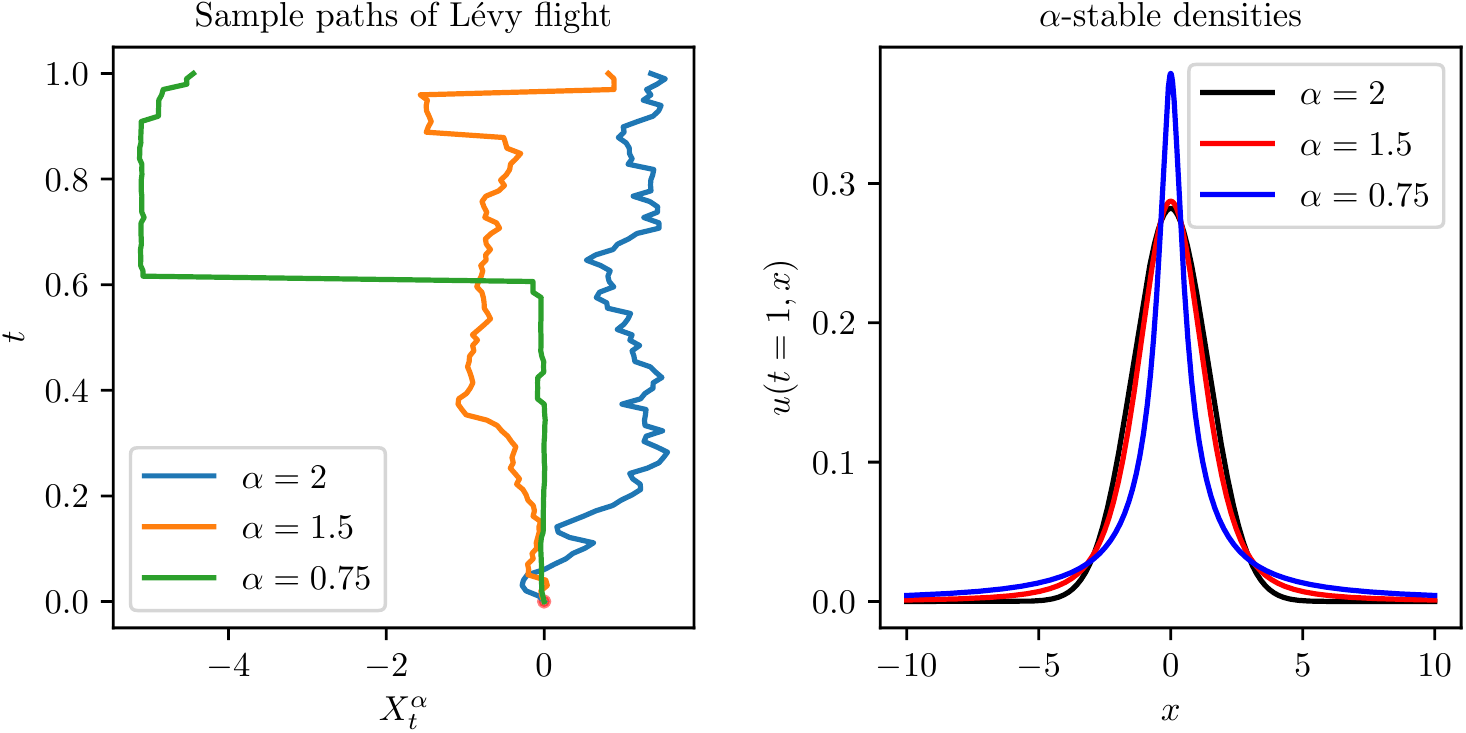}
\caption{The seemingly innocuous heavy tails of the $\alpha$-stable density, signifying non-vanishing probability of long jumps, are responsible for the striking properties of $\alpha$-stable L\'evy flights. As $\alpha$ decreases from $2$, more mass in the middle region of the density is lost and is transferred towards the tails and the center, so that the relative probability of very small movements and very long movements increases (\textit{right}). This is evident in the sample paths of the process (\textit{left}).}
\label{fig:symmeric_stable_paths_comparison}
\end{figure}

\begin{figure}[t]
\centering
\includegraphics[width=0.75\textwidth]{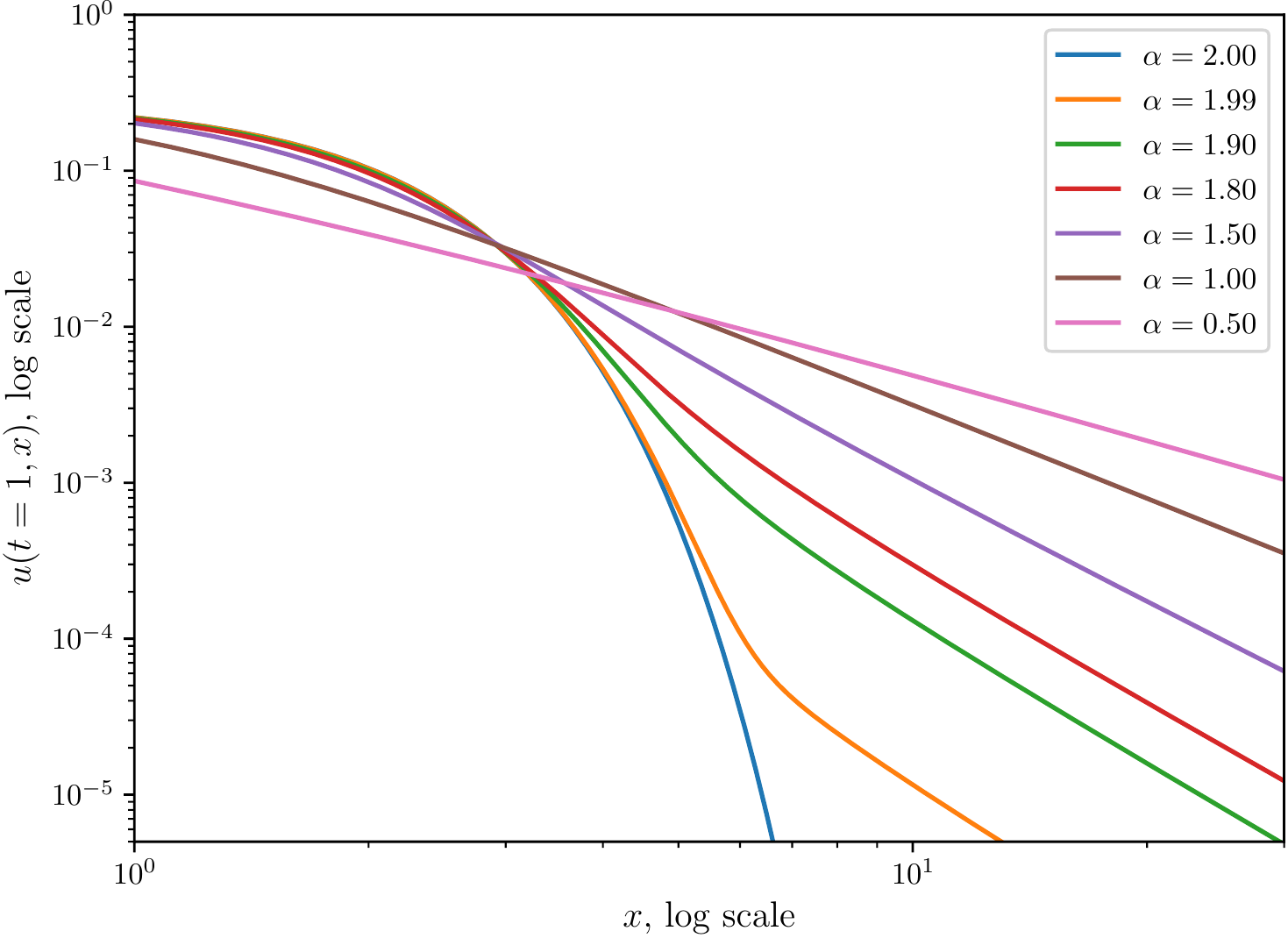}
\caption{A plot of the $\alpha$-stable densities in a log-log scale that illustrates the tail behavior asserted in \eqref{eq:pareto}. While $\alpha$-stable densities do not have a closed form expression for all $x$, their simple, asymptotic inverse power-law behavior is an important heuristic.}
\label{fig:symmeric_stable_loglog}
\end{figure}

Using the isotropic distribution introduced above, we introduce the isotropic $\alpha$-stable L\'evy flight $X^\alpha_t$ by providing the corresponding discrete stochastic process. This is given for $t =k\Delta t$ with integer $k$ by $X^\alpha_0 = 0$ and the rule \citep{Meerschaert2012}
\begin{equation}\label{eq:isotropic_levy_flight_discretization}
X_{t+\Delta t}^{\alpha} = X^{\alpha}_{t} + \Delta X^\alpha; \quad  \Delta X^\alpha \sim S_\alpha(\gamma = 0,\sigma = k(\Delta t)^{1/\alpha}, \mu = 0).
\end{equation}
The continuous-time stochastic process $X^\alpha_t$ for $t \ge 0$ can be thought of as a scaling limit as $\Delta t \rightarrow 0$ of the above random walk, and enjoys several theoretical properties such as stability and an extended central limit theorem \citep{Meerschaert2012, meerschaert2001limit}. However, has the property that for $\alpha < 2$, the paths of $X^\alpha_t$ are almost surely discontinuous, in contrast to Brownian motion -- hence the name L\'evy ``flight''.
Given an initial distribution $u_0(x)$ of particles in $\mathbb{R}$ which undergo $\alpha$-stable L\'evy flight, the evolution of the distribution $u(x,t)$ for $t > 0$ is governed by the \textit{space-fractional} diffusion equation \citep[Section 1.2]{Meerschaert2012}
\begin{align}
\label{eq:fractional_laplacian_equation}
\begin{split}
\frac{\partial u}{\partial t}(x,t) &= -k^\alpha (-\Delta)^{\alpha/2} u(x,t) \\
u(x, t = 0) &= u_0 (x),
\end{split} 
\end{align}
as illustrated in Figure \ref{fig:symmeric_stable_paths_density}. 
The fractional negative Laplacian $(-\Delta)^{\alpha/2}$ is defined for $0 < \alpha < 2$ and for any dimension $d$ as
\begin{equation}\label{eq:fractional_laplacian_definition}
\MG{
(-\Delta)^{\alpha/2} u(\xb) = C_{d,\alpha}
\text{ p.v.}\int_{\mathbb{R}^d} \frac{u(\xb) - u(\yb)}{|\xb-\yb|^{d+\alpha}} d\yb, \quad \xb \in \mathbb{R}^d,}
\end{equation}
with
\begin{equation}
C_{d,\alpha}=
\frac{4^{{\alpha/2}} \Gamma\left({\alpha/2}+\frac{d}{2}\right)}
{\pi^{d/2}|\Gamma(-{\alpha/2})|};
\end{equation}
see \citet{lischke2020fractional}.
We have defined this operator in any dimension for future reference, although our present discussion only requires the case $d = 1$. 
Perhaps the simplest characterization of the fractional Laplacian is the Fourier representation, 
\begin{equation}\MG{
\mathcal{F} \left[ (-\Delta)^{\alpha/2} u \right] (\bm{\xi})
=
|\bm{\xi}|^{\alpha} \mathcal{F}[u](\bm{\xi}),
\quad
\bm{\xi} \in \mathbb{R}^d}.
\end{equation}

The simplest case of \eqref{eq:fractional_laplacian_equation} is the initial condition $u_0(x) = \delta_0(x)$, in which case the solution is
\begin{equation}\label{eq:symmetric_fundamental_solution}
u(x,t) = f_\alpha(x;\gamma = 0,\sigma = k t^{1/\alpha},\mu = 0). 
\end{equation}
This is known as the fundamental solution. 
Although this solution cannot be written in closed form, it satisfies
\begin{equation}
u(x,t_2) = \left(\frac{t_2}{t_1}\right)^{-{1}/{\alpha}} \, u\left(\frac{x}{\left({t_2}/{t_1}\right)^{-{1}/{\alpha}}}, t_1\right), \quad t_2 > t_1 > 0,
\end{equation}
as shown in \citet[Section 1.2]{Meerschaert2012}.
This illustrates that a plume of particles undergoing isotropic $\alpha$-stable L\'evy flight spreads by a factor of $(t_2/t_1)^{1/\alpha}$ as time elapses from $t_1$ to $t_2$, a faster rate when $\alpha < 2$ than the normal rate
$t^{1/2}$ . Thus, $\alpha$-stable L\'evy flight is an example of \textit{superdiffusion}. The dependence of the above solution as well as sample paths on $\alpha$ is shown in Figure \ref{fig:symmeric_stable_paths_comparison}.

However, since $\alpha > 0$, the tail behavior of the isotropic $\alpha$-stable density implies that the second moment of $X^\alpha_t$ diverges for $\alpha < 2$, 
\begin{equation}
\langle X^\alpha_t \rangle^2 = \infty,\quad 0 < \alpha < 2. 
\end{equation}
with the first moment (the mean) diverging also when $\alpha \le 1$ \cite{nolan2020univariate,zaburdaev2015levy}. This implies that the variance of $\alpha$-stable motion is not a useful statistic for parameterizing $\alpha$-stable L\'evy flight; it bears no useful relationship to $\alpha$. This aspect can be tackled in several ways, motivating the introduction of further fractional-order operators, such as tempered operators and fractional material derivatives \MG{discussed below}.   

We point out several important properties of the fractional Laplacian. From the definition \eqref{eq:fractional_laplacian_definition}, it is clear that $(-\Delta)^{\alpha/2} c = 0$, $c$ being a constant. The fractional Laplacian also satisfies the semigroup property $(-\Delta)^{\alpha/2} (-\Delta)^{\beta/2} = (-\Delta)^{(\alpha+\beta)/2}$ \citep{Samko1993}. However, one property that is apparent from the definition is that, unlike integer-order derivatives, the fractional Laplacian is a \textit{nonlocal} operator, i.e. the value of $(-\Delta)^{\alpha/2} u (\MG{\xb})$ depends on the values of $u$ in all of $\mathbb{R}$ (or $\mathbb{R}^d$, for $d>1$). In contrast, the value of any integer-order derivative of $u$ at $\MG{\xb}$ depends only on the values of $u$ in an infinitesimal neighborhood of $\MG{\xb}$.   

\subsubsection{The Riemann-Liouville fractional derivatives and asymmetric $\alpha$-stable L\'evy flight}
\label{sec:RL}
The fractional Laplacian \eqref{eq:fractional_laplacian_definition} was introduced in the previous section as a symmetric or rotation invariant operator for describing the symmetric or isotropic $\alpha$-stable L\'evy flight. \MG{That} model introduced a stability parameter $0 < \alpha \le 2$ allowing it to generalize normal diffusion, with the scale $\sigma$ and center $\mu$ playing similar roles as the standard deviation and mean of the normal distribution. However, the stable distribution also allows for a skewness parameter $\gamma \in [-1,1]$, with $\beta = 0$ in the symmetric case, which has no analogue in the normal distribution or for Brownian motion. This is due to the central limit theorem, which states that the use of any finite-variance distribution for the i.i.d. increments $\Delta X$ in \eqref{eq:iid_rule}, no matter how asymmetric, leads to $X_t$ being normally distributed, \MG{and therefore} necessarily symmetric about the mean. 
In this section, we introduce the one-sided Riemann-Liouville fractional derivatives as appropriate operators for modeling asymmetric $\alpha$-stably L\'evy flights, which are defined by 
\eqref{eq:isotropic_levy_flight_discretization} with  $\Delta X^\alpha \sim S_\alpha(\gamma,\sigma = k(\Delta t)^{1/\alpha}, \mu = 0)$ for nonzero $\beta$. 

The left-sided and right-sided Riemann-Liouville derivatives in $\mathbb{R}$ are defined, for $n = \lceil \alpha \rceil$, as 
\begin{align}
\phantom{}^{\phantom{}\text{RL}}_{\phantom{R}a} \mathbb{D}^{\alpha}_x u(x) &= 
\frac{1}{\Gamma(n-\alpha)}
\left[
\frac{\partial^n}{\partial z^n}
\int_{a}^z \frac{u(y)}{|z-y|^{\alpha-n+1}} dy
\right]_{z = x}
, 
\label{eq:leftRL}\\
\phantom{}^{\text{RL}}_{\phantom{R}x} \mathbb{D}^{\alpha}_{b} u(x) &=
\frac{(-1)^n}{\Gamma(n-\alpha)}
 \left[
\frac{\partial^n}{\partial z^n}
\int_z^{b}\frac{u(y)}{|z-y|^{\alpha-n+1}} dy
\right]_{z=x}
.
\label{eq:rightRL}
\end{align}
The texts of \citet{podlubny1998fractional}, \citet{oldham1974fractional}, and \citet{Meerschaert2012} discuss these operators in detail. 
These derivatives are frequently used in models with $a = -\infty$ and $b = \infty$. In connection with initial value problems, the left-sided Riemann-Liouville derivative in time, 
$\phantom{}^{\phantom{}\text{RL}}_{\phantom{R}0} \mathbb{D}^{\alpha}_t u(t)$,
is sometimes used with $a = 0$.
We have written the definitions \eqref{eq:leftRL} and \eqref{eq:rightRL} to avoid ambiguities in notation, and clearly show that substitution of the variable $x$ occurs after integration and differentiation.
An alternative approach is to define Riemann-Liouville fractional integrals separately, as in the right-hand sides of \eqref{eq:leftRL} and \eqref{eq:rightRL}; see \citet{Samko1993}.

\begin{figure}[t]
\includegraphics[width=\textwidth]{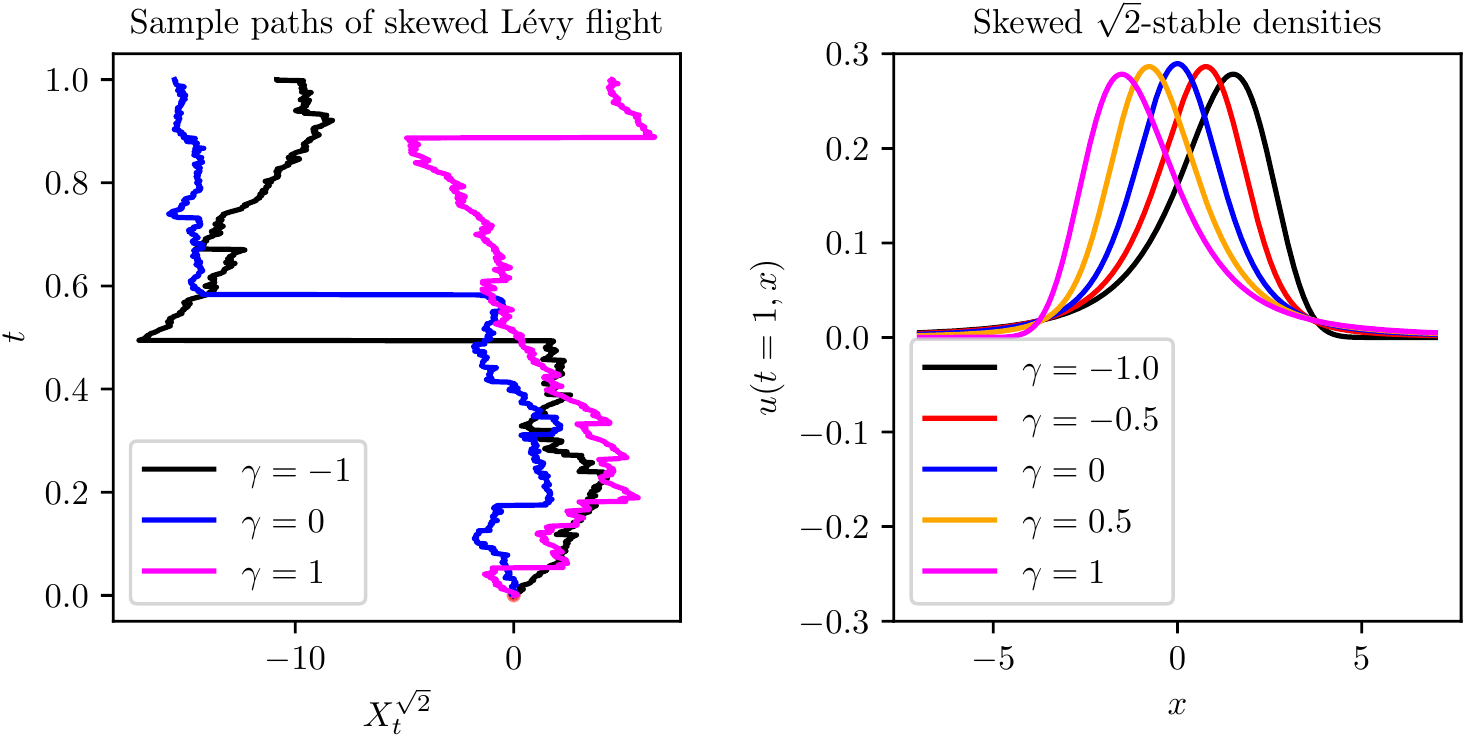}
\caption{$\alpha$-stable L\'evy flights allow for asymmetric diffusion, which has no analogue within the classical diffusion framework. The $\alpha$-stable density \eqref{eq:asymmetric_stable_distribution} admits a skewness parameter $\beta$, ranging from $-1$ to $1$, which can adjust the relative probability of long jumps in a given direction (\textit{right}), a statistical property that is evident in the sample paths (\textit{left}). Such models are governed by the fractional-order diffusion equation involving Riemann-Liouville derivatives, as in \eqref{eq:asymmetric_stable_diffusion}.}
\label{fig:asymmeric_stable_paths_comparison}
\end{figure}

One quirk of the notation for Riemann-Liouville derivatives in \eqref{eq:leftRL} and \eqref{eq:rightRL} is the writing of the upper and lower limits of integration $[a,x]$ and $[b,x]$, respectively, as subscripts. While this is suggestive, the result is that the variable of evaluation $x$ occurs twice in the notation for each operator. If these derivatives are evaluated at any numerical value of $x$, this value should be substituted in both locations; thus, $\phantom{}^{\phantom{}\text{RL}}_{{a}} \mathbb{D}^{\alpha}_5 u(5)$
represents a valid evaluation of the derivatives, but
$\phantom{}^{\phantom{}\text{RL}}_{{a}} \mathbb{D}^{\alpha}_x u(5)$ and
$\phantom{}^{\phantom{}\text{RL}}_{{a}} \mathbb{D}^{\alpha}_5 u(x)$ do not.

With $a = -\infty$ and $b = \infty$,
the Riemann-Liouville derivatives can be represented in frequency space by
\begin{align}\label{eq:RL_fourier_transforms}
\mathcal{F}\left[
{\phantom{}^{\phantom{,}\text{RL}}_{{-\infty}} \mathbb{D}^{\alpha}_x u} 
\right](\xi)
&= (-i\xi)^\alpha \mathcal{F}[u](\xi), \\
\mathcal{F}\left[
{\phantom{}^{\text{RL}}_{\phantom{R}x} \mathbb{D}^{\alpha}_\infty u}
\right](\xi)
&= (i\xi)^\alpha \mathcal{F}[u](\xi).
\end{align}
In one dimension, these can be used in the asymmetric diffusion model
\begin{equation}\label{eq:asymmetric_stable_diffusion}
\begin{aligned}
\frac{\partial u}{\partial t}(x,t) 
&= 
\frac{-k^\alpha}{\cos(\pi\alpha/2)}
\left[
p \,
\left( \phantom{}^{\phantom{,}\text{RL}}_{{-\infty}} \mathbb{D}^{\alpha}_x u(x,t) \right)
+
(1-p) \,
\left( \phantom{}^{\text{RL}}_{\phantom{R}x} \mathbb{D}^{\alpha}_{\infty} u(x,t) \right)
\right] \\
u(x, t = 0) &= u_0 (x),
\end{aligned}
\end{equation}
which describes anomalous diffusion of independent particles. Here, the positions of each particle at time steps of $k \Delta t$ for integer $k$ are governed by \eqref{eq:iid_rule} with increments $\Delta X$ being drawn from the \textit{asymmetric} $\alpha$-stable distribution
\begin{equation}\label{eq:asymmetric_stable_distribution}
X^{\alpha, p, \Delta t} \sim S_\alpha(\gamma = 2p-1,\sigma = k (\Delta t)^{1/\alpha},\mu = 0).
\end{equation}
Thus, the skewness ranges from $\gamma = -1$ when $p = 0$ to $\gamma = 1$ when $p=1$. The fundamental solution of \eqref{eq:asymmetric_stable_diffusion} is
\begin{equation}\label{eq:asymmetric_fundamental_solution}
u(x,t) = f_\alpha(x;\gamma = 2p-1,\sigma = k t^{1/\alpha},\mu); 
\end{equation}
cf. equation \eqref{eq:symmetric_fundamental_solution}.

Sample paths of the process just described are illustrated in Figure \ref{fig:asymmeric_stable_paths_comparison}.
Note that when $p = 1/2$, the distribution reverts to the symmetric $\alpha$-stable
distribution, and it can be shown in this case that equation \eqref{eq:asymmetric_stable_diffusion} reduces to
\eqref{eq:fractional_laplacian_equation}; more specifically, 
\begin{equation}\label{eq:riesz_operator}
\frac{1}{\cos(\pi\alpha/2)}
\left[
\frac{1}{2} \,
\left( \phantom{}^{\phantom{,}\text{RL}}_{{-\infty}} \mathbb{D}^{\alpha}_x u(x) \right)
+
\frac{1}{2} \,
\left( \phantom{}^{\text{RL}}_{\phantom{R}x} \mathbb{D}^{\alpha}_\infty u(x) \right)
\right]
=
(-\Delta)^{\alpha/2} u(x).
\end{equation}

The Fourier representation \eqref{eq:RL_fourier_transforms} suggests that the left-sided
Riemann-Liouville derivative  ${\phantom{}^{\phantom{,}\text{RL}}_{{-\infty}} \mathbb{D}^{\alpha}_x u}$
should be thought of as a fractional power of the operator $\partial / \partial x $. However, the correspondence between \eqref{eq:asymmetric_stable_diffusion} and \eqref{eq:asymmetric_stable_distribution} makes it clear that to obtain a complete description of 
$\alpha$-stable L\'evy flights in one dimension necessitates two operators, a left-sided and a right-sided
operator, which agree with one another when $\alpha = 2$. Our interest is these models lies in the fact that 
an extended centralized limit theorems hold for processes with i.i.d. increments drawn from distributions with infinite variance, but
for which the tails of the density function satisfy Pareto-type conditions as in \eqref{eq:pareto}. For such processes, $\alpha$-stable distributions
play an analogous role to the normal distribution in the classical central limit theorem; unlike the classical theorem, for full generality, skewed 
$\alpha$-stable distributions must be included in such a result. See \citet{meerschaert2001limit}
 or \citet{Meerschaert2012} for a treatment of these results. 
 
We mention how the Riemann-Liouville derivative can be utilized in dimensions $d > 1$. An anisotropic diffusion operator was introduced by \citet{Meerschaert1999} and \citet{Benson2000} as
\begin{equation}\label{eq:directional_laplacian}
-(-\Delta)_M^{\alpha/2}u(\MG{\xb})= C_{\alpha,d} \int_{|\MG{\bm{\theta}}| = 1} D^{\alpha}_{\MG{\bm{\theta}}} u(\MG{\xb}) M(d\MG{\bm{\theta}}), \quad 
C_{\alpha,d} = \frac{\Gamma(\frac{1-\alpha}{2}) \Gamma(\frac{d+\alpha}{2})}{2 \pi^{\frac{1+d}{2}}}.
\end{equation}
Here, $M(d\MG{\bm{\theta}})$ denotes a nonnegative measure on the angle $\MG{\bm{\theta}}$ in the unit sphere $\{ | \MG{\bm{\theta}} | = 1\}$ in $\mathbb{R}^d$, and
the Riemann-Liouville directional derivative is given by
\begin{equation}
D^{\alpha}_{\MG{\bm{\theta}}} u(\MG{\xb}) = 
\phantom{}^{\phantom{,}\text{RL}}_{{-\infty}} \mathbb{D}^{\alpha}_t v (t) \big|_{t = 0}, 
\quad \text{where } v(t) = u(\MG{\xb + t\bm{\theta}}). 
\end{equation}
\citet{Benson2000} showed that when the measure $M$ is uniform, the operator \eqref{eq:directional_laplacian} reduces to the fractional Laplacian \eqref{eq:fractional_laplacian_definition}. In higher dimensions and for general measures $M$, the operator \eqref{eq:directional_laplacian} plays an analogous role to the operator in the right-hand side of \eqref{eq:asymmetric_stable_diffusion}, which is in fact a special case of it for $d = 1$. As such, it is used in models of anistropic multivariate $\alpha$-stable L\'evy diffusion. 

\subsubsection{Subdiffusion and the Caputo fractional derivatives}
\label{sec:subdiffusion}
The superdiffusive model introduced above, in which a plume of particles spreads out in space with \MG{higher-order rate} $t^{1/\alpha}$ for $0 < \alpha < 2$, raises the question of whether a process can be constructed which results in diffusion \MG{{characterized by a \textit{lower-order} rate}} than the Brownian rate $t^{1/2}$. In this section, we introduce such a model, constructed as Brownian motion with random waiting times drawn from a skewed stable distribution, supported over positive real numbers with a power-law tail. Here, we step away from the framework of the SDE given by \eqref{eq:iid_rule}. Rather than being defined by a simple time-stepping scheme with i.i.d. increments, the paths of the process are defined by a transformation, or ``postprocessing'', of Brownian paths $B_t$. 

\begin{figure}[t]
\includegraphics[width=\textwidth]{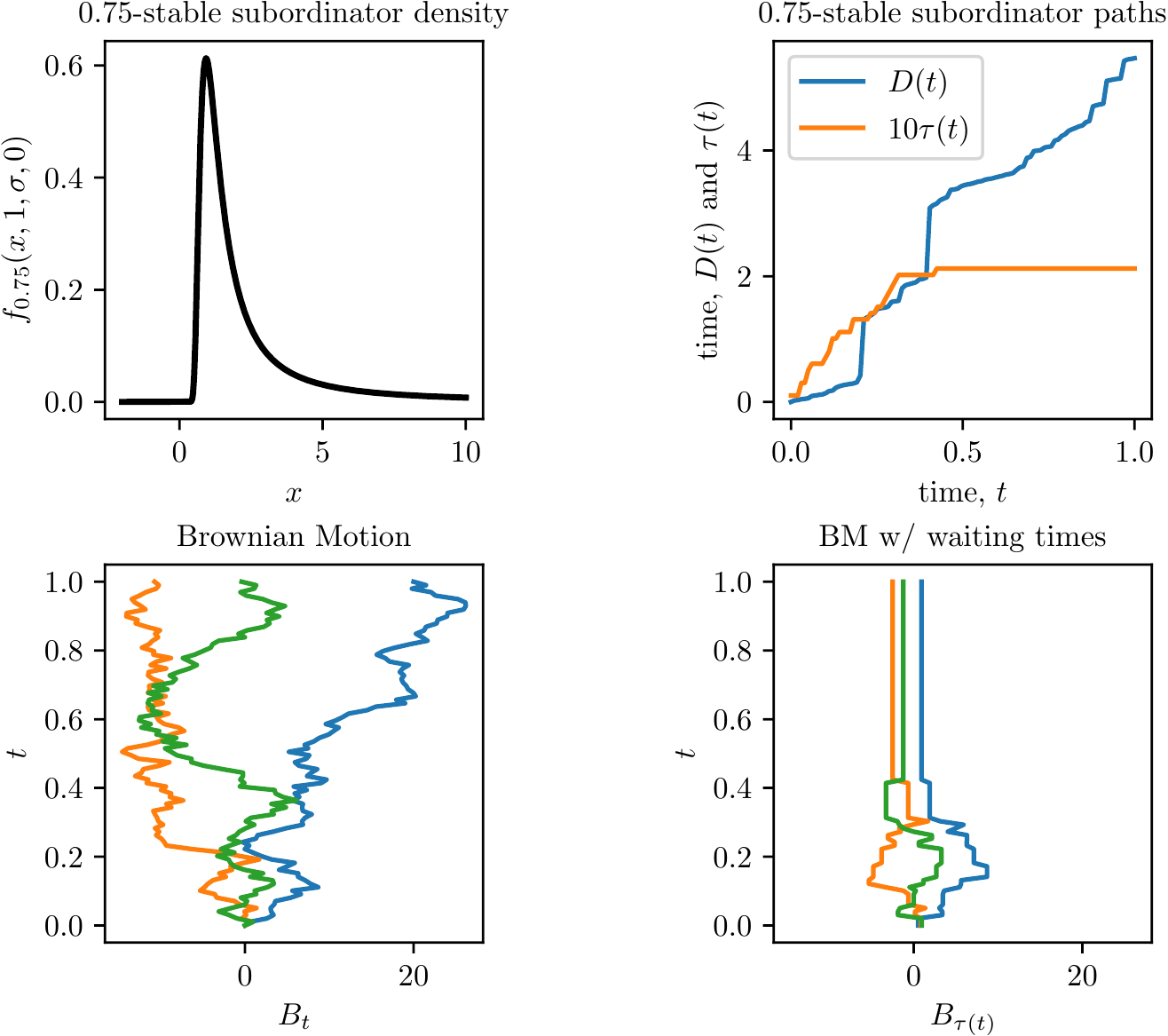}
\caption{
(\textit{Top left})
Example of an $\alpha$-stable subordinator density function, representing the density for random waiting times for the processes corresponding to the time-fractional diffusion equation \eqref{eq:time_fractional_diffusion}.
(\textit{Top right})
Sample path of the subordinator (cumulative waiting time) $D(t)$, the parent path, and the inverse subordinator (operational time) $\tau(t)$ given by \eqref{eq:subordinator_definition}. Note that as $t$ increases, $\tau(t)$ need not advance.
(\textit{Bottom left})
Three sample paths of Brownian motion.
(\textit{Bottom right})
Three sample paths of Brownian motion with waiting times, constructed from the Brownian paths in the bottom left panel. The particles trace out the same Brownian paths in space, but now wait for potentially several time steps at each location, as specified by the operation time $\tau(t)$. 
}
\label{fig:subordination}
\end{figure}

We introduce Brownian motion with waiting times, denoted by $B_{\tau(t)}$. The intuition is that the particle paths traced out in space by a discretization of $B_{\tau(t)}$ are paths of discretized Brownian motion $B_t$, but the particles wait at each point of the path for a random time drawn from the totally skewed stable distribution. The operational time $\tau(t)$, which introduces waiting and replaces linear time $t$, is an \textit{inverse stable subordinator}. This is a stochastic process in the variable $t$, although we write $\tau(t)$ rather than using a subscript for typographical reasons. This process is constructed by first defining the \textit{stable subordinator} $D(t)$, and defining $\tau(t)$ to be the inverse process\footnote{The definition of $\tau$ in terms of $D$ is an example of a right-continuous inverse of an increasing functions.
Paths of $D$, thought of as functions of $t$, are nondecreasing, so that each path of $\tau$ constructed in this way is a continuous-from-the-right inverse of the parent path of $D$ used to construct it.}  of $D(\tau)$. Both $D(t)$ and $\tau(t)$ are nondecreasing processes with units of time. In terms of paths, $\tau(t)$ arises from $D(t)$ as
\begin{equation}\label{eq:subordinator_definition}
\tau(t) = \inf \{ \tau \text{ such that } D(\tau) > t \}.
\end{equation}
Intuitively, $D(t)$ represents a cumulative waiting time process, keeping track of the total time waited by a particle throughout a path, while the inverse $\tau(t)$ represents an operational time, i.e., the time spent traveling. The increments of $D(t)$ represent the time waited at each location of a particle before the jump to the next location. 
More specifically, $D(t)$ is a totally skewed $\beta$-stable L\'evy process \eqref{eq:asymmetric_stable_distribution} with stability index $\beta \in (0,1)$, $\gamma = 1$, scale $\sigma = \cos(\pi \beta/2)$, and center $\mu = 0$; see \citet{Meerschaert2012}, Example 5.14. The construction of sample paths of $B_{\tau(t)}$ is demonstrated in Figure \ref{fig:subordination}.

 The resulting probability density function of $D(t)$,
 \begin{equation}\label{eq:subordinator_density}
 \psi_\beta(t) = 
 f_\beta(t;\gamma=1,\sigma= \cos(\pi \beta/2),\mu=0)
 \end{equation}
 for waiting times is supported in nonneagative real numbers.
 Due to the nonnegative support of the waiting time density, the characteristic function \eqref{eq:alpha_stable_characteristic} yields the Laplace transform of the waiting time density as
 \begin{equation}\label{eq:subordinator_laplace_transform}
\mathcal{L}\left[
\psi_\beta
\right](s)
=
e^{-s^\beta}, 
 \end{equation}
 where the Laplace transform is defined as
 \begin{equation}
\mathcal{L}\left[u\right](s)
=
\int_0^\infty
e^{-st} u(t) dt
.
\end{equation}
See \citet{Meerschaert2012}, p. 108 and p. 156 for a discussion. 
The variance of the process $B_{\tau(t)}$ is given by
\begin{equation}
\left\langle B_{\tau(t)} \right\rangle^2 = \frac{2}{\Gamma(\beta+1)} t^{\beta}, \quad 0 < \beta < 1,
\end{equation}
which is the desired subdiffusive property. Note that the finiteness of the variance does not imply that the normal central limit theorem applies to $B_{\tau(t)}$, which is not equal in distribution to Brownian motion nor to any L\'evy process. In fact, $B_{\tau(t)}$ is not a Markov processes.


The probability density of Brownian motion with waiting times $B_{\tau(t)}$ is governed by the time-fractional diffusion equation, 
\begin{align}\label{eq:time_fractional_diffusion}
\phantom{}^{\text{C}}_{{0}} \mathbb{D}^{\beta}_t u(t) 
&= k^2 \Delta u(x,t) \\
u(x, t = 0) &= u_0 (x)
.
\end{align}
Here, the Caputo derivative is defined for $0 < \beta < 1$ by
\begin{equation}\label{eq:caputo_definition}
\phantom{}^{\text{C}}_{a} \mathbb{D}^{\beta}_t u(t) = 
\frac{1}{\Gamma(1-\beta)}
\int_{a}^{t}
\frac{d u}{dt}(s) 
\frac{1}{|s-t|^{\beta}} ds
.
\end{equation} 
For $a = 0$, this operator is characterized by the simple Laplace transform representation (see \citet{Meerschaert2012}, page 111)
\begin{equation}\label{eq:laplace_transform_caputo}
\mathcal{L}\left[ 
\phantom{}^{\text{C}}_{{0}} \mathbb{D}^{\beta}_t u
\right](s) 
=
s^\beta \mathcal{L}[u](s)
-
s^{\beta - 1} u(0)
.
\end{equation}
Higher order Caputo derivative can be defined, although the Laplace transforms of the resulting operators involve initial conditions for derivatives of $u$; see Section 2.3 of \citet{Meerschaert2012}. The Caputo derivative is most frequently utilized as a derivative in time for initial-value problems, with the fractional order $0 < \alpha < 1$.

Before introducing the fundamental solution to the time-fractional diffusion, we introduce the Mittag-Leffler function \citep{mainardi2007sub,mainardi2020mittag}
\begin{equation}
E_\theta(z) = \sum_{\ell=0}^{\infty} \frac{z^\ell}{\Gamma(\theta \ell + 1)}, \quad \theta > 0. 
\end{equation}
This Mittag-Leffer $E_\theta(z)$ reduces to the exponential function $e^z$ when $\theta = 1$, and has Laplace transform property
\begin{equation}
\mathcal{L}\left[ E_\theta(-k^2 t^\theta) \right](s) = \frac{s^{\theta-1}}{s^\theta + k^2}, 
\end{equation}
which immediately implies that $E_\beta(-k^2 t^\beta)$ solves the fractional ordinary differential equation
\begin{equation}
\phantom{}^{\text{C}}_{0} \mathbb{D}^{\beta}_t u = k^2 u. 
\end{equation}
Returning to the diffusion equation \eqref{eq:time_fractional_diffusion} with initial condition $u(x,t=0) = \delta(x)$, applying the Fourier transform in space implies that
\begin{equation}
\mathcal{F}\left[ u(\cdot,t) \right](\xi) 
=
E_\beta(-k^2 \xi^2 t^\beta),
\end{equation}
which, as shown by \citet{mainardi2007sub}, yields a solution that can be written
\begin{equation}
u(x,t) = t^{-\beta/2} U(|x|/t^{\beta/2}); 
\quad
\end{equation}
with
\begin{equation}
U(x) = \frac{1}{2} 
\sum_{k=0}^\infty\frac{(-x)^k}{k!\Gamma[-(\beta/2) k + 1 - (\beta/2)]}.
\end{equation}
being a special case of the Fox-Wright function.
Note that $U(x) = u(x,t=1)$.
While the fundamental solution above is transcendental, it has the following properties: for $\alpha = 1$, it reduces to the solution
\eqref{eq:classical_diffusion_fundamental_solution} of the classical diffusion equation; for $0 < \alpha < 1$, the solution decays faster than exponential and slower than Gaussian; and the second moment of the solution 
is
\begin{equation}
\sigma^2(t) = 2 \frac{t^\beta}{\Gamma(\alpha+1)}
\end{equation}
Note that the $t^\beta$ scaling of this second moment is consistent with the scaling of the fundamental solution above. 

\subsubsection{Continuous time random walks and space-time fractional diffusion}\label{sec:CTRWs}
{
Both the $\alpha$-stable L\'evy flight $X^{\alpha,p}_t$, which led to the space-fractional diffusion equation discused in Section \ref{sec:RL}, and Brownian motion with $\beta$-stable subordinator operational time $B_{T^{\beta}(t)}$, which led to the time-fractional diffusion equation discussed in Section \ref{sec:subdiffusion}, are examples of \textit{continuous time random walks} \citep{Meerschaert2012}. A continuous time random walk (CTRW) allows for a general family of processes in space to be time-changed by a general family of waiting-time processes. To illustrate this concept, we consider the process $X^{\alpha,p}_{T^\beta(t)}$, which is $\alpha$-stable L\'evy flight $X^{\alpha,p}_t$ defined at the discrete level by \eqref{eq:asymmetric_stable_distribution} time-changed by the $\beta$-stable subordinator process $t \mapsto T^\beta(t)$ introduced in Section \ref{sec:subdiffusion}. This models a particle that performs independent jumps drawn from the $\alpha$-stable process, waiting at each point for a random time drawn independently from the $\beta$-stable subordinator process. As shown by, e.g., \citet{Meerschaert2012} (Section 4.5), the probability density of this particle position is then governed by a differential equation that is fractional in both time and space,
\begin{equation}\label{eq:spacetime_asymmetric_stable_diffusion}
\begin{aligned}
\phantom{}^{\text{C}}_{0} \mathbb{D}^{\beta}_t u(t,x)
&= 
\frac{-k^\alpha}{\cos(\pi\alpha/2)}
\left[
p \,
\left( \phantom{}^{\phantom{,}\text{RL}}_{{-\infty}} \mathbb{D}^{\alpha}_x u(x,t) \right)
+
(1-p) \,
\left( \phantom{}^{\text{RL}}_{\phantom{R}x} \mathbb{D}^{\alpha}_{\infty} u(x,t) \right)
\right] \\
u(x,t = 0) &= u_0 (x),
\end{aligned}
\end{equation}
}

{While intuitive, this result deserves a more detailed outline within the general theory of CTRWs.  In the standard CTRW model, particles wait at a location for time drawn from a density function $\psi$, and jump to a new location by an increment drawn from a density function $\phi$. The waiting time and jump samples are assumed to be i.i.d., and uncoupled from each other \citep{zaburdaev2015levy, MeKl00, scalas2004uncoupled,torrejon2018generalized}.  Thus, the densities $\psi$ and $\phi$ completely determine the CTRW. From the waiting time density $\psi$, the probability that a particle will remains at any given position for time $t$ is
\begin{equation}\label{eq:survival_probability}
\Psi(t) = 1 - \int_0^t \psi(t) dt;
\end{equation}
this is referred to as the \textit{survival probability} of a CTRW particle. Then, given an initial probability density of a particle $u_0(x) = u(x,t=0)$, which can also be thought of as an initial distribution of an ensemble of independent particles, the following equation was derived by \citet{Montroll1965} for the density at later times:
\begin{equation}\label{eq:ctrw_master}
u(x,t) = \Psi(t) u_0(x)  - \int_0^t \psi(t-\tau) 
\int_{-\infty}^{\infty} \phi(y) u(x-y,\tau) dy d\tau .
\end{equation}
This equation is central to the CTRW theory\footnote{This equation was also derived by \citet{scher1973stochastic1, scher1973stochastic2} and is referred to as the CTRW equation of Scher and Lax by \citet{klafter1980derivation}. Other authors, such as \citet{torrejon2018generalized} refer to this as the \textit{master equation} of a CTRW.}. Taking the Laplace transform in time, the Fourier transform in space, and solving for $\mathcal{F} \left[\mathcal{L} [u]\right] (\xi,s)$ yields the Montroll-Weiss equation \citep{Montroll1965},
\begin{equation}\label{eq:montroll_weiss}
\mathcal{F}\left[  \mathcal{L} [u] \right](\xi,s)
=
\frac{1-\mathcal{L}[\psi](s)}{s}
\frac{\mathcal{F}[u_0](\xi)}{1-\mathcal{L}[\psi](s) \mathcal{F}[\phi](\xi)}. 
\end{equation}
In the case that $\phi$ is the $\alpha$-stable density \eqref{eq:asymmetric_stable_distribution} and $\psi$ is the $\beta$-stable subordinator density \eqref{eq:subordinator_density}, then $\mathcal{F}[\phi]$ is given by the analytical formula \eqref{eq:alpha_stable_characteristic} and $\mathcal{L}[\psi]$ by \eqref{eq:subordinator_laplace_transform}, so that the Montroll-Weiss equation represents a closed-form solution of $u$ in $(\xi,s)$-space. Unsurprisingly, it is impossible to perform inverse transforms and obtain $u$ itself analytically, but $u$ can be shown to satisfy \eqref{eq:spacetime_asymmetric_stable_diffusion} using the representations \eqref{eq:laplace_transform_caputo} and \eqref{eq:RL_fourier_transforms} \citep{Meerschaert2012}. 
}

\subsubsection{L\'evy walks and fractional material derivatives}\label{sec:levy_walks}
{Superdiffusive $\alpha$-stably L\'evy flight exhibits infinite MSD, which is a drawback for certain applications. Related to this is the infinite speed of propagation intrinsic to L\'evy flights, i.e., the fact that particles have a nonzero probability of traveling an arbitrary large distance in a unit of time. Brownian motion also suffers from this feature, although this probability of large excursions is so low that MSD remains finite. A prototypical model of superdiffusion that cannot be described by a L\'evy flight is ballistic motion, in which particles simply move from an initial configuration in fixed random directions with speed $v$, for all time $t$. A ballistic particle travels a distance $vt$ in time $t$ from an initial position $x_0$. If reorientations are allowed, then the positions of these so-called sub-ballistic particles in space-time are confined to a ballistic cone 
\begin{equation}
\left\{ (x,t) \text{ such that } x \in [x_0 - vt, x_0+vt], t \ge 0 \right\}.
\end{equation}
Because the density function of the particle positions is compactly supported, all moments of the position are finite. Such a process cannot be described by L\'evy flights.
}

\begin{figure}[t]
\centering
\includegraphics[width=0.75\textwidth]{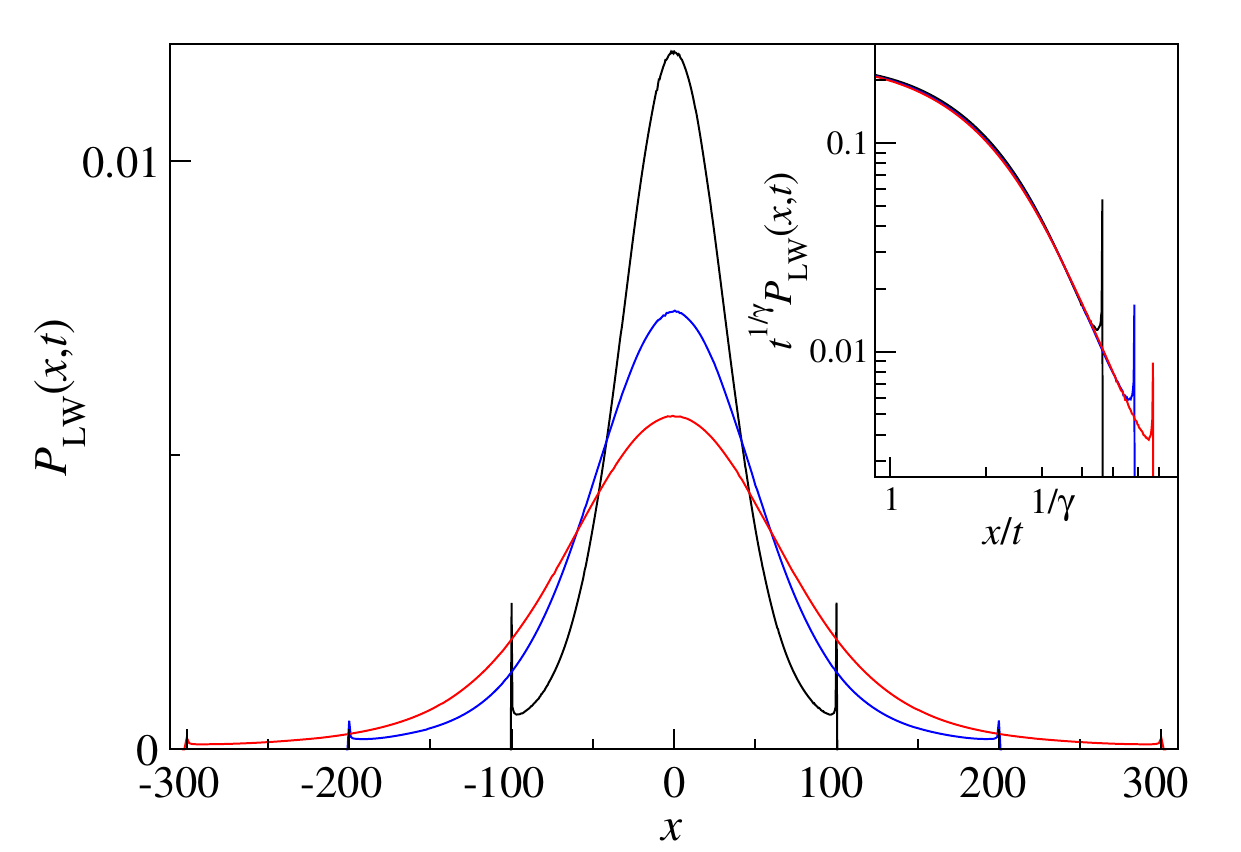}
\caption{The evolution of the probability density function (denoted $P_{\text{LW}}$ in the figure) of a L\'evy walk, reproduced from \cite{zaburdaev2015levy}. Here, $\gamma = 3/2$ and the density is plotted for $t = 100$ (black), $t = 200$ (blue), and $t = 300$ (red). The density mimics the density of a $\gamma$-stable L\'evy flight in an interior region of the ballistic cone, scaling outwards as $t^{1/\gamma}$, supported inside the ballistic front (consisting of two points in one dimension) that scales outwards as $t$.}
\end{figure}

{
To capture such behavior, we introduce the L\'evy \textit{walk} model, following \citet{zaburdaev2015levy}. Such models are based on continuous-in-time motion of particles, rather than instantaneous jumps. A speed $v$ of particles in a medium is specified; each particle moves with speed $v$ in a chosen direction, before a reorientation event occurs in which the direction changes instantaneously and the particle continues to move with speed $v$ before the next direction. Assuming the direction at reorientation is sampled uniformly on the unit sphere, such a walk is determined by a probability density function for the duration of movement $\psi(\tau)$. This leads to a survival probability $\Psi(t)$ given by \eqref{eq:survival_probability}, with $\psi$ now representing the duration density. Thus, $\Psi(t)$ returns the probability that a particle has persisted in a given direction for time $\tau$, i.e., has \textit{not} experienced reorientation for time $\tau$. Similar to the CTRW case, a master equation can be derived for the probability density $u(x,t)$ of the location of the particle in Laplace-Fourier space:
\begin{equation}\label{eq:levy_walk_master_eq}
\mathcal{F}\left[  \mathcal{L} [u] \right](\xi,s) 
=
\frac{\mathcal{L}[\Psi](s + i v \xi) + \mathcal{L}[\Psi](s - i v \xi ) }{2 - \mathcal{L}[\psi](s + i v \xi) + \mathcal{L}[\psi](s - i v \xi)} \mathcal{F}[u_0](\xi). 
\end{equation}
Unlike the master equation for CTRWs, this equation exhibits coupling in Fourier and Laplace variables, representing coupling in space-time\footnote{A L\'evy walk may be compared to a non-standard CTRW in which waiting times prior to jumps are correlated to the jump length, e.g., propertional to the jump length, so that long excursions are penalized by long waiting times. See \citet{}}. 
This results in governing equations that are considerably more complex than those of a standard CTRW. For a L\'evy walk, $\psi$ is taken to be a Pareto-type distribution, 
\begin{equation}
\psi(\tau) = \frac{1}{\tau_0} \frac{\gamma}{(1 + \tau/\tau_0)^\gamma}, \quad \tau_0 > 0, \gamma > 0.
\end{equation}
An asymptotic expansion of $\mathcal{L}[\psi]$ and $\mathcal{L}[\Psi]$ substituted in \eqref{eq:levy_walk_master_eq} yields the following approximation for the evolution of the density function of a L\'evy walk in Fourier-Laplace space:
\begin{equation}\label{eq:levy_walk_master_eq}
\mathcal{F}\left[  \mathcal{L} [u] \right](\xi,s) 
\approx
\frac{(s + i v \xi)^{\gamma-1} + (s - i v \xi )^{\gamma-1} }{(s + i v \xi)^\gamma + (s - i v \xi)^\gamma} \mathcal{F}[u_0](\xi). 
\end{equation}
Given $v$ and $u_0$, this equation can be inverted to compute $u(x,t)$, but obtaining a governing equation in $(x,t)$ is less straightforward from this point on, due to space-time coupling. \citet{sokolov2003towards} suggest defining a fractional \textit{material} or \textit{substantial} derivative
\begin{equation}
(v^{-1} \partial_t \pm \partial_x)^{1/\gamma} u := \mathcal{F}^{-1} \mathcal{L}^{-1}
\left[
(s + i v \xi)^\gamma + (s - i v \xi)^\gamma \mathcal{F}[u_0](\xi)
\right],
\end{equation}
in order to obtain a governing equation for $u(x,t)$. Recent works, such as those of \cite{chen2015discretized}, have explored numerical discretizations for these operators.
}

{
Despite the greater mathematical difficulties related to governing equations, as compared to other fractional models,
L\'evy walks have been widely used due to the physical nature of finite speed of propagation and finite MSD; see \cite{zaburdaev2015levy} for a survey. When $1 < \gamma < 2$, by numerical approximations, it can be seen that $u(x,t)$ evolves from a $\delta$-distribution with ``a central part of the profile approximated by the Lévy distribution sandwiched between two ballistic peaks'' that propogate at speed $v$, with an MSD and self-similarity property for large $t$ that features a superdiffusive scale factor of $t^{1/\gamma}$ \cite{zaburdaev2015levy}. 
}

\subsubsection{Variable-order fractional derivatives} \label{sec:variable-order}
Given the physical meaning within stochastic models of the fractional order $\alpha$ in derivatives such as \eqref{eq:fractional_laplacian_definition}, \eqref{eq:leftRL}, and \eqref{eq:caputo_definition}, it is reasonable to expect that these parameters may vary in space and time. Variable-order fractional models are convenient to describe anomalous diffusion in the case of heterogeneous materials or media, or, more generally, when the nature of the diffusion process (subdiffusive, superdiffusive, and classical) changes with space and time. While models with constant fractional order are the simplest and most widely used, some of the model descriptions we discuss in the following sections are improved by the use of a variable fractional order. 
In recent years, with the purpose of increasing the descriptive power of fractional operators, new models characterized by a variable fractional order have been introduced for both space- and time-fractional differential operators \cite{AntilRautenberg2019_SobolevSpacesWithNon,Darve2021fractional,DElia2021fractional,Razminia2012,Zheng2020} and several discretization methods have been designed \cite{Chen2015,SchneiderReichmannEtAl2010_WaveletSolutionVariableOrderPseudodifferentialEquations,Zeng2015,Zheng2020optimal-order,Zhuang2009}. 
The improved descriptive power of variable-order fractional operators has been demonstrated in some recent works on parameter estimation \cite{Pang2020,Pang2019fPINNs,Wang2020Inverting}.

Given a function 
\begin{equation}
\alpha: \mathbb{R}^d \times \mathbb{R} \rightarrow \mathbb{R}, 
\end{equation}
i.e., a function $\alpha(\xb,t)$ of space and time, we define variable-order operators as follows. For a function $u(\xb,t)$ with $\xb \in \mathbb{R}^d$ and $t \in \mathbb{R}$, we define the variable-order fractional Laplacian\footnote{{For more recent works and novel definitions of variable-order fractional Laplacians we refer the reader to \cite{AntilRautenberg2019,DElia2021fractional}.}} as 
\begin{equation}
\mathfrak{L}^{\alpha(\cdot,\cdot)} u(\xb,t)
=
C_{d,\alpha(\xb,t)} \text{ p.v.} 
\int_{\mathbb{R}^d} \frac{u(\xb,t) - u(\yb,t)}{|\xb-\yb|^{d+\alpha(\xb,t)}} d\yb.
\end{equation}
Here, $\alpha(\xb,t)$ is restricted to take values in $(0,2)$. Note that for constant $\alpha$,  $\mathfrak{L}^{\alpha} = (-\Delta)^{\alpha/2}$.
For $d = 1$ and $\alpha(x,t)$ restricted to $(0,1)$, we define the variable-order left-sided Riemann-Liouville fractional derivative as
\begin{equation}
\phantom{}^{\phantom{,}\text{RL}}_{{-\infty}} \mathbb{D}^{\alpha(\cdot,\cdot)}_x u(x,t) 
= 
\frac{1}{\Gamma(1-\alpha(x,t))}
\frac{\partial}{\partial x}
\int_{-\infty}^x \frac{u(y,t)}{|x-y|^{\alpha(x,t)}} dy,
\end{equation}
The right-sided Riemann-Liouville may be defined for variable order in an analogous way. 
We define the variable-order Caputo fractional derivative, again for $\alpha(x,t)$ taking values in $(0,1)$, as
\begin{equation}\label{eq:variable_order_caputo}
\phantom{}^{\text{C}}_{0} \mathbb{D}^{\alpha(\cdot,\cdot)}_t u(x,t) = 
\frac{1}{\Gamma(1-\alpha(x,t))}
\int_{-\infty}^t
\frac{d u}{dt}({s}) 
\frac{1}{|s-t|^{\alpha(x,t)}} {ds}.
\end{equation}

\subsubsection{Relationships between processes, fractional models, and applications}\label{sec:proc-oper-appl}
To summarize and offer a quick look-up of anomalous diffusion processes, their corresponding fractional models, and applications of each process/model, we have included these relationships in Table \ref{table1}. This table includes references to the previous sections where each process and model is described, as well as pointers to the applications in the following sections where the models are utilized. We have limited references to applications to only those three areas that we focus on in this article. 
\begin{table}[t]
\begin{tabular}{
  |p{.225\textwidth}|
   p{.25\textwidth}|
   p{.075\textwidth}|
   p{.3\textwidth}|}
\hline
\textbf{Process}                                   & \textbf{Fractional Model}                            & \textbf{Section}         & \textbf{Application}                                                                                \\ \hline
Symmetric $\alpha$-stable L\'evy Flight            & Space-fractional diffusion with fractional Laplacian & \ref{sec:frac_lapl}      & Scalar turbulence \textsection \ref{sec: nonlocality}.                                                                                                                                                   \\ \hline
Asymmetric $\alpha$-stable L\'evy Flight           & Space-fractional diffusion with R.L. derivatives     & \ref{sec:RL}             & Subsurface flows through fractured media \textsection \ref{sec:subsurface_space_frac}.                                                                                                                  \\ \hline
Brownian Motion with stable L\'evy waiting times & Time-fractional diffusion with Caputo derivative     & \ref{sec:subdiffusion}   & Rheology and failure of viscoelastic materials \textsection \ref{Sec:VE}, \ref{sec:plasticity}, \ref{sec:damage}, mobile-immobile subsurface flow dynamics \textsection \ref{sec:time_fractional_subsurface}. \\ \hline
CTRW with stable L\'evy jumps and waiting times                                       & Space-and-time-fractional diffusion                  & \ref{sec:CTRWs}          & Subsurface flows through fractured media with trapping zones \textsection \ref{sec:subsurface}.                                                                                                          \\ \hline
L\'evy Walk                                        & Material derivative with coupled Fourier-Laplace space  & \ref{sec:levy_walks} & Superdiffusion as for L\'evy flights, but with finite MSD and confined to ballistic cone; see \citet{zaburdaev2015levy}.                                                                                \\ \hline
\end{tabular}
\caption{Relationships between diffusion models, fractional models, and applications discussed in this article.}
\label{table1}
\end{table}

\subsection{Connection to Nonlocal Calculus}\label{sec:nonlocal_connection}

Fractional-order differential operators can be viewed as a special case of nonlocal models \citep{Defterli2015,d2020unified,d2013fractional}. The intrinsic nonlocality of fractional operators has been illustrated in the previous section; this property describes the fact that fractional-order derivatives of a function at a point $\xb \in \mathbb{R}^d$ typically depend on values of the same function at all points $\yb \in \mathbb{R}^d$, no matter how large the distance between $\xb$ and $\yb$ may be.  An example of this is the formula \eqref{eq:fractional_laplacian_definition} for the fractional Laplacian. 

General nonlocal diffusion (or Laplace) models include integral operators of the form \cite{Du2012,Du2013} 
\begin{equation}\label{eq:nonlocal-laplacian}
\mathfrak{L}[u](\xb) =
\int_{\mathbb{R}^d}
\gamma(\xb,\yb)[u(\xb) - u(\yb)] d\yb
  \end{equation}
with kernels $\gamma$ having support in $\{|\xb-\yb| \le \delta\}$, where the so-called {\it interaction radius} $\delta$ is such that $\delta\in(0,\infty]$. A quick comparison with the integral formula \eqref{eq:fractional_laplacian_definition} shows that when the kernel $\gamma$ is properly selected and $\delta=\infty$, then the fractional Laplacian is formally equivalent to \eqref{eq:nonlocal-laplacian} (see \cite{d2013fractional} for a rigorous derivation and a discussion). 

Nonlocal Laplace operators featuring kernels with bounded support may be preferred to fractional operators for physical reasons when modeling short-range interactions \cite{Askari2008,Silling2000} as well as mathematical convenience when posing volume conditions, the nonlocal counterpart of classical boundary conditions \cite{DEliaNeumann2019,Du2012}. The latter reason gives rise to truncated fractional-order derivatives \cite{Burkovska2021,d2020unified}.

General nonlocal models also allow for more flexibility with regards to regularity. Considering diffusion or Poisson's problems, fractional-order problems exhibit regularity explicitly parametrized by the fractional order \cite{Samko1993}; in contrast, nonlocal models involving nonsingular kernel operators lead to problems that impose no regularity on the solution \cite{Du2012} and can be naturally utilized to model fracture dynamics \cite{Silling2000,Silling2007}. Finally, we remark that the relationship between fractional and nonlocal models extends to more general operators than those of diffusion/Laplace type. There is indeed a well-established nonlocal vector calculus \citep{Du2013,Gunzburger2010}, of which fractional-order vector calculus is a special case (see \citep{d2020unified} for rigorous results where the convergence of truncated fractional gradient and divergence is proven in norm and pointwise). 

\subsection{A Remark about Numerical Methods for Fractional-Order Models}\label{sec:numerical_methods}

Over the past two decades, a significant amount of progress has been made in developing numerical methods, ranging from finite-difference/volume schemes to finite-element methods, in addition to a variety of new spectral theories for single and multi-domain spectral methods, obtaining efficient and easy-to-construct smooth/non-smooth basis and test functions. Performing a thorough and inclusive review of all the contributions made in this direction is nearly impossible and out of the scope of present work. Interested readers can find a wide spectrum of research carried out in the context of numerical analysis of fractional models in \cite{almeida2015computational,baleanu2012fractional,chen2016generalized,DElia2020Acta,DElia2020cookbook,diethelm2005algorithms,furati2018advances,gorenflo1997fractional,guo2015fractional,li2019numerical,samiee2021unified,tarasov2019handbook,zayernouri2013fractional_SLPs,zayernouri2015tempered,Zayernouri2022-CUPbook,zhou2020implicit}, and references therein. 

We restrict ourselves to discussing one aspect related to numerical methods, on the computational feasibility of solving fractional models. 
In the time-fractional case, efficient long-time numerical integration is of interest to capture inherent long time far-from-equilibrium dynamics and to enable the full convolution computations for large-scale systems. To this end, a number of fast time-stepping schemes have been developed during the last 20 years, which greatly reduce the cost of solving fractional models, making them quite comparable to clasical models. These include the \textit{fast convolution method} by \citet{Lubich2002}, which reduced the computational complexity of direct finite-difference discretizations of time-fractional models from $\mathcal{O}(N^2)$ to $\mathcal{O}(N \log N)$, and memory requirements from $\mathcal{O}(N)$ to $\mathcal{O}(\log N)$, where $N$ denotes the number of time steps. High-order extensions of the method were developed \cite{Yu2016JCP,Zeng2018stable} and applied to three-dimensional simulations of fluid-structure interactions in cerebral arteries and aneurysms \cite{Yu2016JCP}. Among a vast number of works in the literature, we also briefly outline matrix-based schemes, such as fast-inversion approaches \cite{Lux2015} and kernel compression methods \cite{Baffet2017} for time-fractional problems. For space-fractional FPDEs, adaptive methods and hierarchical matrices approaches have accomplished similar, dramatic reductions in computational complexity and memory costs for solving models \cite{zhao2017adaptive, xu2018efficient, li2020fast}. Efficient solvers and preconditioners for the fractional Laplacian were also developed by \citet{ainsworth2017aspects}.
The point we make is that from two decades of numerical methods development in the field, the current state-of-the-art numerical methods for fractional models produce computational costs comparable to integer-order cases, therefore being timely computational tools to be readily employed in large-scale systems modeled by FPDEs.

\section{Anomalous subsurface transport} \label{sec:subsurface}

The accurate prediction at large scales of contaminant transport in both surface and subsurface water is fundamental for efficient management of water resources and hence critical for environmental safety. However, the explicit description of the systems where transport takes place is extremely challenging, especially at large scales, due to the complexity of surface and subsurface environments. In fact, the latter \MG{f}eature heterogeneities that are either hard or impossible to measure and, hence, cannot be described with certainty at all scales and locations of relevance. On the other hand, even when the environment's microstructure can be captured, numerical simulations of PDE models such as the advection-diffusion equation (ADE) may be prohibitively expensive if conducted at small scales. Furthermore, \MG{these} same equations that are accurate at small scales \MG{f}ail to predict solutes' behavior at larger scales, due to the appearance of ``anomalous'', or ``non-Fickian'' behavior \cite{Levy2003}.

Still, in the past, the classical ADE has been broadly utilized as a model for solute transport \cite{Decker1999,Erel1998,Matthess1991}. As thoroughly explained in \cite{Neuman2009}, in the presence of heterogeneous media, ADEs fail to be accurate at large scales.
\MG{Such classical models at the} coarse-grained \MG{scale} can be considered \MG{accurate} only when media properties do not vary rapidly in the neighborhood of a point. \MG{Even with mild heterogeneities, quantities defined at large scales vary rapidly enough to be treated as random functions of space and/or time, in which case the ADE becomes an SDE \cite{Benson2001}.} Interestingly, when treating the ADE's parameters as stochastic, the ensemble mean concentration through randomly heterogeneous media is generally non-Fickian, i.e. non-classical. This can be observed in a simple manner by performing Monte Carlo numerical simulations. After generating several random realizations of the underlying velocity field, the ADE is numerically solved for each field and the concentration is averaged over all realizations, revealing \MG{such} non-classical behavior \cite{Neuman2009}.

In view of the following section where fractional behavior is discussed in the context of turbulence, we point out that the above stochastic theories are closely related to those governing turbulent diffusion. However, while transport in porous media takes place at small Reynolds numbers, the latter take place at large ones. Furthermore, porous velocities depend on hydraulic properties in a known manner, whereas turbulent velocities fluctuate randomly in space–time, making the first uncertainty epistemic (e.g. incomplete knowledge of
medium properties) and the second aleatory (i.e. controlled by chance).
This makes it easier to reduce the uncertainty in solute transport models by tuning them using hydrogeologic data (see e.g. \cite{Tartakovsky2007}).

In this section we show that Fractional ADEs (FADEs) are appropriate models to describe non-Fickian transport of solutes without the prohibitive burden of resolving the heterogeneities at the small scales explicitly thanks to their integral nature that allows to embed length-scales in the definition of the operator. Before reporting on early works featuring a \MG{simple} fractional Laplacian model and later works where variable fractional orders are introduced, we dedicate a few words to another nonlocal model, also popular in the literature: the continuous time random walk (CTRW) approach. As we point out later on, these models have similarities and share advantages, \MG{the most important of which is perhaps} the strong connection to stochastic processes that makes them easier to analyze and interpret.

\paragraph{Fractional subsurface models based on continuous time random walks.} 
{In Section \ref{sec:CTRWs}, we discussed the basic concepts of CTRWs, introduced by \citet{Montroll1965}.
We now explain how these models arise in subsurface transport and lead to fractional equations, following \citet{Berkowitz2006}; further relevant works in the literature include \MG{\citet{Berkowitz1995,Berkowitz1998,Boano2007,Dentz2003}, and \citet{Valocchi1989}.}}
 
To analyze subsurface particles, we begin by examining the solute concentration $C(x,t)$ for a given configuration of particles; $C(x,t)$ refers to the number of particles at a site $x$, normalized by the total number of particles in the system. 
In the absence of sinks and sources, the solute concentration $C(x,t)$ varies with time $t$ at the site $x$ by following a stochastic mass balance expression, i.e.
\begin{equation}\label{eq:master-equation}
\dfrac{\partial C(x,t)}{\partial t}=
-\sum\limits_y \big(w(y,x)C(x,t) 
- w(x,y)C(y,t)\big).
\end{equation}
The expression above is known in the literature as (discrete) {\it master equation}  \cite{Oppenheim1977}. Here, $w(x,y)$ is the transition rate at which a particle moves from $y$ to $x$, the first term in the sum represents the normalized rate of solute outflow from site $x$ to all sites $y$, whereas the second term represents the normalized rate of solute inflow from all sites $y$ to $x$.
We further assume that the transition
rates corresponding to different sites or displacements are statistically independent, i.e. hydraulic and transport properties of porous media and system states (e.g. hydraulic fluxes) lack spatial correlations. This is referred to as statistical incoherence; under this assumption, the ensemble mean concentration $c(x,t)=\langle C(x,t) \rangle$, where $\langle \cdot \rangle$ refers to an average over all possible configurations of the particle system,
satisfies the so called {\it generalized master equation}, i.e.
\begin{equation}\label{eq:g-master-equation}
\dfrac{\partial c(x,t)}{\partial t}=
-\sum\limits_y \int\limits_0^t 
\Big(\theta(y-x,t-\tau)c(x,\tau)
- \theta(x-y,t-\tau)c(y,\tau)\Big)\,d\tau.
\end{equation}
%
%
As discussed in \citet{Berkowitz2006}, this equation is equivalent to a spacetime coupled CTRW equation
\begin{equation}
c(\MD{x},t) 
=
\sum_y
\int_0^t
\chi(s'-s,t'-t)c(s',t') dt' + \delta(s) \delta(t - 0^+),
\end{equation}
with an explicit correspondence between the function $\theta$ and the space-time density function $\chi(s,t)$; see also \citet{klafter1980derivation}. If the CTRW is uncoupled, i.e., 
\begin{equation}
\chi(s,t) = \psi(s) \phi(t),
\end{equation}
then this equation is equivalent to the CTRW equation \eqref{eq:ctrw_master} discussed in Section \ref{sec:CTRWs}.
As a result, $c(x,t)$, in absence of advection, and with $\psi$ and $\phi$ given by the stable distributions specified in Section \ref{sec:CTRWs}, is governed by the FPDE \eqref{eq:spacetime_asymmetric_stable_diffusion}. 
This cements the importance of FPDEs in subsurface transport, although in some cases, the incoherence assumption that is required to derive the generalized master equation may not be valid. 

A simple and fairly general FADE for subsurface transport under the influence of both advection and anomalous diffusion is the one-dimensional advection and space-time-fractional diffusion equation with constant coefficients (see, e.g., \cite{Sun2020}):
\begin{equation}
\label{eq:simple-1d-fade}
\begin{aligned}
\phantom{}^{\text{C}}_{{0}} \mathbb{D}^{\beta}_t c(x,t)
&= -V\dfrac{\partial c}{\partial x} 
-D \left[p \,
\left( \phantom{}^{\text{RL}}_{\phantom{R}{-\infty}} \mathbb{D}^{\alpha}_x c(x,t) \right)
+
(1-p) \,
\left( \phantom{}^{\text{RL}}_{\phantom{R}x} \mathbb{D}^{\alpha}_{\infty} c(x,t) \right)
\right] \\
c(x,t = 0) &= c_0 (x),
\end{aligned}
\end{equation}
where $c$ is the solute concentration, $V$ a constant velocity, $D$ a constant diffusion coefficient and $\alpha$ the fractional order. 
In Section \ref{sec:CTRWs}, we presented equation \eqref{eq:spacetime_asymmetric_stable_diffusion}, which is identical to the above equation except for the advection term $-V \partial c / \partial x$, as the governing equation for the probability density of a continuous-time random walk. As discussed by \citet{Meerschaert2012}, the inclusion of the advection term corresponds to a stochastic model in which the particle drifts with constant velocity and jumps to the left or the right with density specified by the diffusion term. 
Thus, when $\beta = 1$ in \eqref{sec:CTRWs}, the FADE \eqref{eq:simple-1d-fade} governs the evolution of the probability density
\begin{equation}
f_\alpha(x;\gamma = 2p-1,\sigma = k (\Delta t)^{1/\alpha},\mu = Vt)
\end{equation}
of the skewed $\alpha$-stable process $S_\alpha(\gamma = 2p-1,\sigma = k (\Delta t)^{1/\alpha},\mu = Vt)$.
Comparing to the asymmetric diffusion model in Section \ref{sec:RL}, with fundamental solution \eqref{eq:asymmetric_fundamental_solution},
this density differs only in that the center drifts with velocity $\MG{V}t$. 
This describes a particle that drifts with velocity $Vt$ and makes jumps to the left or right drawn from the stable distribution \eqref{eq:asymmetric_stable_distribution}. 
More specifically, when particle path ware discretized in steps of $\Delta t$, the position of the particle increments by $V \Delta t + X^{\alpha, p, \Delta t}$, where $\Delta X^{\alpha, p, \Delta t}$ is given by \eqref{eq:asymmetric_stable_distribution}, at each time step. 
When $0 < \beta < 1$, similar to the CTRW model described in Section \ref{sec:CTRWs}, this equation governs the probability density of a particle undergoing the process just described, time-changed by the inverse $\beta$-stable subordinator, again introducing waiting times to the process. 

As pointed out by \citet{Neuman2009}, when $\beta = 1$, \eqref{eq:simple-1d-fade} corresponds to a Markovian random walk processes of statistically independent and identically distributed non-Gaussian displacements, and, as such, they can only occur in an uncorrelated velocity field; in hydrology, this can be viewed as a limitation of both CTRW and plain FADEs. 
Instead, it is possible that variable-coefficient or variable-order models may be able to describe processes associated with statistically non-homogeneous velocity fields. 
However, we are not aware of a specific theoretical framework that relates variable-coefficient and variable-order FADEs and CTRWs. 
Nor are we aware of a framework that relates such variable parameters to physical properties of the medium.
At present, the only way of estimating such parameters is by fitting the models to observed concentrations and/or mass fluxes, and not by hydraulic data such as hydraulic conductivity, advective porosity or flow parameters such as hydraulic gradients, fluxes and advective porosities \cite{Neuman2009}.

As discussed in Section \ref{sec:CTRWs}, limits of a CTRW with infinite and statistically independent waiting times leads to time-fractional FPDEs. A physical mechanism that would result in time-fractional derivatives in a FADE is particle trapping due to media heterogeneities \cite{Compte1996,Giona1992}. Such models are discussed in Section \ref{sec:time_fractional_subsurface}.

We conclude this section with advantages in using FADEs as opposed to more general CTRW models. First, it is well-known \cite{Meerschaert2001} that FADEs can account for source and boundary terms and velocity dynamics can be easily included by an additional velocity equation, which leads to a velocity-concentration coupled system. \MG{Furthermore, even though not thoroughly explored, model fitting for FADEs is a computationally less challenging task than for general CTRWs, due to the smaller number of parameters to fit.}

\subsection{Evidence of fractional behavior in the presence of heterogeneity}\label{sec:subsurface-evidence}

In this section we provide two examples of fractional behavior of solute concentration. We start by considering a highly heterogeneous environment and then we show that even in circumstances where a classical behavior is expected, i.e. in the absence of heterogeneities, the macroscopic solute concentration behaves nonlocally \MG{and is described by a FADE.}

\begin{figure}[t]
\centering
\includegraphics[width=0.9\textwidth]{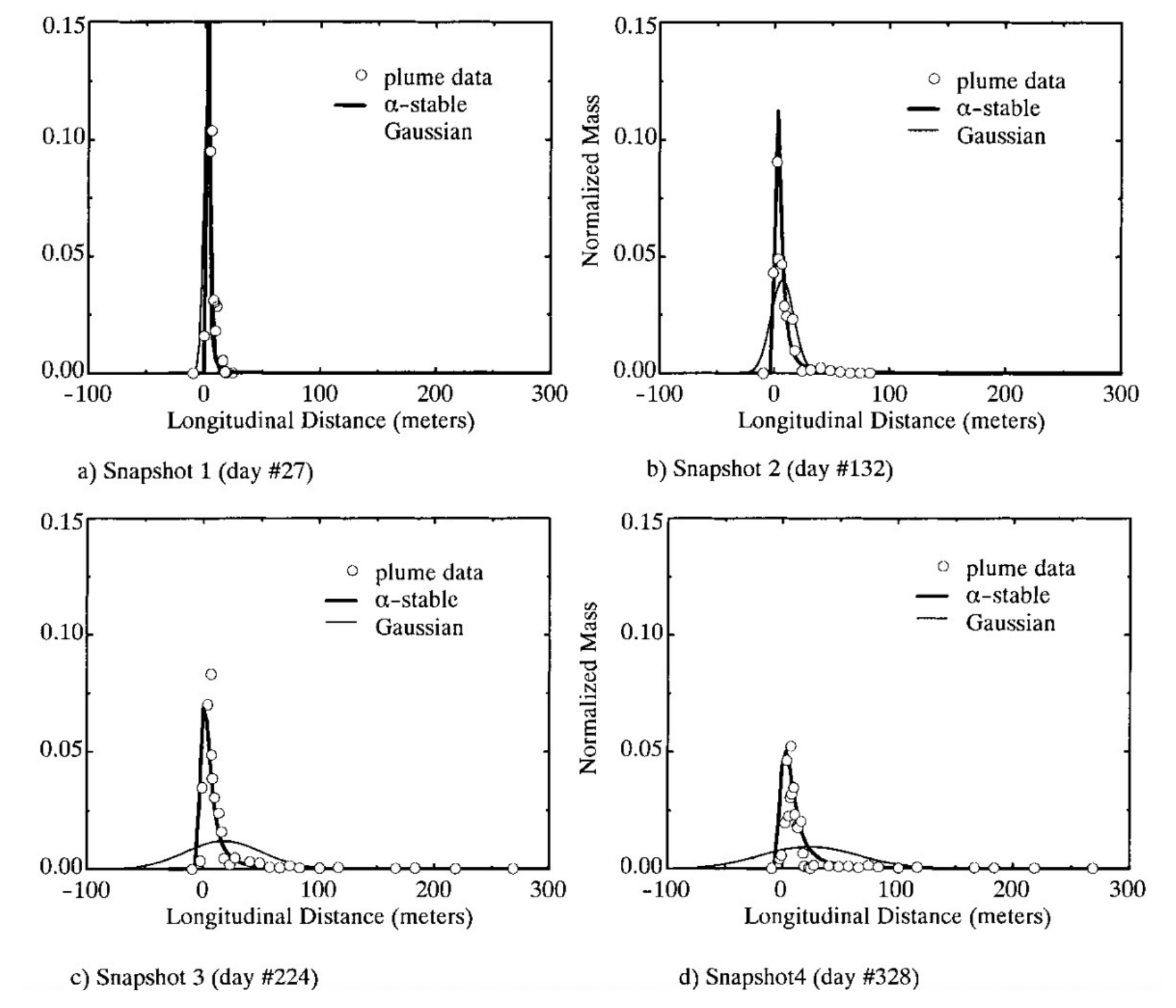}
\caption{A comparison of classical (Gaussian) and fractional ($\alpha$-stable) predictions of the normalized mass as a function of space at specific time instants for the MADE data set. The data points represent the maximum concentration measured in vertical slices perpendicular to the direction of the plume. These maxima were then integrated versus the travel distance. Source: \cite{Benson2001}.}
\label{fig:subsurface-evidence}
\end{figure}

Fractional behavior is most readily seen in transport through heterogeneous media. The first experiment we discuss studied subsurface transport of tritium in a highly heterogeneous environment such as the MADE site, located on the Columbus Air Force Base in northeastern Mississippi. This unconfined, alluvial aquifer consists of generally unconsolidated sands and gravels with smaller clay and silt components. Irregular lenses and horizontal layers were observed in an aquifer exposure near the site \cite{Rehfeldt1992}. Detailed studies characterizing the spatial variability of the aquifer and the spreading of the conservative tracer plume for the experiment conducted at the beginning of the 90's can be found in \cite{Boggs1993}. \citet{Benson2001} used \eqref{eq:simple-1d-fade} to model particle concentration; model parameters \MG{were} determined a priori by tuning them on the basis of measurements (we refer to \cite{Benson2001}, Sections 4.2 and 4.3, for a detailed description of the calibration process). In Figure \ref{fig:subsurface-evidence} we report four snapshots of the normalized longitudinal tritium mass distribution. These plots are obtained by numerical integration of the analytic solutions of both the classical ADE and the FADE.  These distributions clearly indicate that the fractional model outperforms the classical one. 

\begin{figure}[t]
\centering
\includegraphics[width=1\linewidth]{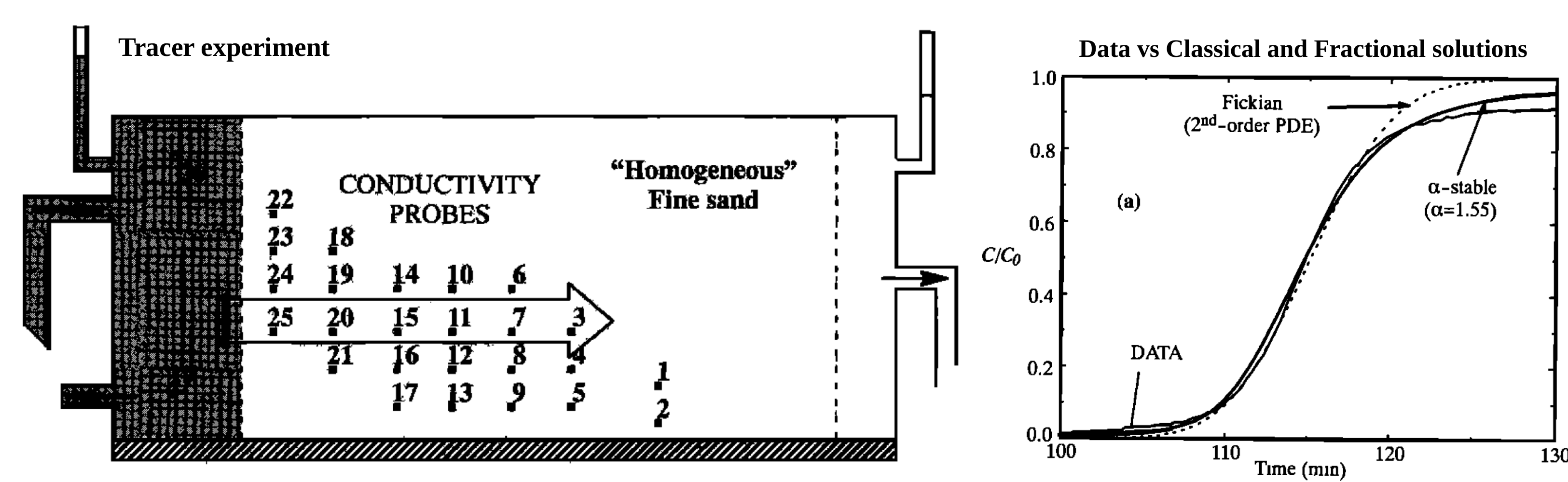}
\caption{Tracer transport in homogeneous sand shows evidence of anomalous behavior which can be reproduced by a fractional diffusion-advection equation. On the left, the set up of the homogeneous sand tracer experiment (as described in \cite{Burns1996}). On the right, a comparison of the corresponding classical (Fickian) and fractional ($\alpha$-stable) predictions. Source: \cite{Benson2000}.}
\label{fig:homo-sand}
\end{figure}

Strong heterogeneity, however, is not necessary to observe fractional behavior. Increasing experimental evidence suggests that in laboratory experiments where the media is ``constructed'' as nearly homogeneous, the observations \MG{are nevertheless} consistent with anomalous transport, see, e.g., \cite{Levy2003,Zhang2007}. In fact, some authors even \MG{argue} that \MG{transport in the subsurface is always anomalous} \cite{Zhang2013}. \citet{Benson2000} analyzed a test case where the tracer's concentration was intuitively expected to follow a classical ADE. They considered a one-dimensional tracer test in a laboratory-scale, 1 \MG{meter}, sandbox, constructed with very uniform sand in an effort to minimize heterogeneity, see Figure \ref{fig:homo-sand}, left. In other words, the sandbox was designed and built using as homogeneous a porous medium as possible by following the set up in \cite{Burns1996}. \MG{S}imple tracer tests, conducted to estimate the transport characteristics of the sand, indicated the appearance of non-classical breakthrough curves \MG{(BTCs, i.e., plots of the concentration in time at a fixed point downstream of a source)} with heavy tails, similar to $\alpha$-stable solutions. This behavior was likely due to channeling within smaller and smaller grains that resulted from sand emplacement through standing water and from cracked and intact surface clays on the sand particles \cite{Benson2000}. In Figure \ref{fig:homo-sand}, right, a comparison conducted in \cite{Benson2000} between BTCs obtained with the classical ADE and the FADE equation shows the agreement of the latter with measured BTCs at a specific location. While in this Figure the differences between classical and fractional behavior are not striking as in Figure \ref{fig:subsurface-evidence}, they are still noticeable.

We also mention that evidence of anomalous behavior and its successful description by FADEs has been observed in unsaturated soils \cite{Zhang2005}, saturated porous media \cite{Zhou2003}, streams and rivers \cite{Deng2004,Deng2006}, and overland solute transport due to rainfall \cite{Deng2006overland}.

\subsection{State of the art: a progression of fractional models for subsurface transport}\label{sec:subsurface-sota}

As described at the beginning of this section, classical diffusion does not take into account long-distance spatial and time correlations. The anomalous movement of particles in the subsurface, however, depends on both far upstream/downstream concentrations (resulting in space-fractional equations \citep{cushman1991diffusion,cushman1994nonequilibrium,Schumer2001,Zhang2007}) or past conditions (resulting in time-fractional equations \citep{cushman1994nonequilibrium,dentz2006delay,dentz2004time,Schumer2003}). Considering only the movement of solute particles in an infinitesimal neighborhood, \MG{as} in the classical diffusion model \MG{of} Brownian motion, is too restrictive for the complexities of groundwater pore spaces or trapping zones in natural streams. More specifically, the presence of preferential paths in hydrologic domains results in high-velocity zones (superdiffusion), whereas the presence of trapping regions results in low-velocity zones where the particles ``wait'' before they return to the higher velocity zone (this concept is also known in the literature as the distinction between immobile and mobile zones) \citep{Sun2020}.

In this section we review fractional models of increasing complexity for anomalous subsurface transport. While the simpler models are viable choices in the presence of a low degree of heterogeneity, as this degree increases, more sophisticated models are required to obtain reliable predictions. We first present early works featuring a one-dimensional space FADE with constant coefficients and constant fractional order. Next, we extend this model to the case of variable coefficients and generalize it to the multi-dimensional setting. We then present two types of one-dimensional time FADEs and conclude the section with a very general model featuring both space and time fractional derivatives of variable order. For all these models, we refer to Section \ref{sec:classification} for their mathematical details and interpretation in the context of stochastic processes.

\subsubsection{Spatial fractional derivatives}\label{sec:subsurface_space_frac}
We introduce the constant-coefficients, constant-order spatial FADE in one dimension introduced in \cite{Benson2000} and provide details regarding its parameters in relation to solute transport. The solute concentration at point $x$ and time $t$, $c(x,t)$, satisfies the equation
\begin{equation}\label{eq:space-FADE-const-1d}
\frac{\partial c}{\partial t}(x,t) 
= -V\dfrac{\partial c}{\partial x} 
-D \left[p \,
\left( \phantom{}^{\text{RL}}_{{-\infty}} \mathbb{D}^{\alpha}_x c(x,t) \right)
+
(1-p) \,
\left( \phantom{}^{\text{RL}}_{\phantom{R}x} \mathbb{D}^{\alpha}_{\infty} c(x,t) \right)
\right]
\end{equation}
where $V$ is the average plume velocity, $D$ is a fractional diffusion coefficient\footnote{The units of $D$ are given by ${\rm L}^\alpha/{\rm T}$ where $\alpha$ is the fractional order, L indicates space and T indicates time.} that controls the rate of spreading, $1\leq\alpha\leq 2$  (dimensionless) is the fractional order, and $0\leq p \leq 1$ determines the skewness $\gamma = 2p-1$. Solutions can be positively ($p=0$) or negatively ($p=1$) skewed, whereas they are symmetric when $p=0.5$, for which the sum of the Riemann-Liouville derivatives results in the fractional Laplacian 
The fractional
order $\alpha$ codes for the heterogeneity of the velocity field, with a higher probability
of large velocities as it decreases towards one \cite{Clarke2005}. We recall that for $\alpha=2$ the FADE reduces to the traditional advection-diffusion equation
(ADE) for groundwater flow and transport. The FADE above was introduced for the first time by \citet{Benson2000levy} to model scale-dependent dispersivity in
fitted groundwater plumes. In this paper the authors observed that, given a data set of solute concentration, the fitted parameter $D$ grows with time when the classical ADE is used; such evidence of superdiffusion is an indicator that a space-fractional model is preferable. Indeed, in subsequent works, see e.g. \cite{Benson2001}, the same authors show that the FADE allows the same data set to be fit with a constant-coefficient model such as \eqref{eq:simple-1d-fade}, where $D$ does not vary over time. From a particle perspective, the combination of left-sided and right-sided RL derivatives allows a solute particle to jump to any point in the domain; this simple concept was used by \citet{Schumer2001} to provide a derivation of \eqref{eq:simple-1d-fade} using an Eulerian interpretation of the particles' behavior.

\paragraph{The Gr\"unwald-Letnikov discretization technique.}
A standard discretization technique used in the FADE community for the approximation of the left-sided and right-sided RL derivatives \eqref{eq:leftRL} and \eqref{eq:rightRL} in \eqref{eq:simple-1d-fade}  is the shifted Gr\"unwald-Letnikov (GL) finite difference formula introduced by \citet{Meerschaert2004}. The GL scheme is based on the following identities:
\begin{equation}\label{eq:GLformula}
\begin{aligned}
^{\text{RL}}_{-\infty} \mathbb{D}^{\alpha}_x u(x,t) &=
\lim_{h\to 0} h^{-\alpha} \sum_{j=0}^\infty g_j^{\alpha} u(x+(j-1)h,t) \\[2mm]
\phantom{}^{\text{RL}}_{\phantom{R}x} \mathbb{D}^{\alpha}_\infty u(x,t) &=
\lim_{h\to 0} h^{-\alpha} \sum_{j=0}^\infty g_j^{\alpha} u(x-(j-1)h,t),
\end{aligned}
\end{equation}
where the GL weights are given by
\begin{equation}
g_j^\alpha = (-1)^j \dfrac{\Gamma(\alpha+1)}{\Gamma(j+1)\Gamma(\alpha-j+1)}.
\end{equation}
The GL approximation of the one-dimensional FADE is obtained by truncating the summation in \eqref{eq:GLformula}. The temporal derivative and the classical first-order spatial derivative can be obtained by standard time discretization schemes for PDEs. \MG{The discretizations \eqref{eq:GLformula} lead to stable time-stepping schemes and clearly highlight the nonlocal nature of fractional derivatives. Such discretizations also illustrate the potential for higher computational and memory cost entailed by solving FPDEs using direct approaches, which are greatly reduced using more the advanced numerical methods described in Section \ref{sec:numerical_methods}.}

\paragraph{FADEs with variable coefficients on bounded domains.}
In a heterogeneous porous medium, at a scale where the geological character of the medium changes with location, the constant-coefficient model \eqref{eq:simple-1d-fade} is insufficient for accurate and reliable predictions. A first step towards a more accurate model is introducing space dependence in the material parameters $V$ and $D$. Furthermore, in practical settings, simulations of solute transport must be confined to bounded domains, so that it becomes mandatory to establish ways to prescribe {\it nonlocal} boundary conditions that guarantee existence and uniqueness of solutions. In the literature there are at least three variants of the FADE with space-dependent coefficients \cite{Kelly2019chapter}: the fractional-flux ADE (FF-ADE), the fractional-divergence ADE (FD-ADE), and the fully fractional divergence ADE (FFD-ADE). In this review we focus on the former because of its resemblance \MG{to} classical advection-diffusion equations\footnote{We refer to \cite{Kelly2019chapter} for more details regarding the FD-ADE and FFD-ADE models.}.
\MG{For this model, we formulate the associated equation} on bounded domains.

The FF-ADE model in the one-dimensional domain $(-L,L)$ is derived from the classical conservation of mass equation
\begin{equation}\label{eq:FF-ADE}
\dfrac{\partial }{\partial t}c(x,t) + \dfrac{\partial }{\partial x}q(x,t) = 0, \;\;{\rm for}\;x\in(-L,L),
\end{equation}
where the {\it flux} $q$ is given by the following constitutive equation \cite{Schumer2001}
\begin{equation}\label{eq:q-flux}
q(x,t)= V(x) c(x,t) + D(x)
\Big[
p \big(\phantom{}^{\text{RL}}_{-L} \mathbb{D}^{\alpha-1}_x c(x,t)\big)
-(1-p)\big(
\phantom{}^{\text{RL}}_{\phantom{R}x} \mathbb{D}^{\alpha-1}_L c(x,t)\big)
\Big].
\end{equation}
Here, the first term is the advective flux that models the average drift of contaminant particles, whereas the second and third terms are the dispersive fluxes, which model large particle jumps in the left and right directions, respectively. Note that, because we consider the bounded domain $(-L,L)$, the integrals in the the left- and right-sided derivatives are ``truncated'' at $-L$ and $L$, respectively. Furthermore, since $\partial(\phantom{}^{\text{RL}}_{-L} \mathbb{D}^{\alpha-1}_x c(x,t))/\partial x = \phantom{}^{\text{RL}}_{-L} \mathbb{D}^\alpha_x c(x,t)$ the RL derivatives in the definition of the flux $q$ have exponent $\alpha-1$. The resulting FF-ADE corresponds to the models proposed in, e.g., \cite{Zhang2006}. We point out that, as described in detail in \cite{Kelly2019JCP}, Caputo derivatives as the ones introduced in Section \ref{sec:classification} can also be used in place of RL derivatives in the definition of the flux (leading to what is referred to as Caputo flux). 

\smallskip
The restriction of the FADE to a bounded domain requires the prescription of appropriate boundary conditions to guarantee that equation \eqref{eq:FF-ADE} is well-posed. We consider two types of boundary conditions: reflecting and absorbing. Using the flux function defined in \eqref{eq:q-flux}, we can identify a reflecting (or no-flux) condition by setting the diffusive part of the flux $q$ equal to zero at the boundary, i.e. $x=\pm L$. As an example, the reflecting boundary condition on the right boundary corresponds to  $$\phantom{}^{\text{RL}}_{\phantom{R}0} \mathbb{D}^{\alpha-1}_L c(L,t)-\phantom{}^{\text{RL}}_{\phantom{R}L} \mathbb{D}^{\alpha-1}_1 c(L,t)=0.$$ 

Instead, absorbing boundary conditions correspond to prescribing a zero ``Dirichlet'' condition at the boundary, i.e., 
$$c(\pm L,t)=0.$$ 
Clearly, these conditions can be mixed resulting in absorbing/reflecting boundary conditions on either the left or right boundary of the domain. It is important to note that, in the absence of advection, the no-flux (reflecting) condition implies that the total mass is conserved, see Proposition 2.3 \MG{by \citet{Kelly2019JCP}}.
We also mention that a new space-fractional model with variable advection and diffusion coefficients for anomalous, anisotropic transport has been proposed by \MG{\citet{d2021analysis}.}

\paragraph{Multidimensional FADEs}
The multidimensional version of equation \eqref{eq:simple-1d-fade} was proposed by \citet{Meerschaert1999} and further analyzed in \cite{Benson2000}. For $(-\Delta)_M^{\alpha/2}$ defined as in \eqref{eq:directional_laplacian}, we have that for $\xb\in\mbRd$ the concentration of a solute is described by the following law:
\begin{equation}\label{eq:multid-FADE}
\dfrac{\partial}{\partial t}c(\xb,t) 
+{\bf V}\cdot\nabla c(\xb,t)
-D(-\Delta)_M^{\alpha/2}c(\xb,t)= 0,
\end{equation}
where $\bf V$ is the average solute velocity and $D$ is the fractional diffusion coefficient. In \cite{Meerschaert1999} the operator $(-\Delta)_M^{\alpha/2}$ corresponding to  \eqref{eq:directional_laplacian}, is introduced via inverse Fourier transform, i.e. 
$$
(-\Delta)_M^{\alpha/2}c(\xb,t) = \mathcal F^{-1} \left\{ \int_{|\thetab|=1} (i {\bf k}\cdot\thetab)^\alpha \widehat c({\bf k},t)M(d\thetab)\right\}.
$$
Here, $\thetab$ is a $d$-dimensional unit vector, $\bf k$ is the wave vector and $\widehat c$ is the spatial Fourier transform of $c$. Note that the coefficient $D$ can be embedded in the measure $M$ (even when it depends on the space variable). As for the one-dimensional constant-coefficient equation \eqref{eq:simple-1d-fade}, the multidimensional FADE can also be extended to the variable-coefficient case. Furthermore, in the special case of jumps occurring only along the standard coordinate vectors ${\bf e}_j$, it is possible to derive fundamental solutions to \eqref{eq:multid-FADE}. Finally, the special case of uniform measure over the $d-1$ unit sphere corresponds to an advection-diffusion equation where the diffusion term is given by the standard fractional Laplacian operator $(-\Delta)^{\alpha/2}$.

\subsubsection{Temporal fractional derivatives}\label{sec:time_fractional_subsurface}

The time-FADE, used to model particle trapping in heterogeneous porous media, is characterized in a jump process perspective by long waiting times between jumps. \MG{As discussed in Section \ref{sec:subdiffusion}, this equation} replaces the first-order time-derivative in an ADE with a time-fractional derivative of either RL or Caputo type. In this section, we review two popular time-FADEs: the time-FADE (with RL derivatives) and
the fractional mobile-immobile equation (with Caputo derivatives), also known as FMIM.

\paragraph{Time-fractional advection-diffusion equation}
The time-fractional advection-diffusion equation (time-FADE) was introduced in the works by \citet{Zaslavsky1994} and, independently, by \citet{Liu2003}. In one dimension, it is given by
\begin{equation}\label{eq:time-FADE}
\phantom{}^{\text{C}}_{{0}} \mathbb{D}^{\alpha}_t c(x,t) = -v\dfrac{\partial}{\partial x}c(x,t)+D \dfrac{\partial^2}{\partial x^2}c(x,t),
\end{equation}
where the first term is the Caputo derivative defined in \eqref{eq:caputo_definition} on the half-axis. The units of the velocity parameter $v$ are ${\rm L}/{\rm T}^\alpha$ and the ones of the diffusion coefficient $D$ are ${\rm L}^2/{\rm T}^\alpha$, where L denotes units of space and T units of time. Note that, in the literature, $\phantom{}^{\text{C}}_{0} \mathbb{D}^{\alpha}_t f(t)$ is often denoted by $\frac{\partial^\gamma}{\partial t^\gamma}f(t)$, where $\gamma$ plays the same role as $\alpha$.
Furthermore, as pointed out at the beginning of this section, the time-FADE can be seen as the scaling limit of a CTRW. It is possible to obtain representations of solutions to \eqref{eq:time-FADE} by subordination, i.e. via randomization of the time variable by the inverse stable subordinator \cite{Meerschaert2013}.

\paragraph{Fractional mobile-immobile equation}
The fractional mobile-immobile (FMIM) model proposed by \citet{Schumer2003FMIM} is a generalization of the classical mobile-immobile (MIM) model \cite{Coats1964}. The latter, in its classical definition, partitions the solute concentration into a mobile phase, $c_{\text{m}}$, and an immobile phase, $c_{\text{im}}$ and equates the divergence of the total flux of the mobile concentration to a weighted sum of the time rate of change of each phase, i.e.
\begin{equation}\label{eq:MIM}
\dfrac{\partial}{\partial t}c_{\text{m}}(x,t) +
\beta \dfrac{\partial}{\partial t}c_{\text{im}}(x,t) = 
-v\dfrac{\partial}{\partial x}c_{\text{m}}(x,t)+D \dfrac{\partial^2}{\partial x^2}c_{\text{m}}(x,t),
\end{equation}
where $\beta=\eta_{im}/\eta_m$, being $\eta_{im}$ and $\eta_m$ the porosities of the immobile and mobile phases. The relationship between $c_{\text{m}}$ and $c_{\text{im}}$ is then given by one or more coupled mass transfer equations, resulting in the following relationship
\begin{equation}\label{eq:MIMrelation}
\dfrac{\partial}{\partial t}c_{\text{im}}(x,t)=
f(t)\star \dfrac{\partial}{\partial t}c_{m}(x,t)
+ f(t) (c_{\text{m}}(x,0)-c_{\text{im}}(x,0)),
\end{equation}
where $\star$ indicates the convolution operation and $f(t)$ is a memory function. 
The FMIM model in \cite{Schumer2003FMIM} defines $f(t)$ as the power function $f(t)= t^{-\alpha}/\Gamma(1-\alpha)$ with $0<\alpha<1$. By noting that
$$
f(t)\star \dfrac{\partial}{\partial t}c_{m}(x,t) = \phantom{}^{\text{C}}_{{0}} \mathbb{D}^{\alpha}_t c_{\text{m}}(x,t),
$$
the combination of \eqref{eq:MIMrelation} and \eqref{eq:MIM} results in the time-FADE
\begin{equation*}
\dfrac{\partial}{\partial t}c_{\text{im}}(x,t)=
\phantom{}^{\text{C}}_{{0}} \mathbb{D}^{\alpha}_t c_{\text{m}}(x,t)
+ f(t) (c_{\text{m}}(x,0)-c_{\text{im}}(x,0)),
\end{equation*}

A CTRW model for the FMIM model, was developed by \citet{Benson2009}; here, waiting times experienced by solute particles in the immobile phase are modeled by a power law (as for the time-FADE). 
\MG{As regards the realism of the assumptions underlying this model,} power-law waiting times have also been observed in river transport studies by \citet{Haggerty2002} and \citet{Schmadel2016}.

\subsubsection{Variable-order FADEs}

Constant-coefficient and constant-order models are invaluable basic tools for the analysis of complex engineering systems such as the flow through the subsurface; however, they are unable to evolve between different physical behaviors, i.e. they cannot capture transitions between diffusive regimes. These transitions are caused by the fact that solutes in the subsurface diffuse through porous, fractured, layered and heterogeneous aquifers, whose structure changes with space as well as time. This leads to anomalous diffusion characterized by a variable-order scaling of the MSD. 
\MG{While the variable-coefficient models described above allow for more flexibility with regards to scaling properties, a more direct way to address this need is with variable-order models. Several recent works have explored the use of variable-order operators in the context of subsurface modeling and other applications}. The use of these operators becomes particularly important in the presence of complex media that feature a hybrid anomalous mechanism \cite{patnaik2020applications}. As an example, we can exploit variable-order fractional operators, like the ones introduced in Section \ref{sec:variable-order}, when the nature of the transport processes transitions across very different underlying physical phenomena such as transitions from sub-diffusive flow to diffusive flow, and from diffusive flow to super-diffusive flow \cite{Atangana2014,Kobelev2003,Sokolov2006,SUN2009,SUN2010,SUN2014}. Note that these complex transport processes have been observed experimentally in various fields; for fluid flow through porous media we mention, e.g., \MG{the works of \citet{Gerasimov2010} and \citet{Obembe2017}}. 
 
A complete variable-order fractional model was proposed \MG{by \citet{SUN2014} and further explored by \cite{Pang2017}} for the description of the same MADE data set introduced at the beginning of this section. The one-dimensional variable-order time-space FADE is given by
\begin{equation}\label{eq:FADE-VO-full}
    \phantom{}^{\text{C}}_{{0}} \mathbb{D}^{\beta(x,t)}_t c(x,t) = 
    -V\dfrac{\partial c}{\partial x} 
-D^-  \,
\left( \phantom{}^{\text{RL}}_{-\infty} \mathbb{D}^{\alpha(x,t)}_x c(x,t) \right)
- D^+ \,
\left( \phantom{}^{\text{RL}}_{\phantom{R}x} \mathbb{D}^{\alpha(x,t)}_{\infty} c(x,t) \right),
\end{equation}
where the variable-order derivatives are defined as in Section \ref{sec:variable-order}.

\begin{figure}[t]
\centering
\includegraphics[width=0.5\textwidth]{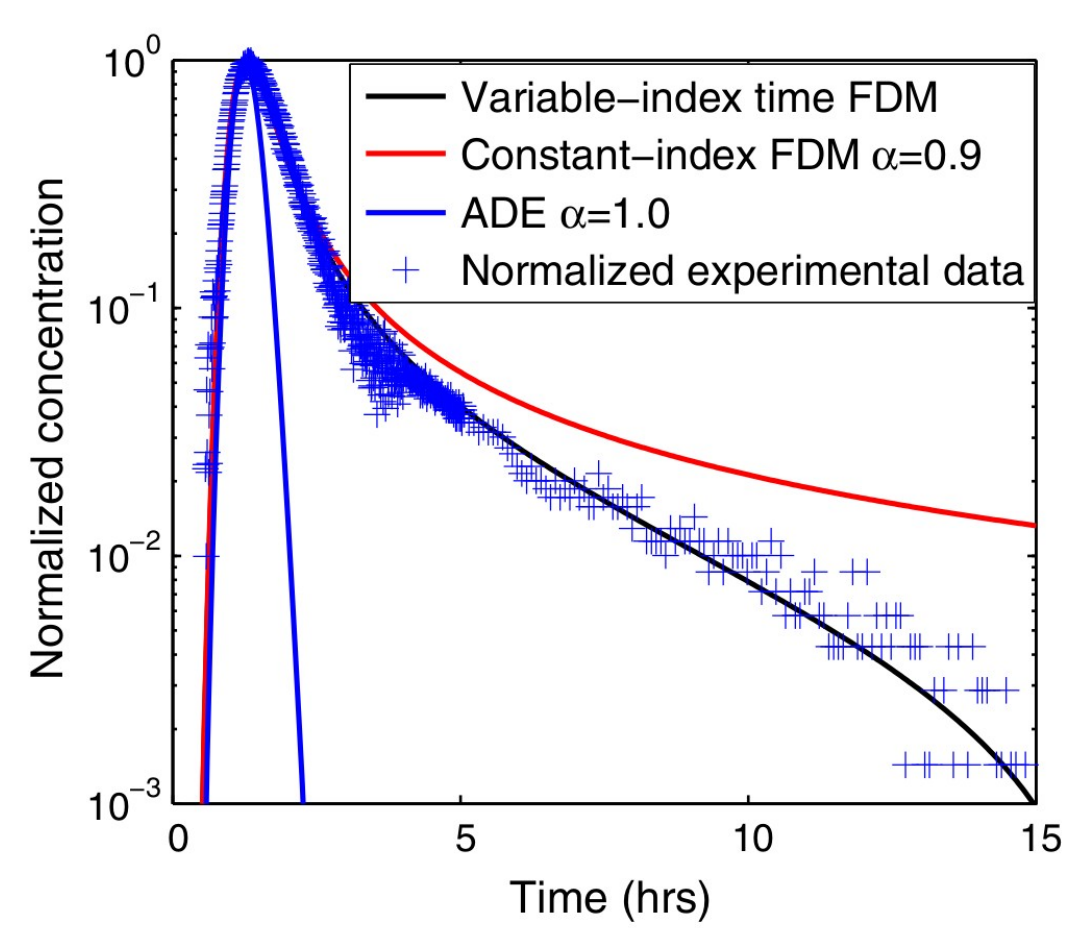}
\caption{A comparison in semi-log scale of the normalized concentration at the extraction point for the Grimsel test site \cite{Berkowitz2008} obtained using the classical advection-diffusion equation, the constant-order time FADE, and variable-order time FADE (ADE, Constant-index FDM, and Variable-index time FDM, respectively in the legend) together with the normalized experimental data. Source: \citet{SUN2014}. Note that their use of $\alpha$ and $\beta$ is switched from our use of the same symbols in the text; thus, in the above legend, $\alpha$ denotes the order of the time-fractional derivative.}
\label{fig:variable-order}
\end{figure}

To confirm the improved accuracy of models such as the one in \eqref{eq:FADE-VO-full} we report in Figure \ref{fig:variable-order} a comparison, conducted \MG{by \citet{SUN2014}}, of a classical model, a constant-order fractional model and a variable-order fractional model. Here, the authors consider concentration data from the field experiment conducted at the Grimsel test site
\cite{Berkowitz2008} where uranine, a fluorescent dye, was injected into a shear zone as a tracer and its concentration was measured at an extraction well away from the injection site. The BTC of uranine, measured at the extraction well corresponds to the blue crosses in the figure. The authors compare the following models: the classical advection-diffusion equation, corresponding to $\beta=1$ and $\alpha=2$ in \eqref{eq:FADE-VO-full}, the constant-order time FADE with $\beta=0.9$ and $\alpha=2$, and the variable-order time FADE with $\beta(t)=0.9 + t/150$, $t\in(0,15]$, and $\alpha=2$. 
BTCs in the figure show that the classical ADE model is not capable to depict the tailing/subdiffusive behavior, whereas the constant-order time FADE underestimates the late time decay, which features classical behavior. The choice of $\alpha$ and $\beta$ in the variable-order time FADE is based on the following considerations: first, the measured BTC has a fast-increasing early time tail, implying a Gaussian-type of particle jump that corresponds to $\alpha=2$. Second, the heavy late-time tail suggests a time-dependent $\alpha$ that should be less than 1 at early times (subdiffusive) and should slowly convergence to 1 at late times (classical diffusion). The corresponding solid black BTC clearly captures the variable diffusion behavior of the normalized concentration.

\subsection{Future directions in anomalous subsurface modeling}\label{sec:subsurface-future}
In the previous sections we provided evidence of the occurrence of anomalous behavior in subsurface transport even for a low degree of heterogeneity and we have shown that FADEs can be accurate models when properly tuned. However, the identification of an optimal fractional model for a specific setting (e.g. for specific hydraulic properties) is not trivial and has not been thoroughly explored in the literature. One of the main challenges in this context is the fact that model parameters cannot be directly related to media properties, as carefully explained \MG{by \citet{Neuman2009}}. Furthermore, often times, it is hard or nearly impossible to collect solute measurements, so that only a very small set of data that are sparse in time and space and potentially affected by noise is available. Yet, in this context, FADEs have the advantage, compared to other models for subsurface transport, of having only a handful of parameters to tune, i.e. the identification problem consists in discovering a small set of parameters such as the diffusivity and the fractional order. 

Only a few works in the literature have addressed this problem. In the context of highly heterogeneous settings, we mention the work by \citet{Pang2017} where the authors propose to use multi-fidelity Bayesian optimization to discover  variable-order fractional operators for the advection–diffusion equation \eqref{eq:FADE-VO-full} from field data in the MADE data set mentioned at the beginning of this section. Other recent works addressing a similar learning problem for fractional operators include optimization-based approaches such as the one used \MG{by \citet{Burkovska2021}}, fractional/nonlocal physics-informed neural network approaches such as \MG{those of \citet{Pang2019,Pang2020}} and operator-regression techniques such as the one developed \MG{by \citet{You2020Regression}}. \MG{It is important to keep in mind \MD{that} the inverse problem of parameter identification may multiply the cost of numerical \MD{computations}, which is already higher in direct approaches than classical methods due to the integral nature of the operators involved and to the strong singularities that require adaptive quadrature rules. Thus, together with the development \MD{of} frameworks for learning models, the efficient numerical methods mentioned in Section \ref{sec:numerical_methods} take on an additional importance in model discovery.}


\section{Turbulence}\label{sec:turbulence}

Richard P. Feynman described turbulence as the most important unsolved problem in classical physics \cite{Eames}, a problem that stands today. 
By ``turbulence'', we refer to the three-dimensional and highly vortical fluid motions characterized by stochastic perturbations in pressure and flow velocity, and caused by excessive kinetic energy in areas of fluid flow that overcome the ``damping effects'' of the fluid's viscosity. The onset of turbulence can be predicted by the dimensionless Reynolds number \textit{Re}, a ratio of kinetic energy to viscous damping in the fluid flow. Yet, the question remains of what mathematically governs the evolution of a turbulent flow and whether it is feasible to fully simulate turbulent flows by means of numerical methods.

In 1970, \citet{emmons1970critique} reviewed the possibilities for computational fluid dynamics, concluding: ``\textit{... the problem of turbulent flows is still the big holdout. This straightforward calculation of turbulent flows -- necessarily three-dimensional and unsteady -- requires a number of numerical operations too great for the foreseeable future."} After almost a decade, however, the field of direct numerical simulation (DNS) of turbulence was established with successful numerical simulations of wind-tunnel flows at moderate \textit{Re} by \citet{hussaini2012instability,karniadakis1993nodes,KMM1987,orszag1972numerical}. These early computational developments were based on employing a Newtonian fluid assumption and applying the principles of conservation of mass, momentum, and energy to an infinitesimally small fluid element or parcel; see e.g., \cite{chorin1990mathematical, temam2001navier, white2006viscous}. This led to the derivation of the Navier-Stokes and energy equations, emerging as a set of convective nonlinear PDEs that govern the evolution of fluid velocity/temperature fields in turbulence. In this context, assuming some proper (random) initial/boundary conditions, one can discretize the governing equations and solve for the ``entire degrees of freedom of turbulence'' in the physical and parametric (stochastic) space. 

The great challenge is that in practice, DNS becomes prohibitively expensive, especially at high \textit{Re}, more so in complex geometries. Hence, one of the main goals in \textit{turbulence modeling} has been to systematically lower the total number degrees of freedom to a manageable level, at the cost of reducing the accuracy of turbulence predictions. This approach has been mainly centered around the overarching theme of ensemble averaging the set of PDEs representing the various scalar and vector turbulence fields; see e.g., \cite{lumley1995turbulence, wilcox1998turbulence}. This gives rise to new mathematical terms in the averaged or \textit{filtered} governing equations, known as \textit{turbulence closure} terms that can only be modeled as they are essentially unknown. When the entire time and length scales of turbulence are averaged, the Reynolds Averaged Navier-Stokes (RANS) equations are obtained, solely describing the \textit{mean-flow dynamics} of turbulence. Alternatively, if one applies a mathematically well-defined low-pass filter to the Navier-Stokes equations, the resulting filtered governing equations describe the \textit{large-eddy} dynamics of turbulence, where only small-scale subgrid dynamics need be modeled; this is referred to as large eddy simulation (LES). The common approach in the literature for modeling closure terms of any kind has been based on the use of classical local differential operators. More specifically, the majority of turbulence models have been constructed based on Boussinesq’s turbulent viscosity concept \citep{pope2001turbulent}, in which one assumes that the turbulent stress tensors are proportional to the \textit{local gradient} of mean velocity at any point. The proportionality coefficient, referred to as \textit{turbulent viscosity}, is to be inferred from data. 

The impetus for the fractional models we describe in this section is that the small-scale dynamics of turbulence are statistically anomalous, i.e., non-Markovian and non-Fickian, so that nonlocal closure models emerge as appropriate tools. Employed at the continuum level, fractional models therefore capture anomalous features in the small-scale stochastic subgrid dynamics of turbulence.
The mathematical modeling of turbulence must address the fact that nonlinear interactions between the turbulence structures and motions create statistically complex phenomena that lead to a variety of anomalous features, including multi-power-law scalings in space-time, rare events, short-to-long range coherent motions, and enhanced turbulent mixing. These features urge better and novel understanding of the underlying nonlocal closure terms that appear as a result of the ensemble averaging or filtering of the governing equations. The nonlocal mode of thinking has the potential to shift the turbulence modeling paradigm and achieve a new level of physical and statistical consistency compared to classical approaches. This is especially true at high \textit{Re}, for which a proper and efficient framework that unites computational, mathematical, and statistical aspects was not available until recently.

\subsection{Evidence of Fractional Behavior in Turbulence}
An intuitive concept of nonlocality and memory effects was been established by \citet{Eringen}, where a point within a fluid field (medium) is influenced by all points of the body at all past times. Coherent random motions and the spatially turbulence spots structures inherently give rise to \textit{intermittent} signals with \textit{self-similarities}, \textit{sharp peaks}, \textit{heavy-skirts}, and \textit{skewed} distributions of velocity increments. Such statistical features have been well-observed experimentally even in the context of most canonical problems, e.g., grid turbulence, in which the skewness factor negatively appears and the Kurtosis factor strongly exceeds three, emphasizing the non-Gaussian character of statistics (see e.g., \cite{davidson2015turbulence}). Moreover as demonstrated by \citet{Egolf_book} (page 92), nonlocal effects appear even in the context of turbulent fields obtained numerically solving the Navier-Stokes equations.  Such widespread statistical measures indicate the non-Markovian and non-Fickian nature of turbulence, and they are the consequence of nonlinear and coherent vortical effects that occur in a wide spectrum of length and time scales. Therefore, nonlocal interactions cannot be ruled out of modern turbulence physics. These considerations are particularly timely; in fact, we can now benefit from the spectrum of modern nonlocal and fractional modeling tools reviewed in Section \ref{sec:classification}, equipped with well-established mathematical/statistical theories, that enable us to take such nonlocal/history effects into account with physical consistency and mathematical rigor.  

Furthermore, averaging entire spatial scales as in RANS models or applying a spatial filter to the Navier-Stokes and energy equations as in LES models would make the underlying physical nonlocality in the corresponding closure terms in RANS models and the subgrid turbulent fluctuations in LES models even more pronounced. This sheds lights on why turbulence modeling is a nonlocal task and further motivates the development of ``nonlocal closure models'' that can properly address and incorporate the underlying memory and long-range effects. Specifically, in what follows, we present a DNS study, recently presented by \citet{Scalar_FSGS}, that introduces new statistical measures and highlights the nonlocal character of subgrid scale dynamics in the context of scalar turbulence.

\subsubsection{The case of scalar turbulence subgrid dynamics}

An ideal LES is such that the true, filtered turbulent intensity is captured accurately through a robust subgrid scale (SGS) modeling that is physically and mathematically expressive. In fact, the LES equations include closure terms that directly link the correct evolution in time of turbulent intensity to the nature of the SGS closure and its modeling. Here, as a canonical problem, we consider the advection-diffusion (AD) equation 
\begin{eqnarray}\label{GE-1-3}
    \frac{\partial \phi}{\partial t}+\frac{\partial}{\partial x_i}\left( \phi \, V_i \right) = -\theta \, V_2+\mathcal{D} \, \frac{\partial^2 \phi}{\partial x_i \partial x_i}, \quad i=1,2,3,
\end{eqnarray}
in which $\mathcal{D}$ denotes the molecular diffusion coefficient of the passive scalar, and the imposed mean scalar gradient is taken to be uniform as $\nabla  \langle \Phi \rangle = \left( 0,\theta,0\right)$, where $\theta$ is a real-valued constant. In the LES representation of the scalar turbulence, multiplying both sides of the ``filtered'' AD equation by $\widetilde{\phi}$, the filtered scalar field ${\phi}$, yields the time-evolution of the filtered turbulent \textit{intensity} as
\begin{equation}\label{eqn: scalar_var1}
\frac{1}{2} \frac{\partial}{\partial t}\left( \widetilde{\phi} \, \widetilde{\phi} \right) + \widetilde{\phi} \, \frac{\partial}{\partial x_i} \left( \widetilde{\phi} \, \widetilde{V}_i \right)= -\theta \, \widetilde{\phi} \, \widetilde{V}_2 + \mathcal{D} \, \widetilde{\phi} \, \frac{\partial^2 \, \widetilde{\phi}}{\partial x_i \partial x_i} - \widetilde{\phi} \, \frac{\partial \, q^R_i}{\partial x_i}.
\end{equation}
Here, $q_i^R$ denotes the $i$-th component of the residual, SGS scalar flux defined as
$q_i^R = \widetilde{\phi V_i} - \widetilde{\phi} \widetilde{V_i} $.
Employing the filtered continuity equation $\nabla \cdot \boldsymbol{\widetilde{V}}=0$ and the chain rule for differentiation, we obtain
\begin{eqnarray}\label{eqn: scalar_var2}
\frac{1}{2} \frac{\partial}{\partial t}\left( \widetilde{\phi} \, \widetilde{\phi} \right) + \widetilde{\phi} \, \widetilde{V}_i \, \frac{\partial \widetilde{\phi}}{\partial x_i} &=& -\theta \, \widetilde{\phi} \, \widetilde{V}_2 + \mathcal{D} \, \frac{\partial}{\partial x_i}\left( \widetilde{\phi} \, \frac{\partial \widetilde{\phi}}{\partial x_i} \right) - \mathcal{D} \, \frac{\partial \, \widetilde{\phi}}{\partial x_i} \, \frac{\partial \, \widetilde{\phi}}{\partial x_i} 
\\
\nonumber
&\phantom{=}& - \frac{\partial}{\partial x_i}\left( \widetilde{\phi} \, q^R_i \right) + q^R_i \, \frac{\partial \widetilde{\phi}}{\partial x_i}.
\end{eqnarray}
Applying the ensemble-averaging operator, $\langle \cdot \rangle$, on \eqref{eqn: scalar_var2} returns a transport equation for the \textit{filtered scalar variance}, $\left\langle \widetilde{\phi} \, \widetilde{\phi} \right\rangle$. \citet{Scalar_FSGS} considers the case of homogeneous turbulent velocity and scalar fields, in which $\left\langle \frac{\partial}{\partial x_i}\left(\cdot \right) \right\rangle = \frac{\partial}{\partial x_i}\langle (\cdot) \rangle = 0$. Defining the filtered scalar gradient as $\widetilde{\boldsymbol{G}}(\xb) = \nabla \widetilde{\phi}(\xb)$, the time-evolution of the filtered scalar variance takes the following form 
\begin{align}\label{eqn: scalar_var3}
\frac{1}{2} \frac{d}{d t}\left\langle \widetilde{\phi} \, \widetilde{\phi} \right\rangle &= -\widetilde{\mathcal{T}} + \widetilde{\mathcal{P}} - \widetilde{\chi} + \Pi, \\
\widetilde{\mathcal{T}} = \left\langle \widetilde{\phi} \, \widetilde{V}_i \, \widetilde{G}_i \right \rangle, \quad \widetilde{\mathcal{P}} = -\theta \left\langle \widetilde{\phi} \, \widetilde{V}_2 \right \rangle&, \quad \widetilde{\chi} = \mathcal{D} \, \left\langle \widetilde{G}_i \, \widetilde{G}_i \right\rangle, \quad \Pi = \left\langle q^R_i \, \widetilde{G}_i \right\rangle. \nonumber
\end{align}
%
\begin{figure}[t!]
	\begin{minipage}[b]{.48\linewidth}
		\centering
		\includegraphics[width=1\textwidth]{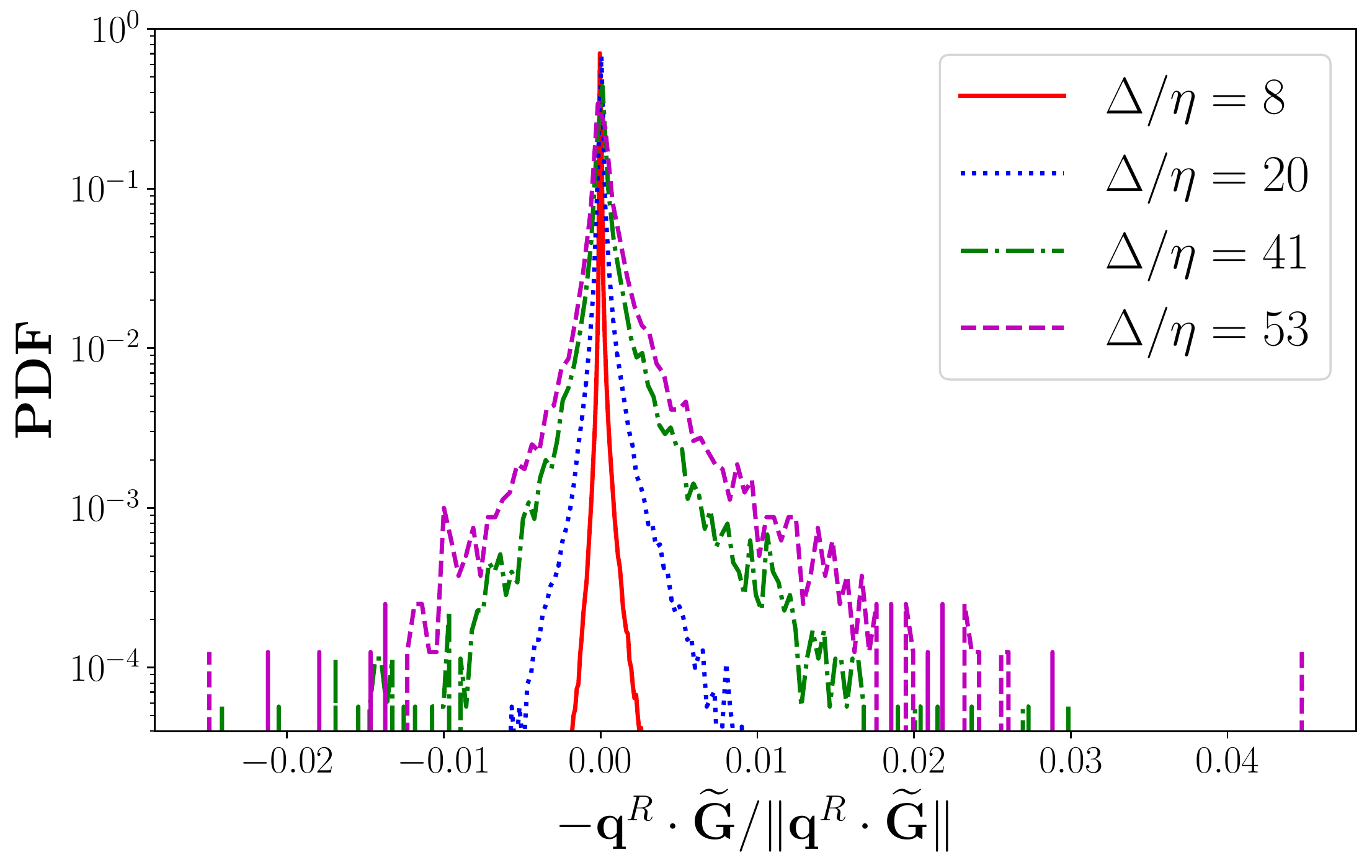}
		\subcaption{\footnotesize}
	\end{minipage}
	\begin{minipage}[b]{.01\linewidth}
		~
	\end{minipage}
	\begin{minipage}[b]{.49\linewidth}
		\centering
		\includegraphics[width=1\textwidth]{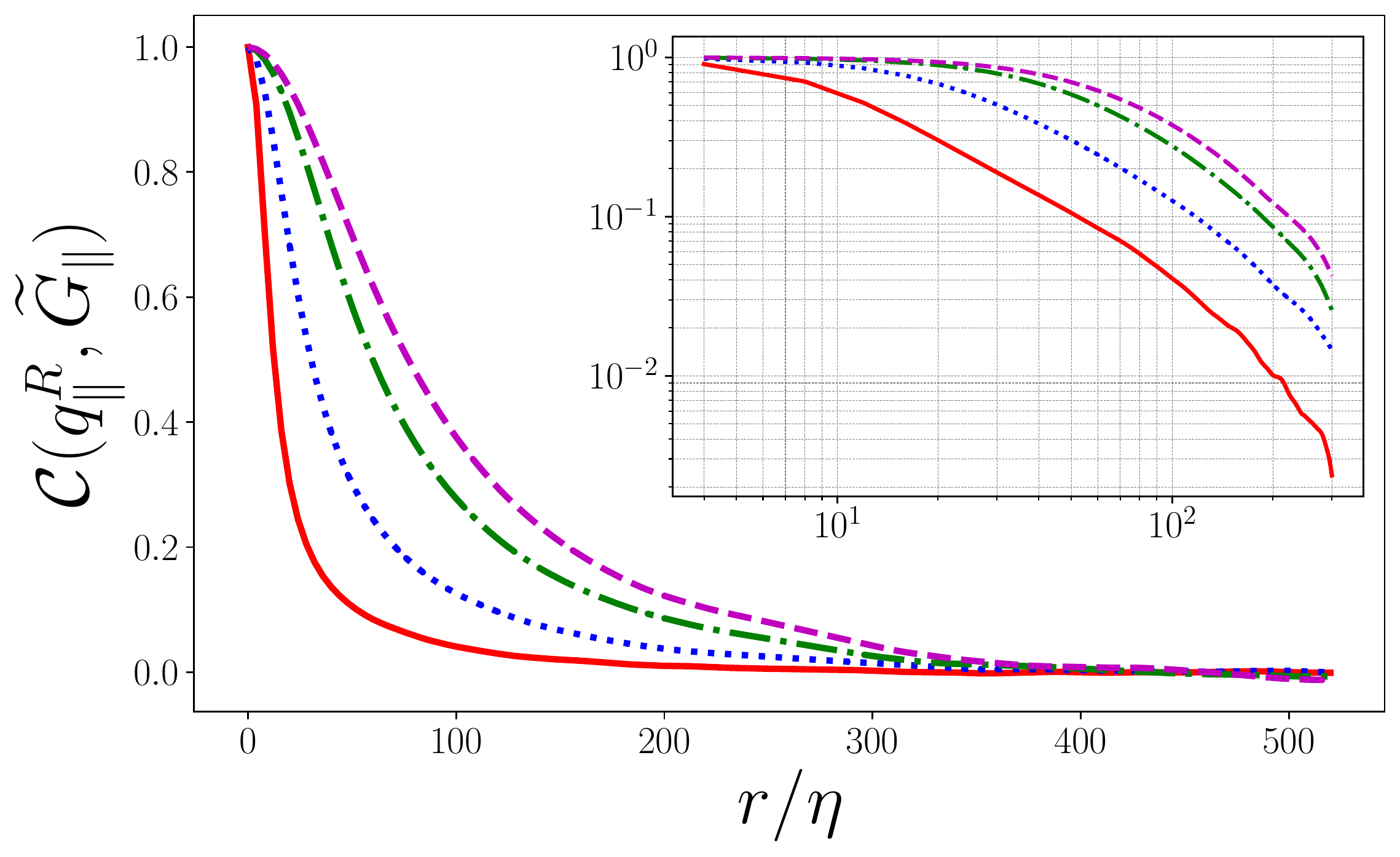}
		\subcaption{\footnotesize}
	\end{minipage}   
	\caption{Statistics of true subgrid-scale contribution to the filtered scalar variance rate. (a) PDF of normalized SGS dissipation of filtered scalar variance, $-\boldsymbol{q}^R \cdot \widetilde{\boldsymbol{G}}$, computed over a sample space of \MG{$10 \,T_{\text{LE}}$} of statically stationary turbulence. (b) Time-averaged two-point correlation function \eqref{eqn: TPC1} between $q^R_\parallel$ and $\widetilde{G}_\parallel$ with $r=\vert \boldsymbol{r}_\perp \vert$. Source: \cite{Scalar_FSGS}.}
	\label{fig: Nonlocality}
\end{figure}

In \eqref{eqn: scalar_var3}, $\widetilde{\mathcal{T}}$ denotes the turbulent transport of filtered scalar variance while $\widetilde{\mathcal{P}}$ represents the production of resolved scalar variance by the uniform mean scalar gradient, and $\widetilde{\chi}$ is the resolved scalar variance dissipation due to the molecular diffusion. Unlike these three terms, $\Pi$ (representing the SGS production of resolved scalar variance) is the only contributing term in \eqref{eqn: scalar_var3} that contains the effects of the SGS scalar flux. Therefore, as pointed out earlier, understanding the true statistical nature of $\boldsymbol{q}^R \cdot \widetilde{\boldsymbol{G}}$ is essential for the SGS modeling and precise evaluation of the resolved scalar variance in the LES. This examination of $\boldsymbol{q}^R \cdot \widetilde{\boldsymbol{G}}$ might be viewed both from single-point and two-point statistics as discussed by \citet{meneveau1994statistics} in the context of the LES for homogeneous isotropic turbulent flows. We focus on the two-point statistics of the SGS production of resolved scalar variance. This quantity is well represented in terms of the following normalized two-point correlation function
\begin{align}\label{eqn: TPC1}
\mathcal{C}(q^R_i \, , \, \widetilde{G}_i) = \frac{\left \langle q^R_i(\xb) \, \widetilde{G}_i(\xb+\boldsymbol{r}) \right \rangle}{\left \langle q^R_i(\xb) \, \widetilde{G}_i(\xb) \right \rangle},
\end{align}
where $\boldsymbol{r}=(r_1,r_2,r_3)$ denotes the spatial shift from the location $\xb$. Moreover, the probability density function (PDF) of the SGS production of scalar variance normalized by its $L_2$-norm i.e., $\boldsymbol{q}^R \cdot \widetilde{\boldsymbol{G}}/\Vert \boldsymbol{q}^R \cdot \widetilde{\boldsymbol{G}} \Vert$, is a novel statistical measure for studying the statistical behavior of $\Pi$ and yielding a more comprehensive insight into the SGS modeling. 

Let $T_{\text{LE}}$ be the eddy turnover time. By taking a large sample space over $10 \, T_{\text{LE}}$ of this stationary process (after resolving the passive scalar field for \MG{$15 \, T_{\text{LE}}$}), the PDF of the normalized SGS production of filtered scalar variance is computed for four different filter widths, $\Delta/\eta=8, \, 20, \, 41, \, 53$. These computations, shown in  Figure \ref{fig: Nonlocality}(a), demonstrate that as $\Delta$ becomes larger, the PDF exhibits broader tails. Emergence of this tail behavior implies that as the filter width increases, long-range spatial interactions become stronger and more pronounced \cite{akhavan2020anomalous}. Motivated by this observation, a two-point diagnosis of the SGS scalar production of the filtered variance as defined in equation \eqref{eqn: TPC1} would be another statistical measure shedding light on the long-range interactions in addition to the filter width effects. Considering $\parallel$ as the direction along the imposed mean scalar gradient and $\perp$ representing the directions perpendicular to the imposed mean gradient, we focus on the evaluation of $\mathcal{C}(q^R_\parallel \, , \, \widetilde{G}_\parallel)$.

Here, one case takes $\boldsymbol{r}=(r_1,0,0)$ and $\boldsymbol{r}=(0,0,r_3)$ and takes the average of the resulting two-point correlation functions. Due to the statistically stationary turbulence, such procedure is performed for 20 data snapshots that are uniformly spaced over \MG{$10 \, T_{\text{LE}}$} (on the same spatio-temporal data, used to compute the PDFs); hence, the time-averaged value of $\mathcal{C}(q^R_\parallel \, , \, \widetilde{G}_\parallel)$ is obtained. Figure \ref{fig: Nonlocality}(b) illustrates this two-point correlation function extending over a wide range of spatial shift, $r=\vert \boldsymbol{r} \vert$, and evaluated at four filter widths similar to the ones utilized in Figure \ref{fig: Nonlocality}(a). This plot quantitatively and qualitatively reveals that as $\Delta$ increases, greater correlation values between the SGS scalar flux $q^R_\parallel(\xb)$, and filtered scalar gradient $\widetilde{G}_\parallel(\xb+\boldsymbol{r})$ are observed at a fixed $r$. These spatial correlations are significant both in the \textit{dissipation} and also \textit{inertial} subranges. 

This confirms substantial nonlocal effects in the true SGS dynamics, which need to be carefully addressed in the SGS modeling for LES. 
\begin{figure}[t!]
	\centering
	\includegraphics[width=.5\textwidth]{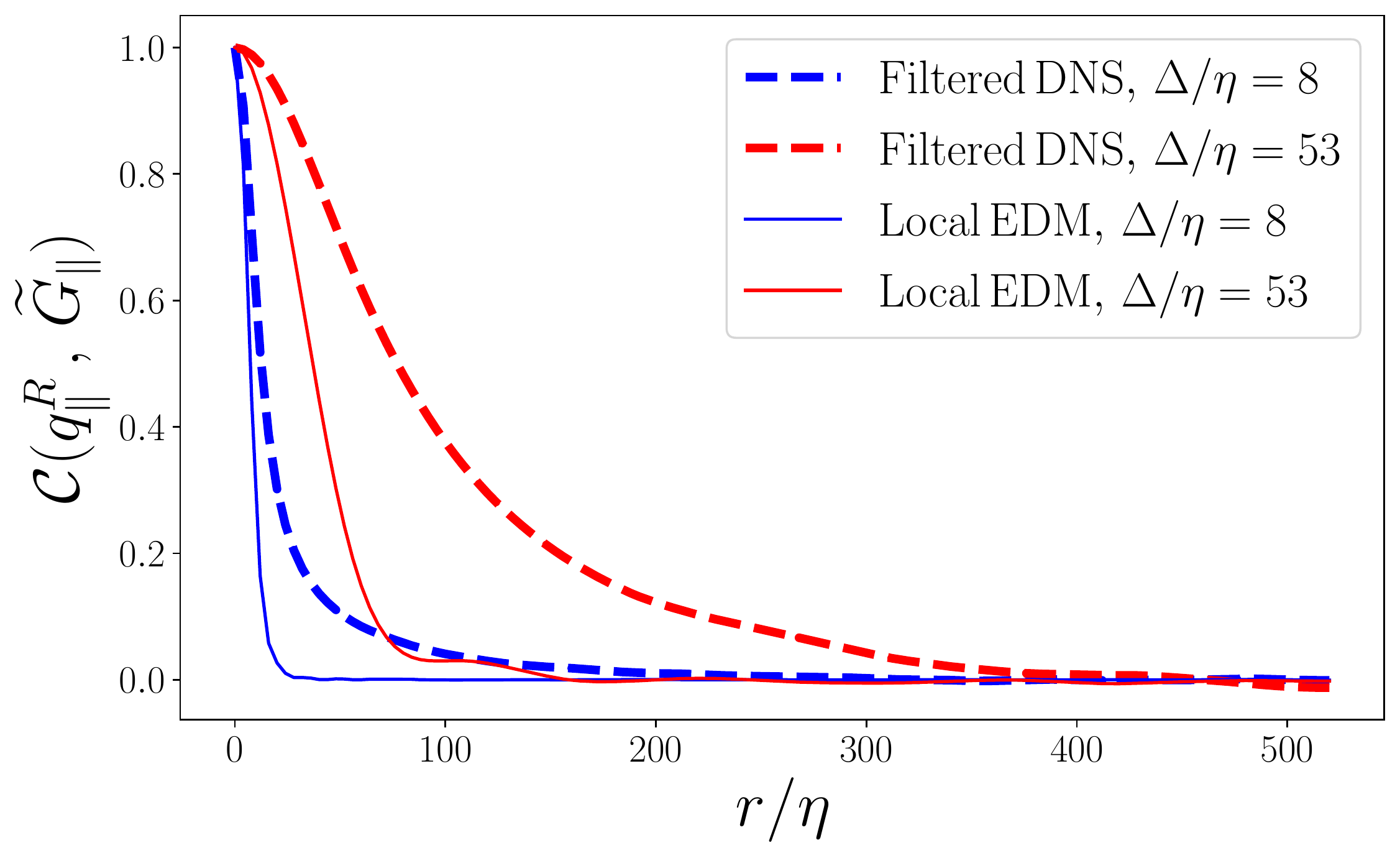}
	\caption{ Comparison between the true values of two-point correlation function given in \eqref{eqn: TPC1} and the ones obtained from the local eddy-diffusivity modeling of the SGS scalar flux given in \eqref{eqn: EDM}. The evaluations are performed at two filter widths of $\Delta/\eta = 8, \, 53$. Source: \cite{Scalar_FSGS}.}
	\label{fig: EDM_TPC}
\end{figure}
A popular and fairly simple approach for modeling the SGS scalar flux is Eddy-Diffusivity Modeling (EDM). In EDM, the main assumption is that the SGS scalar flux is proportional to the resolved scale scalar gradient (i.e., the conventional locality assumption) as
\begin{align}\label{eqn: EDM}
\boldsymbol{q}^R(\xb) \approx -\mathcal{D}_{\text{ED}} \, \widetilde{\boldsymbol{G}}(\xb),  
\end{align}
and $\mathcal{D}_{\text{ED}}$ is the proportionality coefficient. Obviously, EDM is a \textit{local} modeling approach by construction. Computing $\mathcal{C}(q^R_\parallel \, , \, \widetilde{G}_\parallel)$ while $q^R_\parallel$ is approximated with EDM, one can compare it with its true value as shown in Figure \ref{fig: Nonlocality}(b). Figure \ref{fig: EDM_TPC} illustrates such comparison for two filter widths, $\Delta/\eta=8, \, 53$, and it reveals that  
in both of the cases local EDM substantially fails to predict the conspicuous long-range spatial correlations observed in the true two-point correlation values. \textit{This evidence strongly suggests the adoption of more sophisticated, nonlocal mathematical modeling tool that goes beyond conventional SGS modeling.}

\subsection{State of the art: fractional turbulence modeling}
In what follows, we present the history and state of the art in fractional turbulence modeling, including nonlocal RANS, fractional eddy-viscosity modeling, fractional scalar turbulence LES modeling, in addition to tempered fractional LES SGS modeling for turbulence.
\subsubsection{Nonlocal RANS models: a narrative survey}
\MG{Recall} that most of turbulence models are built based on Boussinesq's turbulent viscosity concept. Thus, one conventionally assumes that the turbulent stress tensors $\tau _{ij} ^{R} $ are proportional to the symmetric part of local mean velocity gradient at any point (i.e., strain-rate tensor). Hence, the corresponding proportionality coefficient, known as the turbulent viscosity, emerges as the unknown turbulence model parameter $\mu_{T}$ in 
\begin{equation}
\label{eq.0}
\tau _{ij} ^{R} =  \mu_{T} \left( \frac{\partial \bar u_i}{\partial x_j} + \frac{\partial \bar u_j}{\partial x_i}  \right).
\end{equation}
 \citet{prandtl1942bemerkungen} in 1942 aimed to move beyond this local constraint by introducing the \textit{extended} mixing length concept for the first time. The corresponding new model was a great migration from locality to nonlocality, but did not achieve a remarkable success as it did not significantly improve accuracy. Afterward, he parametrized the primitive model in a way that the mixing length was taken to be greater than the (differential) length-scale of the problem, including a higher (second-order) Taylor expansion term. This strategy was analogous to adding a ``weak sense of (short-range) nonlocality'' to the model. This was regarded as a \textit{weak nonlocal model} in the sense that the stress term  was still in the form of Boussinesq's and the relation with the strain rate tensor in the same point was collinear. However, von Karman insisted on the consideration of the common local mixing length, which is generated by the local flow conditions and suggested considering the mixing length in terms of two succeeding derivatives \cite{Egolf_book}.  \citet{Bradshow}  in 1973 showed that Boussinesq's hypothesis fails over curved surfaces  and noted that form of the stress-strain relations was the main cause of this failure. It should be mentioned that there were some important developments mostly based on polynomial series, compared to the Boussinesq-type modeling including the works done by Spencer and Rivlin \cite{Spencer, Spencer2}, \citet{Lumley} and \citet{Pope}; however, a noticeable lack accuracy both in terms of physics and mathematics emerged as additional second- (and higher) order tensor series developments were demanded, where an interplay between predictability and practicality remained an open question.  

As indicated in Section \ref{sec:classification}, Brownian motion can serve as a statistical model for the spread of a cloud of particles the continuum limit of which is a parabolic integer-order diffusion equation. Generalizing this approach to heavy-tailed processes such as L\'evy processes can model the intermittency in turbulent flow signals and through a heavy-tailed central limit theorem converge to an anomalous diffusion equation with fractional derivatives in space and/or time \cite{West_5}. This suggests that employing fractional-based Reynolds stresses would be a proper model for the turbulent diffusion term. In a pioneering work by \citet{Hinze_1974} in 1974, the authors described the memory effect in a turbulent boundary layer flow. They utilized the experimental data produced downstream of a hemispherical cap, attached to the lower wall of channel geometry. They demonstrated that when one computes eddy-viscosity using Boussinesq's theory in the lateral gradient of the mean flow and turbulence shear-stresses, there is a huge non-uniform distribution that exists in the outer region of the boundary layer. Interestingly, we see a nonlocal expression for the gradient of the transported field in a novel approach by \citet{Kraichnan} in the same year (1974), for the \textit{scalar quantity transport}. Afterward, fractional-order models based on the RANS approach were offered in \cite{Chen, Cushman, Epps2, Hamba1, Hamba2}. Most of these works are using Green's functions based on the residual velocity to provide the expression for the Reynolds stress or scalar fluxes. 

One of the main contributions for the development of nonlocal models has been done by Egolf and Hutter \cite{Egolf1, Egolf_book}. They started from L\'evy flight statistics and generalized the zero-equation local Reynolds shear stress expression to a nonlocal and fractional type. The method is based on Kraichnanian convolution-integral approach and utilizing different weighting functions. Using the mentioned weighting functions, one can make a bridge between the first-order gradient of the common eddy diffusivity models and the mean velocity difference term. The proposed model is based on the four distinct steps that can be followed conveniently to replace a local operator with a nonlocal one. In reality, the final proposed model is a more general and extended version of Prandtl's zero-equation mixing length and shear-layer turbulence models. The proposed model is called \textit{Difference-Quotient Turbulence Model} (DQTM), given by
\begin{equation}
\overline {u'_2 u'_1} = - \sigma \chi _2[\overline {u}_1 (x_1,x_2) - \overline {u}_{1,\text{min}} (x_1)] \frac{\overline {u}_{1,\text{max}} (x_1)- \overline {u}_1 (x_1,x_2)}{x_{2,\text{max}} - x_2} .
\end{equation}
Although well-motivated and presented as somewhat of a generalization of  classical models, such models have not been thoroughly tested against established integer-order models, and their practical efficiency in addressing the nonlocalities have not been adequately examined. 
Recently, a series of remarkable developments in fractional turbulence modeling gave a new and practical perspective on employing fractional calculus in turbulence, and introducing new statistical measures that directly reveal where classical approaches have room for modernization and enhancement. In what follows we review such cutting edge approaches. 

%
\subsubsection{Fractional eddy-viscosity models}

Recently, \citet{Karniadakis_2021} developed a new nonlocal eddy viscosity-based model (see equation \eqref{eq.1} below) that can be applied in both isotropic and anisotropic turbulent flows. They obtained a proper two-point stress-strain rate correlation structure for \textit{a priori} testing the developed model and performed tests based on the high-resolution DNS data set for the homogeneous isotropic turbulence (HIT) and the channel flow canonical test cases.

The investigation of the model performance is set based on the necessary conditions for any LES approach in providing the accurate two-point statistics of the filtered quantities in the terms of correlations and spectra. The proposed model is given by:
\begin{equation}
\label{eq.1}
\tau _{ij} ^{\alpha} =  - \nu_{T} \big( D_i^{\alpha} \bar u_j + D_j^{\alpha} \bar u_i  \big).
\end{equation}
where the derivative operators $D_i^{\alpha}$ and $D_j^{\alpha}$ are both of order $0<\alpha<1$ respectively in $i$ and $j$ directions, however, they are employed as the \textit{truncated} Caputo derivative variations, still being the convolution of the first derivative of velocity with respect to an inverse power-law kernel with index $\alpha$, however, over a truncated (compact) integral support, forming a finite nonlocality horizon, for the purpose of lowering the computational cost. 

\begin{figure}[!]
	\centering
	\includegraphics[scale=0.3]{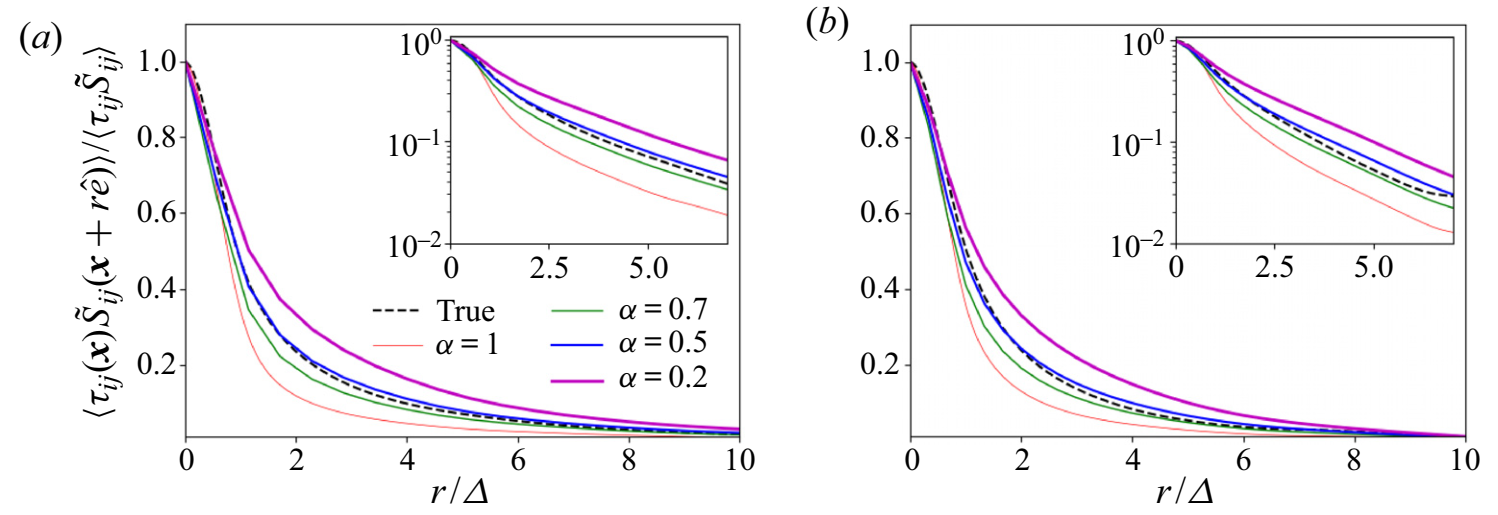}
	\caption{Two-point correlations between SGS stress and filtered rate of strain rate in different scenarios ($\alpha =1$ corresponds to the local model) at the filter size $\Delta=31 \eta$ (a) ,  and $\Delta=53 \eta$ (b). Source: \cite{Karniadakis_2021}.}
	\label{fig:1}
\end{figure}

Several numerical tests conducted in \cite{Karniadakis_2021} indicated that the new model provides a better correlation between the filtered rate of strain rate and subgrid scale stress tensor. Specifically, this model predicts the long tails in the ground-truth subfilter stress–strain-rate correlation functions. However, other conventional local eddy viscosity-based models like classical Smagorinsky, that corresponds to $\alpha = 1$, miss this important feature as they decay faster (see Figure \ref{fig:1}).

In addition to the significantly better capability in the prediction of the long-tail interactions in the new model, the probability density functions of the dissipation quantities for the HIT flow using the box filtering approach, are matching much better than the local model. The local model predictions are purely dissipative and with no tail behavior, which is contrast with the ground-truth DNS data sets. Moreover, effects of the different parameters in the LES procedure have been analyzed including filter size, filter type, wall distance for the channel flow case, and integration radius.    

Alternatively in studying the turbulent transport and mixing, kinetic Boltzmann theory has shown a rich and promising ground based upon principles of statistical mechanics, which by construction is well-suited for the stochastic description of turbulence at microscopic level \cite{harris2004introduction}. In the following, the fundamental sources of nonlocal closure and the SGS modeling for the residual passive scalar flux are studied at the kinetic Boltzmann transport framework. Our objective is to derive a nonlocal eddy-diffusivity SGS model at the continuum level. In what follows we present three recent development of LES SGS where the SGS small motions are modeled by the BGK kinetics transport.

\subsubsection{Fractional LES SGS modeling for scalar turbulence}\label{sec: nonlocality}

Statistical description of LES is well-represented through incorporating a filtering procedure into the kinetic Boltzmann transport. For the purpose of passive scalar transport, applying a spatially and temporally invariant filtering kernel, $\boldsymbol{\mathcal{G}} = \boldsymbol{\mathcal{G}}(\boldsymbol{r})$, onto the distribution  function $g(t,\xb,\boldsymbol{u})$ linearly decomposes that into the filtered, $\widetilde{g}=\boldsymbol{\mathcal{G}} \ast g$, and the residual, $g^{\prime}=g-\widetilde{g}$, components. Therefore, filtering the BGK equation results in the following filtered BTE (FBTE) for the passive scalar:
\begin{equation}\label{eqn: FBTE}
\frac{\partial \widetilde{g}}{\partial t} + \boldsymbol{u}\cdot \nabla \, \widetilde{g} = -\frac{\widetilde{g}-\widetilde{g^{\text{eq}}(\mathcal{B})}}{\tau_g}.
\end{equation}
{where $\mathcal{B}$ represents the generic \textit{Boltzmann  filter size}}.
As elaborated by \citet{girimaji2007boltzmann}, the nonlinear nature of the collision operator, $C_{\mathrm{BGK}}(g)$, prohibits the filtering kernel to commute with; thus, it initiates a source of closure at the kinetic level in FBTE \eqref{eqn: FBTE}. Defining $\widetilde{\mathcal{B}}:=(\boldsymbol{u}-\widetilde{\boldsymbol{V}})^2/c_T^2$, this closure problem is manifested in the following inequality,
\begin{equation}\label{eqn: kinetic-closure}
\widetilde{g^{\text{eq}}(\mathcal{B})} = \frac{\widetilde{\Phi \, \exp(-\mathcal{B}/2)}}{(2\pi)^{3/2} \, c_T^3} \neq \frac{\widetilde{\Phi} \, \exp(-\widetilde{\mathcal{B}}/2)}{ (2\pi)^{3/2} \, c_T^3} = g^{\text{eq}}(\widetilde{\mathcal{B}}).
\end{equation}
The identified closure requires proper means of modeling so that one can numerically solve the FBTE \eqref{eqn: FBTE}. A common practice is to approximate this closure problem with a modified relaxation time approach that is described in detail in \cite{sagaut2010toward}. Despite the success of this approach in some applications, it is not physically consistent with the filtered turbulent transport dynamics \cite{girimaji2007boltzmann}. Nevertheless, here we manage to adjust this inconsistency by looking at the nonlocal effects arising from filtering the Maxwell distribution function, $g^{\text{eq}}(\mathcal{B})$, and model them with proper mathematical tools. Considering the spatial filtering kernel $\boldsymbol{\mathcal{G}}(\boldsymbol{r})$ with the filter-width $\Delta$, and applying it on the Maxwell equilibrium distribution as
\begin{equation}\label{eqn: F-Maxwell}
\widetilde{g^{\text{eq}}(\mathcal{B})} =  \boldsymbol{\mathcal{G}} \ast g^{\text{eq}}\big (\mathcal{B}(t,\boldsymbol{u},\xb)\big ) =  \int_{R_f}^{} \boldsymbol{\mathcal{G}}(\boldsymbol{r}) \,  g^{\text{eq}}\big (\mathcal{B}(t,\boldsymbol{u},\xb-\boldsymbol{r})\big ) \, d\boldsymbol{r},
\end{equation}
where $R_f=[-\Delta/2 \, , \Delta/2]^3$. Subsequently, by rewriting the right-hand side of the passive scalar FBTE \eqref{eqn: FBTE} into the following form
\begin{align}\label{eqn: FBTE-RHS}
-\frac{1}{\tau_g} \left(\widetilde{g} - \widetilde{g^{\text{eq}}(\mathcal{B})} \right) = \underbrace{-\frac{1}{\tau_g} \left(\widetilde{g} - g^{\text{eq}}(\widetilde{\mathcal{B}}) \right)}_{\text{closed}} + \underbrace{\frac{1}{\tau_g} \left(\widetilde{g^{\text{eq}}(\mathcal{B})} - g^{\text{eq}}(\widetilde{\mathcal{B}}) \right)}_{\text{unclosed}},
\end{align}
the unclosed part is structurally multi-exponentially distributed and maybe approximated by a power-law distribution model as we propose 
\begin{align}\label{eqn: Levy-model}
\widetilde{g^{\text{eq}}(\mathcal{B})} - g^{\text{eq}}(\widetilde{\mathcal{B}}) \approx g^\alpha(\widetilde{\mathcal{B}}) = \frac{\widetilde{\Phi}}{c_T^3} \, F^{\alpha}(\widetilde{\mathcal{B}}),
\end{align}
where $F^{\alpha}(\widetilde{\mathcal{B}})$ denotes an $\alpha$-stable L\'evy distribution that is mathematically designed based on heavy-tailed stochastic processes and replicate the power-law behavior \cite{applebaum2009levy, meerschaert2019stochastic}. The corresponding macroscopic continuum variables associated with the filtered \eqref{GE-1-3} are obtained in terms of the filtered distribution functions, $\widetilde{f}$ and $\widetilde{g}$, as
\begin{eqnarray}\label{eqn: continuum-ave-fPhi}
\widetilde{\Phi} &=& \int_{\mathbb{R}^d} \widetilde{g}(t,\xb,\boldsymbol{u}) \,  d\boldsymbol{u},
\\ \label{eqn: continuum-ave-fVPhi}
\widetilde{V}_i &=& \frac{1}{\rho}\int_{\mathbb{R}^d} u_i \, \widetilde{f}(t,\xb,\boldsymbol{u}) \, d\boldsymbol{u}, \quad i=1,2,3.
\end{eqnarray}
According to the microscopic reversibility of the particles that assumes the collisions occur \textit{elastically}, the right-hand side of \eqref{eqn: FBTE} equals zero \cite{saint2009hydrodynamic}. Therefore,
\begin{align}\label{GE-14}
\frac{\partial \widetilde{\Phi}}{\partial t} + \nabla\cdot \int_{\mathbb{R}^d} \boldsymbol{u} \, \widetilde{g} \, d\boldsymbol{u} = 0. 
\end{align}
Since we are working with spatial filtering kernels, $\boldsymbol{\mathcal{G}}=\boldsymbol{\mathcal{G}}(\boldsymbol{r})$,
\begin{align}\label{GE-16}
\int_{\mathbb{R}^d} \boldsymbol{u} \, \widetilde{g} \, d\boldsymbol{u} = \int_{\mathbb{R}^d} (\boldsymbol{u}-\widetilde{\boldsymbol{V}}) \, \widetilde{g} \, d\boldsymbol{u}+ \int_{\mathbb{R}^d} \widetilde{\boldsymbol{V}} \, \widetilde{g}\, d\boldsymbol{u}.
\end{align}
By plugging \eqref{GE-16} into \eqref{GE-14}, we obtain that
\begin{equation}\label{GE-17}
\frac{\partial \widetilde{\Phi}}{\partial t} + \nabla \cdot \left(\widetilde{\Phi} \, \widetilde{\boldsymbol{V}}\right) = -\nabla \cdot \boldsymbol{q},
\end{equation}
where 
\begin{equation}\label{GE_17_2}
q_i=\int_{\mathbb{R}^d} \left(u_i-\widetilde{V}_i\right) \, \widetilde{g} \, d\boldsymbol{u}.
\end{equation} 
The corresponding filtered passive scalar flux is obtained through a sequence of step-by-step derivations as
\begin{align}\label{eqn: flt-flux}
\widetilde{\boldsymbol{q}} = -\mathcal{D} \, \nabla \widetilde{\Phi},
\end{align}
and the divergence of residual scalar flux is derived as the fractional Laplacian of the filtered total scalar concentration,
\begin{align}\label{eqn: res-flux}
\nabla \cdot \boldsymbol{q}^R = -\mathcal{D}_\alpha \, (-\Delta)^{\alpha} \, \widetilde{\Phi}, \quad \alpha \in (0,1],
\end{align}
where $\mathcal{D}_\alpha := \frac{C_\alpha (c_T \, \tau_g)^\alpha}{\tau_g} \, (\alpha+2) \, \Gamma(\alpha)$ is a model coefficient with the unit [$L^\alpha/T$]. The filtered AD equation for the total passive scalar concentration, developed from the filtered kinetic BTE with an $\alpha$-stable L\'evy distribution model, yields a fractional-order SGS scalar flux model at the continuum level. The aforementioned filtered AD equation reads as
\begin{align} \label{eqn: Flt-AD-total}
\frac{\partial \widetilde{\Phi}}{\partial t}+\frac{\partial}{\partial x_i}\left( \widetilde{\Phi} \, \widetilde{V}_i\right) = \mathcal{D} \, \Delta \widetilde{\Phi} +\mathcal{D}_{\alpha} (-\Delta)^{\alpha} \, \widetilde{\Phi}.
\end{align}
Through a proper choice for the fractional Laplacian order $\alpha$, the developed model optimally works in an LES setting. Applying the Reynolds decomposition and considering the passive scalar with imposed uniform mean gradient, equation \eqref{eqn: Flt-AD-total} fully recovers the filtered transport equation for the transport of the filtered scalar fluctuations, $\widetilde{\phi}$.

\subsubsection{Fractional/tempered-fractional LES models for fluid turbulence}
For some pedagogical purposes, we first presented the case of fracitonal LES SGS modeling for scalar turbulence. However, this new paradigm in LES modeling actually began prior to \cite{Scalar_FSGS}. \citet{FSGS}  developed the first ever fractional LES model for homogeneous isotropic turbulent flows as
\begin{equation}
\label{eq.2}
(\nabla . \tau ^{R}) = \mu_{\alpha}  (- \Delta )^{(\alpha)} \bar V,	
\end{equation}
being based upon on the derivation of fractional Laplacian closure term in the spatially filtered Navier-Stokes equations when employing a L\'evy stable distribution as the equilibrium model in the filtered BGK kinetic transport equation. In \cite{Tempered_FSGS}, they later developed a generalized version of this earlier model (suitable for incorporating data with tailored/truncated tails). Employing rather a \textit{tempered L\'evy stable} distribution in the kinetic level this time gave rise to the formulation of the tempered fractional LES closure term as
\begin{equation}
\label{eq.3}
(\nabla . \tau ^{R}) = \nu_{\alpha} \sum _{k=0}^{\kappa} \phi_k ^{\kappa} ( \Delta + \lambda _k)^{(\alpha)} \bar V,	
\end{equation}
forming a novel, data-friendly and expressive tempered fractional Laplacian SGS model for turbulence. They also showed that the newly developed nonlocal models can better recover the non-Gaussian statistics of subgrid-scale stress motions while they are being employed at the continuum level.

\subsection{Future directions in fractional turbulence modeling}\label{sec:turbulence-future}

\citet{Laval}  analyzed the effects of the local and nonlocal interactions on the intermittency corrections in the scaling properties in three-dimensional turbulence. They observed that nonlocal interactions are responsible for the creation of the intense vortices and on the other hand, local interactions are trying to dissipate them. Inspired by the mentioned observations, they came up with a new turbulence model that accounts for both the local and nonlocal interactions for the study of intermittency. In their proposed model, the large and small scales are being coupled by nonlocal interactions using a multiplicative process and additive noise along with a turbulent viscosity model for the local interactions. The results of the new model qualitatively cover the previously observed anomaly and intermittency aspects. 

In the context of nonlocal turbulence modeling, \citet{song2018universal} proposed a variable-order fractional model for wall-bounded turbulent (mean) flows. They represented the Reynolds stresses with a nonlocal fractional derivative of variable-order that decays with the distance from the wall. Interestingly, they found that this variable fractional order has a \textit{universal form} for all {\it Re} and for three different flow types, i.e., channel flow, Couette flow, and pipe flow. In addition to the aforementioned fully-developed flows, they modelled turbulent boundary layers and discussed how the streamwise variation affects the universal curve (see also \citep{song2021variable} for a follow-up work). Later, \citet{pang2020npinns} proposed a nonlocal truncated operator with spatially variable order, which is suitable for modeling wall-bounded turbulence, e.g. turbulent Couette flow. They showed that nonlocal physics-informed neural networks (nPINNs) can jointly infer the variable order, exhibiting a universal behavior with respect to \textit{Re}, a finding that can contribute to better understanding of nonlocal interactions in wall-bounded turbulence.
In terms of memory-effects (i.e., nonlocality in time), Parish and Duraisamy in \cite{parish2017dynamic} developed a dynamic SGS model for LES, based on the Mori--Zwanzig (MZ) formalism. This closure model was constructed by exploiting similarities between two levels of coarse-graining via the Germano identity of fluid mechanics and by assuming that memory effects have a finite-temporal support. This work suggests future studies on using time-fractional derivatives in turbulence models.

The aforementioned developments are practically interesting, mathematically exciting, and algorithmically robust. They 
encourage the field of research in turbulence to \MG{open its} arms towards \MG{the} new wealth of mathematical developments in both theory and practice of fractional modeling. Inevitably, further systematic studies and developments of nonlocal turbulence models are needed (both numerically and experimentally) in order to achieve the charming blend of enhanced accuracy and lowered cost in realistic applications. On this note, we end this section by emphasizing that the non-Markovian/non-Fickian nature does not relax in compressible flows (i.e., variable density problems). Therefore the idea of developing generalized fractional turbulence models for transonic-to-hypersonic flows is a new and \MG{ripe} venue for research, in which the existing sense of classical thermodynamics can become fundamentally non-equilibrium and nonlocal.

\section{Fractional constitutive laws in material science} \label{sec:material}

Accurate modeling of evolving material response and failure across multiple time- and length-scales is essential for life-cycle prediction and design of new materials. While the mechanical behavior of a number of standard engineering materials (e.g., metals, polymers, rubbers) is quite well-understood, a significant modeling effort still needs to be conducted for complex materials, where microstructure heterogeneities, randomness and small scale physical mechanisms (such as collective behavior) lead to non-standard and, at times, counter-intuitive responses. Two examples are bio-tissues and natural materials (e.g. biopolymers), which are multi-functional products of millions of years of evolution, locally optimized for their hosts and environment, and constrained by a limited set of building blocks and available resources \cite{Imbeni2005,Wegst20151}. These materials possess unprecedented properties at low densities, especially due to their hierarchical and multi-scale structure, leading to a wide spectrum of behaviors, such as power-law viscoelasticity, viscoplastic strains under hysteresis loading, damage, failure, fatigue, fractal avalanche ruptures and self-healing mechanisms.

The main motivation for fractional materials modeling is the power-law fingerprint arising in microstructures undergoing anomalous diffusion, observed in a range of complex materials. Such microstructures often display a fractal nature with sub-diffusive dynamics, e.g., of entangled polymer chains, and defect interactions such as dislocation avalanches, cracks and voids. Such non-exponential behavior cannot be accurately modeled by integer-order, linear viscoelastic models, which require arbitrary arrangements of Hookean/Newtonian elements and introduce a limited number of exponential (Debye) relaxation modes that, at most, represent a truncated power-law approximation \cite{bagley1989power}. While these approximations may be satisfactory for short times and engineering precision, they often result in high-dimensional parameter spaces and still lack predictability outside the experimental time/length-scales, often requiring recalibration. In this context, fractional operators become appropriate and natural modeling choices, since their integro-differential operators naturally utilize power-law convolution kernels, coding self-similar microstructural features in a reduced-order mathematical language with smaller parameter spaces (similarly to the case of anomalous transport, see Section \ref{sec:subsurface}). This fact allows accurate and predictive modeling, in an efficient manner, of bio-tissues \cite{Magin2010,Magin2010a,Craiem2008,Suki1994,Djordjevic2003,Nicolle2010,Puig2007} and polymers \cite{Winter1997,Ketz1988,Baghdadi2005,Sollich1998} for multiple time-scales.

In this \MG{s}ection we review fractional models for materials undergoing power-law behaviors, termed \textit{anomalous materials}, in a range of non-equilibrium and path-dependent responses. We start with linear viscoelasticity, introducing the basic modeling building block, known as \textit{Scott-Blair} element that models a single power-law response and can be combined to incorporate more complex behaviors. In harmony with the previous \MG{s}ections, we will emphasize on potential multi-scale connections, stochastic processes, and thermodynamic consistency. After providing evidence of cases where fractional behavior/power-laws appear as intrinsic qualities in a number of systems, we report on the state-of-the-art models incorporating multiple physical mechanisms.

\paragraph{Fractional viscoelasticity: rheological building blocks.} We start with the Boltzmann superposition integral for linear viscoelasticity, obtained from the linear superposition of infinitesimal step strains $\delta \varepsilon(t)$ applied to a viscoelastic material \cite{Mainardi2011}:
\begin{equation}\label{eq:Boltzmann-integral}
    \sigma(t) = \int^t_{-\infty} G(t-\tau) \dot{\varepsilon}(\tau)\,d\tau,
\end{equation}
where $\dot{\varepsilon}$ and $\sigma(t)$ denote, respectively, the strain rate and stress. The convolution kernel $G(t)$, is a relaxation function, directly related to stress relaxation experiments under step strains. It is traditionally modeled through combinations of Hookean springs and Newtonian dashpots, yielding a multi-exponential relaxation in the form $G(t) = \sum^N_{i=1} C_i \exp(-t/\tau_i)$. In this particular choice of kernel, \eqref{eq:Boltzmann-integral} is equivalent to a multi-term ordinary differential equation (ODE). 

Relaxation experiments across multiple time- and frequency-scales indicate that anomalous materials exhibit memory effects in time for stress/strain responses, which translates into a single power-law scaling in the form $G(t) \propto t^{-\alpha}$, with $\alpha \in (0,1)$. This indicates that, contrary to exponential relaxation forms, there is a spectrum of relaxation times arising from the material microstructure \cite{Jaishankar2013}, for which standard ODE models (e.g. generalized Maxwell model in creep/relaxation representations) would require a large number of parameters.

The fundamental fractional rheological building block element, termed \textit{Scott-Blair} (SB) model, is obtained by substituting the power-law kernel $G(t) = Et^{-\alpha}/\Gamma(1-\alpha)$ into \eqref{eq:Boltzmann-integral}, leading to the following form:
\begin{equation}\label{eq:SB-semiinfinite}
    \sigma(t) = \prescript{\text{C}}{-\infty}{}\mathbb{D}^\alpha_t \varepsilon(t) =  \frac{E}{\Gamma(1-\alpha)}\int^t_{-\infty} (t-\tau)^{-\alpha} \dot{\varepsilon}(\tau)\,d\tau,
\end{equation}
which is equivalent to the Riemann-Liouville fractional derivative $\prescript{\text{RL}}{-\infty}{}\mathbb{D}^\alpha_t \varepsilon(t)$ if the function $\varepsilon(t)$ is sufficiently well behaved at $t \to -\infty$ \cite{podlubny1998fractional}. While this equivalence is satisfied for semi-infinite domains, the choice of Riemann-Liouville and Caputo definitions matter when we introduce a causal strain history and switch the lower limit of \eqref{eq:SB-semiinfinite} from $-\infty$ to $0$, which leads to two different fractional Cauchy problems. For the Caputo definition, we have \cite{Mainardi2011}:
\begin{equation}\label{eq:SB-caputo}
  \sigma(t) = E \prescript{\text{C}}{0}{}\mathbb{D}^\alpha_t \varepsilon(t),\quad t > 0, \quad 0 < \alpha < 1,\quad \varepsilon(0) = \varepsilon_0.
\end{equation}
On the other hand, when employing Riemann-Liouville derivatives, we obtain:
\begin{equation}\label{eq:SB-RL}
  \sigma(t) = E \prescript{\text{RL}}{0}{}\mathbb{D}^\alpha_t \varepsilon(t), \quad t > 0, \quad 0 < \alpha < 1,\quad \prescript{\text{RL}}{0}{}\mathbb{D}^{\alpha-1}_t \varepsilon(t)\big|_{t=0} = \varepsilon_0,
\end{equation}
where we remark that problem \eqref{eq:SB-caputo} is more commonly adopted due to the appearance of integer-order ICs, while both aforementioned problems are equivalent in the presence of homogeneous ICs. The SB element provides a constitutive interpolation between a Hookean spring ($\alpha \to 0$) and a Newtonian dashpot ($\alpha \to 1$). The unique parameter pair $(E\,[Pa.s^\alpha],\alpha)$ codes snapshots of a dynamic process instead of an equilibrium state of the system \cite{Jaishankar2013}. Consequently these properties are only associated with equilibrium states in the limit cases for the fractional order $\alpha$. We remark that although the FDE \eqref{eq:SB-caputo} utilizing the Caputo definition is widely employed to represent the SB element in the literature, the pioneering works on anomalous rheology modeling are attributed to \citet{gerasimov1948} in 1948, introducing a similar power-law convolution operator as \eqref{eq:SB-semiinfinite}, which may be referred in the literature as the Gerasimov-Caputo derivative \cite{valerio2014some}. We refer the reader to \cite{valerio2014some,mainardi2012historical} for more details on the historical context of fractional derivatives in visco-elasticity.

\paragraph{Mechanistic and thermodynamic interpretations.} Apart from the Boltzmann integral representation \eqref{eq:Boltzmann-integral}, characterized by an integro-differential nature, the SB element can also be obtained through a continuous arrangement of canonical, Hookean and Newtonian elements, both from their constitutive and free-energy levels \cite{Schiessel1993,Lion1997}, making the notion of SB elements intrinsically incorporating an infinite number of relaxation times more evident. In \cite{Schiessel1993}, a hierarchical ladder-like structure of standard Maxwell viscoelastic elements was employed. This structure led to a coupled system of ODEs, which had an infinite continued fraction (a recursion of fractions) representation in terms of the Maxwell model constants in the Laplace domain. Then, applying an inverse Laplace transform, a fractional stress-strain relationship was recovered for homonegeous initial conditions, therefore equivalent to both forms \eqref{eq:SB-caputo} and \eqref{eq:SB-RL}. In \cite{Lion1997}, an isothermal Helmholtz free-energy density was derived for the SB element from the elastic energies of a discrete-to-continuum arrangement of standard Maxwell branches, obtaining the following form for the free-energy $\psi$ as a function of the strain:
\begin{equation}\label{eq:SB-energy}
    \psi(\varepsilon) = \frac{1}{2} \int^\infty_0 \tilde{E}(z)\left[\int^t_0 \exp(-\frac{t-s}{z})\dot{\varepsilon}(s) ds\right]^2dz, \,\, \tilde{E}(z) = \frac{E z^{-1-\alpha}}{\Gamma(\alpha)\Gamma(1-\alpha)},
\end{equation}
where $\tilde{E}$ denotes the relaxation spectrum. Therefore, \eqref{eq:SB-energy} represents the amount of available elastic energy to perform work from the SB element in the time domain, which cannot be directly inferred from \eqref{eq:SB-caputo} and \eqref{eq:SB-RL}. Naturally, the two limit cases for $\alpha$ are $\psi(\varepsilon) \to E\varepsilon^2/2$ when $\alpha \to 0$, and $\psi(\varepsilon) \to 0$ when $\alpha \to 1$. Furthermore, under suitable thermodynamic constraints, it is shown that the SB element is thermodynamically admissible and that the Caputo representation of \eqref{eq:SB-RL} can be derived from \eqref{eq:SB-energy} under continuum mechanics arguments.

\paragraph{Energy decoupling in the frequency domain.} Similar to the aforementioned representations, power-law structures also appear in viscoelastic dynamic properties and rheological experiments in the frequency domain \cite{Jaishankar2013}, such as the complex shear modulus, defined as the ratio between the Fourier transform of stresses and strains:
\begin{equation}\label{eq:complex-modulus}
    G^*(\omega) := \frac{\mathcal{F}[\sigma](\omega)}{\mathcal{F}[\varepsilon](\omega)} = G^\prime(\omega) + i G^{\prime\prime}(\omega),
\end{equation}
where $\omega\,[s^{-1}]$ denotes the frequency. The term $G^\prime$ is the storage modulus, and $G^{\prime\prime}$ denotes the loss modulus, i.e., the stored and dissipated energy per cycle, respectively. Employing definition \eqref{eq:complex-modulus} into \eqref{eq:SB-RL}, the dynamic modulus of the Scott-Blair element is obtained \cite{Schiessel1995}:
\begin{equation}
    G^\prime(\omega) = \MG{\text{Re}} (G^*) = E\omega^\alpha \cos\left(\frac{\alpha \pi}{2}\right), \quad G^{\prime \prime}(\omega) = \MG{\text{Im}} (G^*) = E\omega^\alpha \sin\left(\frac{\alpha \pi}{2}\right),
\end{equation}
which provides a clear storage/loss decomposition, with the value of $\alpha$ determining whether the material of interest is predominantly dissipative for certain frequency ranges.

\paragraph{Relationships to material microstructure and stochastic processes.} The mechanistic origins of macroscopic power-law behaviors in complex materials are due to spatio-temporal anomalous sub-diffusive processes \cite{Metzler2000} in fractal micro-structures. We focus on the temporal case, in which the MSD of microstructural constituents follows a nonlinear scaling in the form $\langle \Delta x \rangle^2 \propto t^\alpha$. 
\citet{Bagley1983} provided a relationship between the complex shear modulus obtained from the Rouse theory of polymer dynamics. They started with the result of Rouse's theory for the shear modulus, i.e.
\begin{equation}
    G^\prime(\omega) = n k T \sum^N_{p=1} \frac{\omega^2 \tau^2_p}{1+\omega^2\tau^2_p}, \quad
    G^{\prime\prime}(\omega) = \omega \mu_s + n k T \sum^N_{p=1} \frac{\omega \tau_p}{1+\omega^2\tau^2_p},
\end{equation}
where $n$ denotes the number of molecules per unit volume, $N$ is the number of monomers in the polymer chain, $T$ represents the absolute temperature, $k$ is Boltzmann's constant. The term $\tau_p$ denotes the relaxation times of the solution, which was approximated as $\tau_p \approx \tau_1/p^2 = 6(\mu_0 - \mu_s)/(p^2\pi^2nkT)$, which is valid when the number of submolecules $N$ is large. The terms $\mu_0$ and $\mu_s$ denote, respectively, the steady-flow viscosities of the solution and solvent. They further worked on Rouse's results, and by assuming the polymer chains and $\omega \tau_1$ to be sufficiently large, obtained the following power-law form for the dynamic shear modulus:
\begin{equation}
    G^*(\omega) = i \omega \mu_s + \left[\frac{3}{2} (\mu_0 - \mu_s) nkT\right]^{1/2}(i\omega)^{1/2}.
\end{equation}
After applying the inverse Fourier transform, the above relationship leads to a Riemann-Liouville representation between stresses-strains with $\alpha = 1/2$. Similar observations were also reported for $\sigma(t)$ utilizing a Zimm model, where the inclusion of hydrodynamic interactions lead to a fractional order $\alpha = 2/3$.
\citet{glockle1993fox} showed that fractional relaxation can be modeled by a special type of CTRW describing a trapping problem due to entanglements of polymer chains, thus slowing down the relaxation process. In their work, the random walkers, i.e., the particles, are considered as packages of free volume that allow conformational reorientations of chain segments, thus leading to relaxation. They obtained a waiting time distribution of such particles through a Fox-Wright representation in the form:
\begin{equation}
    \chi(t) \sim \frac{{A}}{{\bar{\tau}}} \sum^\infty_{k=0} \frac{(-1)^k}{\Gamma(-\beta k - \beta)} \left(\frac{\bar{\tau}}{t}\right)^{\beta k + \beta + 1},
\end{equation}
for which the leading term indicates that the CTRW waiting time corresponding to fractional relaxation exhibits a Lévy-type decay in the form $\chi(t) \sim t^{-\beta-1}$.

\paragraph{Connecting dynamic viscoelasticity across scales.} A connection between power-laws propagating from micro- to macro-rheology was proposed in \cite{mason2000estimating}, with the use of a Generalized Stokes-Einstein Relation (GSER) for spheres undergoing generalized Langevin dynamics in a viscoelastic medium:
\begin{equation}\label{eq:GSER}
    \lvert G^*(\omega)\rvert \approx \frac{k T}{\pi a \langle \Delta r^2(1/\omega)\rangle \Gamma[1+\alpha(\omega)]}, \quad \alpha(\omega) \equiv \frac{d\, \mathrm{ln} \langle \Delta r^2(1/\omega)\rangle}{d\, \mathrm{ln}\,t}\big|_{t=1/\omega},
\end{equation}
which is valid for spheres of radius $a$ comparable to the length-scale 
of the embedding medium. Here, the dynamic shear modulus $G^*(\omega)$ is related to a velocity memory function from Langevin dynamics. Among a variety of representations for the GSER, \eqref{eq:GSER} assumes a power-law structure of the MSD with exponent $\alpha$, which approaches zero when the sphere is confined by elastic structures present in the complex fluid. Such power-law representation also reduces errors near the frequency extremes when employing Laplace and Fourier transforms.

\paragraph{Physical interpretation of fractional orders.} Despite existing connections between micro- and macro-rheological properties, the physical interpretation of the emerging fractional orders has been elusive. More recently, a connection between the fractional order and the fractal dimension of the material microstructure was made by \citet{Mashayekhi2019Fractal}, where the authors extended the Zimm theory of polymer dynamics to fractal media as a bridge between the meso- and macro-scales. They showed that the fractional order is a rate-dependent material property that is strongly correlated with the fractal and spectral dimensions in fractal media.

\subsection{Evidence of fractional behavior}\label{sec:material-evidence}

We provide a few examples of fractional/power-law behaviors in viscoelasticity and micro/macro-scale plasticity. We start with two examples in viscoelasticity of solid-like and fluid-like natures in which fractional modeling is more appropriate, both with better fits and a reduced number of model parameters.

\begin{figure}[t!]
        \centering
        \begin{subfigure}[b]{0.37\textwidth}
				\includegraphics[width=\columnwidth]{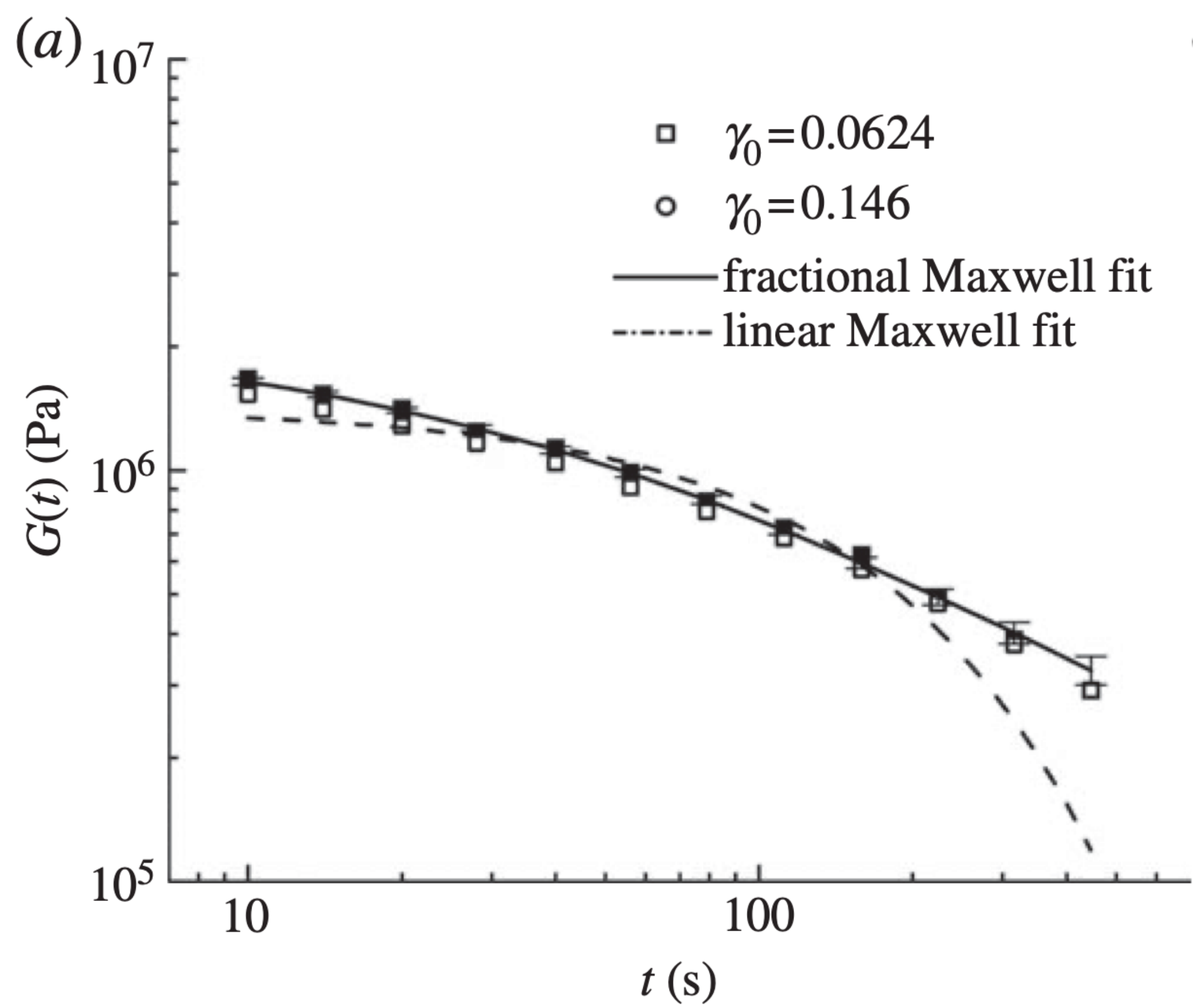}
				\caption{}

        \end{subfigure}%
        ~
        \begin{subfigure}[b]{0.4\textwidth}
                \includegraphics[width=\columnwidth]{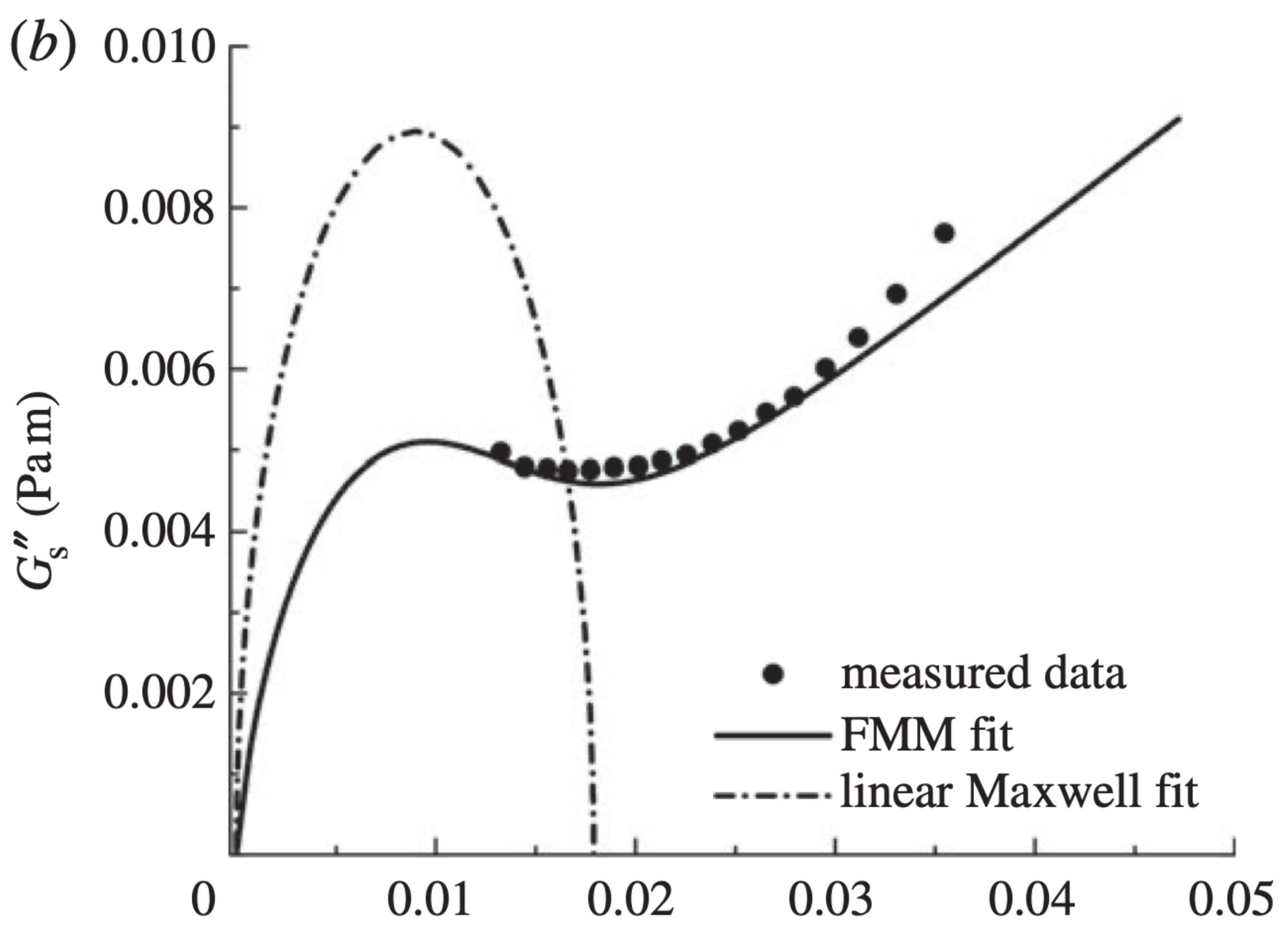}
                \caption{}
        \end{subfigure}
\caption{Comparison between standard and fractional order models. \textit{(a)} Relaxation behavior of Butyl rubber using experimental data from Scott-Blair. \textit{(b)} Cole-Cole plot ($G^\prime$ \textit{versus} $G^{\prime\prime}$) for the dynamic properties of acacia gum. Source: \cite{Jaishankar2013}. \label{fig:mat_fractional_classical}}
\end{figure}

\paragraph{Viscoelastic rheology.} \citet{Jaishankar2013} calibrated classical and fractional Maxwell models to the four orders-of-magnitude relaxation data for highly anomalous butyl rubber data from \citet{blair1947limitations} (Fig.\ref{fig:mat_fractional_classical} (a)), and observed that the three-parameter fractional Maxwell model provided an excellent fit to the experimental data, while a multi-exponential, integer-order Maxwell model required six parameters to provide a satisfactory fit. Moreover, using the calibrated fractional relaxation parameters they obtained an accurate prediction of the creep compliance for the same material, especially for long-time behavior. The second experiment from \cite{Jaishankar2013} concerns the dynamic properties of acacia gum, a commonly used food preservative. In this case, they compared a four-parameter fractional Maxwell model with a single mode (three-parameter) standard Maxwell model (Fig.\ref{fig:mat_fractional_classical}(b)) and demonstrated that while the fractional Maxwell model captures a complex Cole-Cole behavior, its integer-order counterpart is unable to even estimate the qualitative response.
We note that other factors, such as material heterogeneity can introduce multiple power-law relaxation regimes. 

In \cite{Stamenovic2007} Stamenovi\'{c} measured the complex shear modulus $G^*(\omega)$ of cultured human airway smooth muscle and observed two distinct power-law regimes separated by an intermediate plateau. \citet{Kapnistos2008}  found an unexpected tempered power-law relaxation response of entangled polystyrene ring polymers, compared to the usual relaxation plateau of linear chain polymers. Such behavior was interpreted through self-similar conformations of double-folded loops of ring polymers, instead of the reptation observed in linear chains.

\paragraph{Power-law plasticity.} The creep behavior of human embrionic stem cells (ESCs) under differentiation was studied by \citet{Pajerowski2007} through micro-aspiration experiments at different pressures. The cell nucleus demonstrated distinguished visco-elasto-plastic power-law scalings, with $\alpha = 0.2$ for the plastic regime, independent of the applied pressure. It is discussed that such low power-law exponent arises due to the fractal arrangement of chromatin inside the cell nucleus.
\begin{figure}[t]
        \centering
        \begin{subfigure}[b]{0.4\textwidth}
				\includegraphics[width=\columnwidth]{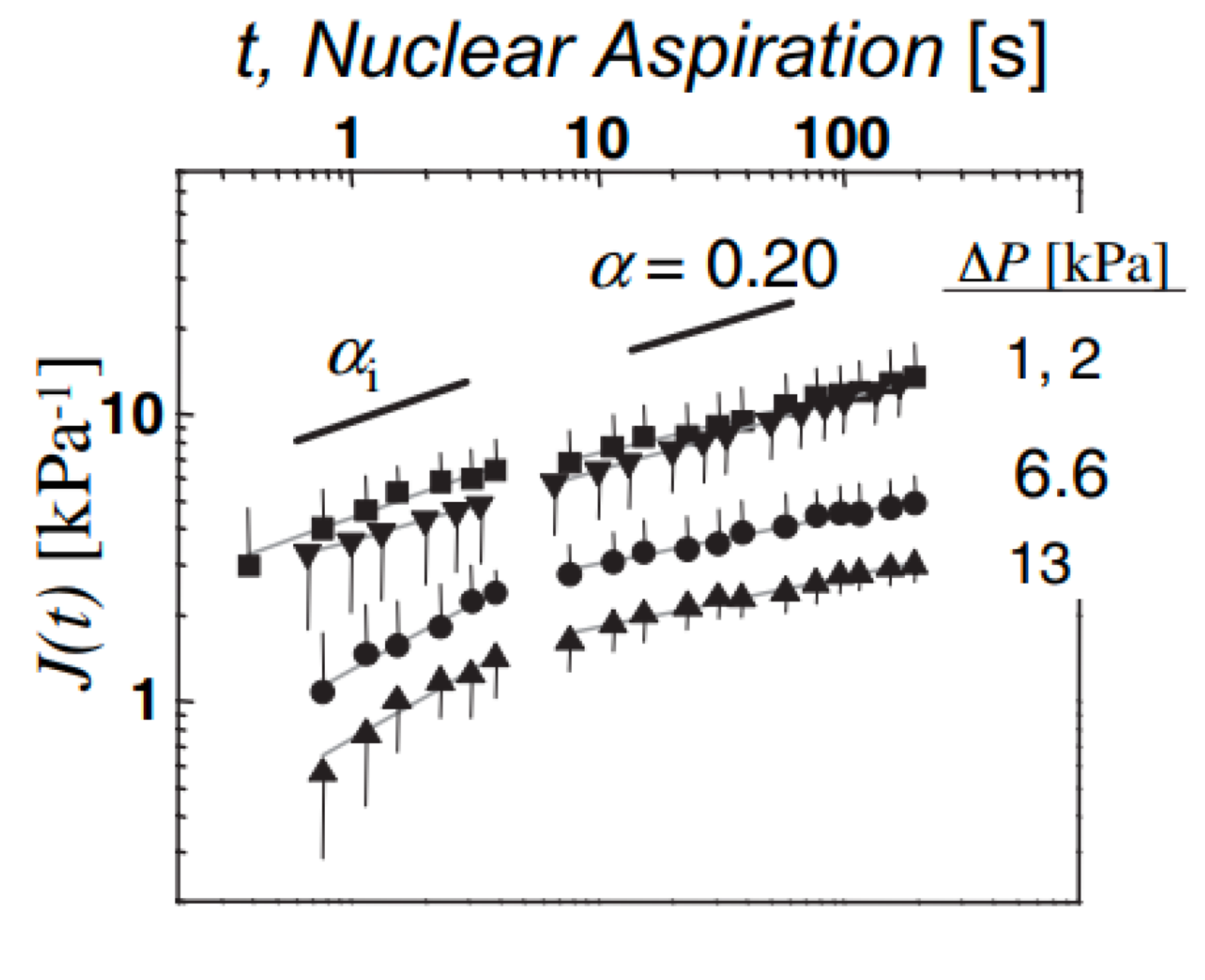}
				\caption{ }

        \end{subfigure}%
        ~
        \begin{subfigure}[b]{0.35\textwidth}
                \includegraphics[width=\columnwidth]{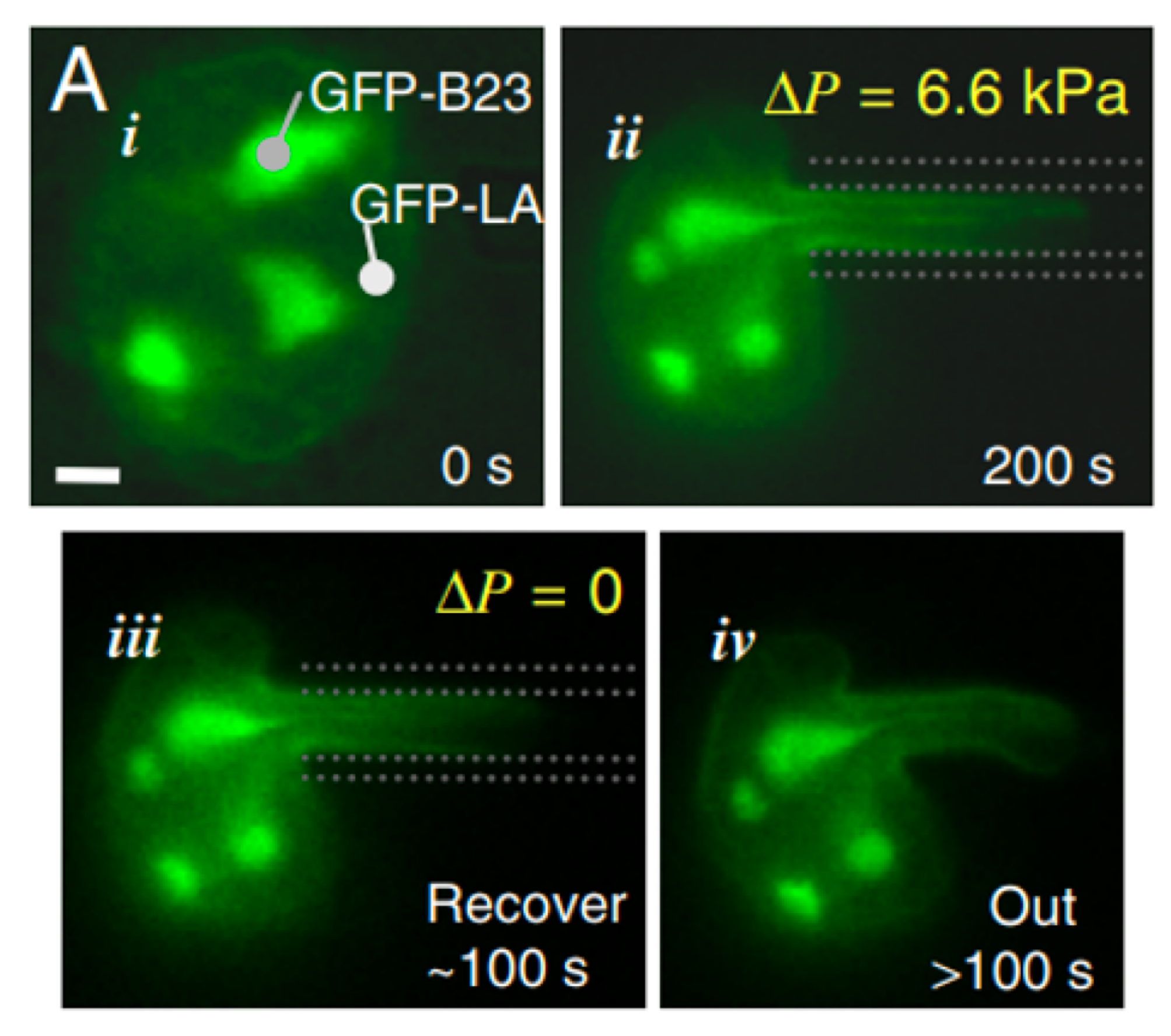}
                \caption{}
        \end{subfigure}
\caption{\textit{(a)} Scale-free creep of ESCs nuclei under aspiration.
Low applied stresses $\Delta P$ yield a single a single power-law creep scaling. For large stresses, a plastic transition is observed at $\tau_{plastic} \approx 8-10\, [s]$, with a creep exponent $\alpha \approx 0.2$, independent of stress values. \textit{(b)} Different stages of nucleus
aspiration, showing a viscoelastic recovery (ii)-(iii), followed by
irreversible plastic deformation (iv). Source: \cite{Pajerowski2007}. \label{fig:Plast_1}}
\end{figure}

Studies on force-induced mechanical plasticity of mouse embrionic fibroblasts were performed by \citet{Bonadkar2016}. It was found that the viscoelastic relaxation and the permanent deformations followed a stochastic, normally-distributed, power-law scaling $\beta(\omega)$, with values ranging from $\beta \approx 0$ to $\beta \approx 0.6$. The microstructural mechanism of plastic deformation in the cytoskeleton is due to the combination of permanent stretching and buckling of actin fibers.

As for evidence of power-laws in failure of crystalline materials, \citet{Richeton2005} investigated the emergence of intermittency and dislocation avalanches in polycrystalline plasticity through acoustic emission experiments on ice under creep compression. Their findings demonstrate that different from the scale-free, close-to-critical dislocation dynamics of single crystals \cite{Miguel2001}, the introduction of average grain sizes $\langle d \rangle$ from the polycrystal microstructure led to a tempered power-law distribution of avalanche sizes. While the exponential tempering cutoff changes with $\langle d \rangle$, the authors observed a constant power-law scaling for all samples.

\paragraph{Connections to stochastic processes.} Although the sub-diffusive MSD coefficient $0 < \alpha < 1$ is observed in a variety of studies for complex materials and fluids, there exist different interpretations on the underlying stochastic processes linked to the sub-diffusive physics, e.g., crowding or caging effects in cells and polymers. \citet{szymanski2009elucidating} utilized fluorescence correlation spectroscopy (FCS) of proteins immersed in crowded dextran solutions and reported a distribution of MSD coefficients with average $\langle \alpha \rangle \approx 0.82$, and compared this experimental finding with recovered distributions of MSD coefficients for simulated fractional Brownian motion (fBm), obstructed diffusion (OD), and CTRWs. Their findings indicated that the recovered distributions for fBm and OD matched the experiments, while the recovered distributions for CTRW-induced diffusion, with average $\langle \alpha \rangle \approx 0.59$ did not agree well with the data due to ergodicity breaking. \citet{weber2010bacterial} studied the subdiffusion of bacterial chromosomal loci in viscoelastic cytoplasm and further concluded that fractional Langevin motion to be more likely than CTRW and OD, due to the presence of ergodicity and a negative velocity auto correlation function. Regarding polymers, \citet{Wong2004} studied the thermal motion of colloidal tracer particles in entangled actin filament (F-actin) networks, under different concentrations and network mesh sizes. They observed a sub-diffusive behavior when the tracer particles radius were comparable to the network mesh size, and suggested that such anomalous behavior happens due intermittent caging behavior, followed by sudden infrequent jumps with a power-law distribution of caging times $\tau_c$ in the form $P(\tau_c) = \tau^{-1.33}_c$.

\subsection{State of the art: Anomalous materials modeling}\label{sec:material-sota}

As observed in Section \ref{sec:material-evidence}, experimental evidence suggests that complex material behavior may possess more than a single power-law scaling in the viscoelastic regime, particularly in multi-fractal structures, which are characteristic of cells \cite{Stamenovic2007} and biological tissues \cite{vincent2012structural}, due to their complex, hierarchical and heterogeneous microstructure. For such cases, a single SB element is not sufficient to capture the observed behavior, even if linear viscoelasticity holds. Furthermore, material nonlinearity due to large strains and additional physics such as plasticity, damage and failure require more advanced rheological models, which could have full or partial fractional nature. In this section we refer to a class of fractional models in the literature, classified by rheology type and nature of the corresponding FDEs. We acknowledge that rheology is a vast field with a large number of different types of material behavior, and here we limit our review to visco-elasto-plasticity, damage mechanics and failure.

\subsubsection{Viscoelasticity}
\label{Sec:VE}

\paragraph{Linear viscoelasticity.} We start by introducing two natural extensions of the SB viscoelastic model through serial and parallel combinations. The first one is the fractional Kelvin-Voigt (FKV) model, which is given by a parallel combination of SB elements, and relates the stresses $\sigma(t)$ and strains $\varepsilon(t)$ in the following additive form \cite{Schiessel1993}:
\begin{equation}\label{eq:FKV}
\sigma(t) = E_1 \prescript{\text{C}}{0}{} \mathbb{D}^{\alpha_1}_t \varepsilon(t) + E_2 \prescript{\text{C}}{0}{} \mathbb{D}^{\alpha_2}_t \varepsilon(t) ,\quad t > 0, \quad \varepsilon(0) = 0,
\end{equation}
where the fractional orders are such that $0 < \alpha_1, \alpha_2 < 1$, and $E_1\,[Pa.s^{\alpha_1}]$, $E_2\,[Pa.s^{\alpha_2}]$ are the associated pseudo-constants. The corresponding relaxation function also assumes additive form of two SB elements:
\begin{equation}
G^{\text{FKV}}(t) := \frac{E_1}{\Gamma(1-\alpha_1)} t^{-\alpha_1} + \frac{E_2}{\Gamma(1-\alpha_2)} t^{-\alpha_2},
\end{equation}
where contrary to the scale-free relaxation behavior of a single SB element, since we assume $\alpha_2 > \alpha_1$, the FKV model possesses two time-scale dependent power-law regimes, given by $G^{\text{FKV}} \sim t^{-\alpha_2}$ as $t\to 0$ and $G^{\text{FKV}} \sim t^{-\alpha_1}$ as $t\to \infty$, which characterizes a transition from faster to slower relaxation regimes. We note that this quality allows the FKV model to describe materials that reach an equilibrium behavior for large times when $\alpha_1 \to 0$, which is intuitive from the mechanistic standpoint as one of the SB elements becomes a Hookean spring. 

Through a serial combination of SB elements, we obtain the fractional Maxwell (FM) model \cite{Jaishankar2013}, given by:
\begin{equation}\label{eq:FM}
\sigma(t) + \frac{E_2}{E_1} \prescript{\text{C}}{0}{}\mathbb{D}^{\alpha_2 - \alpha_1}_t \sigma(t) = E_2 \prescript{\text{C}}{0}{}\mathbb{D}^{\alpha_2}_t \varepsilon(t) , \quad t>0,
\end{equation}
with $0 < \alpha_1 < \alpha_2 < 1$, and two sets of initial conditions for strains $\varepsilon(0) = 0$, and stresses $\sigma(0) = 0$. We note that in the case of non-homogeneous initial conditions, there needs to be compatibility conditions \cite{Mainardi2011} between stresses and strains at $t=0$. The corresponding relaxation function for this building block model assumes the more complex  Miller-Ross form \cite{Jaishankar2013}:
\begin{equation}\label{eq:G_FMM}
G^{\text{FM}}(t) := E_1 t^{-\alpha_1} E_{\alpha_2-\alpha_1, 1-\alpha_1}\left(-\frac{E_1}{E_2} t^{\alpha_2 - \alpha_1}\right),
\end{equation}
where $E_{a,b}(z)$ denotes the two-parameter Mittag-Leffler function, defined as \cite{Mainardi2011}:
\begin{equation}
  E_{a,b}(z) = \sum^\infty_{k=0} \frac{z^k}{\Gamma(a k + b)},\quad \text{Re}(a) > 0,\quad b \in \mathbb{C},\quad z \in \mathbb{C}.
\end{equation}
Interestingly, the presence of a Mittag-Leffler function in (\ref{eq:G_FMM}) produces a stretched exponential relaxation for smaller time-scales and a power-law behavior for larger time-scales. The asymptotic behaviors are given by $G^{\text{FM}} \sim t^{-\alpha_1}$ as $t\to 0$ and $G^{\text{FM}} \sim t^{-\alpha_2}$ as $t \to \infty$, indicating that, contrary to the FKV model, the FM model has a constitutive transition from slower-to-faster relaxation. We refer the reader to \cite{McKinley2013,bonfanti2020fractional} for a number of applications of the aforementioned models. Additionally, we notice that both FKV and FM models are able to recover the SB element with a convenient set of pseudo-constants, or naturally reveal the necessity of standard rheological elements according to available data. Furthermore, we also outline more complex building block models that produce more flexible responses, including three to four fractional orders, such as the fractional Kelvin-Zener (FKZ), fractional Poynting-Thomson (FPT), and fractional Burgers (FB) models. We refer the reader to \cite{bonfanti2020fractional,Schiessel1993} for more details on such models.

\paragraph{Numerical Discretization.} A well known numerical scheme to discretize the time-fractional Caputo derivatives of order $0 < \alpha < 1$ in \eqref{eq:FKV} and \eqref{eq:FM} is the implicit L1-difference scheme by \citet{lin2007finite}. Let points on a uniform time-grid be defined as $t_n = n \Delta t$ with $n=0,\,1,\,\dots,\,N$ time steps of size $\Delta t$. The discrete time-fractional Caputo derivative of a function $u(t)$ evaluated at $t = t_{n+1}$ is given by:
\begin{equation}
\label{eq:L1FDM}
{}^{\text{C}}_0{\mathbb{D}}_t^\alpha u(t) \Big|_{t=t_{n+1}} =
\frac{1}{\Gamma(2-\alpha)} \sum_{j=0}^{n} d_j \frac{u_{n+1-j}-u_{n-j}
}{\Delta t^\alpha} + r^{n+1}_{\Delta t},
\end{equation}
where $r^{n+1}_{\Delta t} \le C_u \Delta t^{2-\alpha}$ with the constant $C_u$ only depending on $u(t)$, and the convolution weights $d_j :=
(j+1)^{1-\alpha}
- j^{1-\alpha}, j=0,1,\dots,n$. The above expression can be rewritten and
approximated as:
\begin{equation*}
{}^{\text{C}}_0{\mathbb{D}}_t^\alpha   u(t)\Big|_{t=t_{n+1}} \approx
\frac{1}{\Delta t^\alpha \Gamma(2-\alpha)} \left[ u_{n+1} - u_n +
\mathcal{H}^{\alpha}u \right],
\end{equation*}
where the so-called \textit{history term} $\mathcal{H}^{\alpha}u$ is given by:
\begin{equation} \label{Eq:History}
\mathcal{H}^{\alpha}u = \sum_{j=1}^{n} d_j \left[ u_{n+1-j}-u_{n-j} \right].
\end{equation}
We note that although the above discretization is of practical and simple implementation, there exist many sophisticated numerical methods for fractional Cauchy equations that employ faster schemes, and also address non-smooth, nonlinear and stiff problems. We also emphasize that employing the kernel $G(t)$ into the Boltzmann representation for the aforementioned models may be impractical, since one would need other specialized numerical methods that are model-dependent, and would require evaluations of Mittag-Leffler functions.

\paragraph{Nonlinear Viscoelasticity.} Fractional linear viscoelastic models are suitable candidates to describe the anomalous dynamics of a number of materials undergoing small strains. However, under large strains, material nonlinearities induce stress/strain dependencies on the relaxation behavior. One alternative to incorporate such nonlinearity is through quasi-linear viscoelasticity (QLV) \cite{fung2013biomechanics}, which replaces $G(t)$ by a multiplicative decomposition between a reduced relaxation function $g(t)$ and an instantaneous, nonlinear elastic tangent response:
\begin{equation}\label{eq:QLV}
		\sigma(t,\varepsilon) = \int^t_0 g(t-s)\frac{\partial \sigma^e(\varepsilon)}{\partial \varepsilon} \dot{\varepsilon}\,ds,
\end{equation}
with $\sigma^e(\varepsilon)$ and $g(0^+) = 1$. Fractional approaches to QLV were developed by \citet{Doehring2005} for arterial valve cusp and by \citet{Craiem2008} for arterial wall viscoelasticity. In the latter, a reduced, fractional Kelvin-Voigt-type relaxation function $g(t) = C + Dt^{-\alpha}$ was employed, with pseudo-constant $D\,[s^\alpha]$, and nonlinear exponential form $\sigma^{e}(\varepsilon) = A \left(e^{B \varepsilon} - 1\right)$, with constant $A\,[Pa]$. Therefore, the fractional QLV formulation is able to capture not only linear/nonlinear instantaneous stress response, due to the rearrangement and alignment of fibers with the load direction, but also the anomalous power-law relaxation of the fractal microstructure. We also mention nonlinear models that take into account the Mittag-Leffler-type relaxation dynamics, such as the fractional QLV model in \cite{Doehring2005} and the fractional K-BKZ model introduced by \citet{Jaishankar2014}.

\subsubsection{Visco-elasto-plasticity}
\label{sec:plasticity}
Several works employed fractional calculus to account for the visco-plastic regimes of several classes of materials. We outline three of them: time-fractional, space-fractional and stress-fractional.

Time-fractional approaches focus on introducing memory effects into internal variables \citep{Suzuki2016,Xiao2017}, and consequently modeling power-laws in both viscoelastic and visco-plastic regimes. This is of interest for polymers, cells, and tissues. In this context, fractional visco-elasto-plastic models provide a constitutive interpolation between rate-independent plasticity and Perzyna's visco-plasticity by introducing a SB model acting the plastic regime \citep{Suzuki2016}, and utilizes a rate-dependent yield function of the form
\begin{equation}
f(\sigma, q) := |\sigma| - \left[\sigma^Y + K {}^C_0{\mathcal{D}}_t^{\alpha_K} q(t) + H q(t)\right], \quad 0 < \alpha_K < 1,
\end{equation}
where $\sigma^Y$ and $q$ denote, respectively, the yield stress and the accumulated plastic strain, with pseudo-constant $K\,[Pa.s^{\alpha_K}]$ and Hookean constant $H$. The above form for the yield function was later proved to be thermodynamically consistent in a further extension of the model to account for continuum damage mechanics \cite{suzuki2021thermodynamically}. 

A three-dimensional space-fractional approach to elastoplasticity was also developed by \citet{sumelka2014application} to account for spatial nonlocalities. The model is based on rate-independent elastoplasticity, and nonlocal effects are accounted for through a fractional continuum mechanics approach, where the strains are defined by a space-fractional Riesz-Caputo derivative of displacements $u(x)$ in the form
\begin{equation}
    \prescript{\text{RC}}{a}{\mathbb{D}}^\alpha_b u(x)= \frac{\Gamma(2-\alpha)}{2} \left( \prescript{\text{C}}{a}{\mathbb{D}}^\alpha_x u(x) + (-1)^n \prescript{\text{C}}{x}{\mathbb{D}}^\alpha_b u(x)\right),
\end{equation}
for left- and right-sided fractional Caputo derivatives \citep{sumelka2014application} with $n = \lceil \alpha \rceil$.

Finally, stress-fractional models for plasticity have found applicability in soil mechanics and geomaterials that follow non-associated plastic flow \cite{Sumelka2014VP,Sumelka2019Soil}, i.e., the yield surface expansion in the stress space does not follow the usual normality rule, and may be non-convex. The work by \citet{Sumelka2014VP} proposed a three-dimensional fractional viscoplastic model, where a fractional flow rule with order $0 < \alpha < 1$ in the stress domain naturally models non-associative plasticity. Interestingly, this model recovers the classical Perzyna visco-plasticity as $\alpha \to 1$, and the effect of the fractional flow rule can be a compact descriptor of microstructure anisotropy. Recently, a similar stress-fractional model was developed \cite{Sumelka2019Soil}, and successfully applied to soils under compression. We refer the reader to the detailed review work by \citet{sun2018new} for a review of uses fractional calculus in plasticity.

\subsubsection{Damage mechanics, ageing and failure}\label{sec:damage}

There have also been recent efforts to include damage, ageing and failure effects into fractional calculus frameworks. Existing formulations are focused on either adding classical failure frameworks into existing fractional constitutive laws, or by developing fractional failure mechanisms. Here, we mostly focus on the latter and start with the work by \citet{caputo2015damage}, that developed a variable-order viscoelastic model in the form:
\begin{equation}
  \sigma(x,t) = g(\alpha(x,t)) A(x) {}^C_{t_0}{\mathbb{D}}_t^{\alpha(x,t)} \varepsilon(x,t),
\end{equation}
where $g(\alpha(x,t)) := (\alpha_C - \alpha(x,t))^2/4$ denotes a material degradation function with critical damage $\alpha_C$, $A(x)$ represents a space-dependent pseudo-property, and $0 < \alpha(x,t) < \alpha_{C}$ is the variable fractional order, also interpreted here as damage. The variable-order Caputo derivative is defined in \eqref{eq:variable_order_caputo}. Interestingly, this mixed interpretation for $\alpha(x,t)$ makes it a multi-physics descriptor for anomalous damage, viscosity, and material ageing. The evolution of $\alpha(x,t)$ is described by an integer-order phase-field equation, and the resulting model is proved to be thermodynamically admissible.

A key aspect to develop failure models relies on consistent forms of damage energy release rates, i.e., on obtaining the compatible operator for the loss of elastic energy, which is a nontrivial task even for the simplest fractional constitutive law (\ref{eq:SB-caputo}). This has been achieved by employing the concept of fractional free-energy densities \cite{Lion1997,Fabrizio2013,Alfano2017}. \citet{Alfano2017}
developed a cohesive zone, damaged fractional viscoelastic Kelvin-Zener model, and studied the influence of integer and fractional damage energy release rates on damage evolution. In this case, integer-order energy loss considers Hookean-type rheology to compute the damage energy release rates, which may be justified when Hookean elements are present in the viscoelastic constitutive law, but incompatible for fully-fractional cases (an arrangement of Scott-Blair elements). The corresponding free-energy for the SB element is given by:
\begin{equation}\label{eq:free-energy-SB-explicit}
  \psi^{SB}(t)  = \frac{E}{2\Gamma(1-\alpha)}\int^t_0 \int^t_0 \left(2t-\tau_1 - \tau_2\right)^{-\alpha}\dot{\varepsilon}(\tau_1)\dot{\varepsilon}(\tau_2)\,d\tau_1 d\tau_2,
\end{equation}
with $0 < \alpha < 1$, which clearly carries a power-law behavior over time. Among their findings, the authors obtained a rate-dependence of the fracture energy in terms of the fractional-order $\alpha$, opening interesting directions towards failure of anomalous viscoelastic media such as polymers. In \cite{suzuki2021thermodynamically} this idea was extended to plasticity, and a fractional visco-elasto-plastic model with memory-dependent damage was developed, with isotropic damage evolution $0 \le D(t) < 1$ given by Lemaitre's approach \cite{Lemaitre1996}:
\begin{equation}\label{eq:lemaitre-damage}
  \dot{D}(t) = \frac{\dot{\gamma}(t)}{1-D(t)}\left(-\frac{Y^{ve}(t)}{S}\right)^s,
\end{equation}
with material damage parameters $s,\,S \in \mathbb{R}^+$, plastic slip $\dot{\gamma}$ and damage energy release rate $Y^{ve}(t) = -\psi^{SB}(t)$. We note that although \eqref{eq:lemaitre-damage} is a nonlinear ODE, the memory is introduced through the power-law form of $Y^{ve}$ \eqref{eq:free-energy-SB-explicit}. In this formulation, the viscoelastic and visco-plastic fractional orders introduce a competition between rate-dependent hardening and damage-induced softening, which could open interesting directions for modeling localized hardening in failing anomalous media. Sumelka \textit{et al.} in \cite{sumelka2020modelling} also developed the idea of memory-dependent damage for soft materials through a stress-driven time-fractional hyperelastic damage model, with evolution equation in the following fractional nonlinear Cauchy form:
\begin{equation}
  {}^C_{t-l_t}{\mathbb{D}}_t^{\alpha} D(x,t) = \frac{1}{T^\alpha}\Phi\big\langle \frac{I_D}{\tau_D} - 1 \big\rangle,
\end{equation}
where $\Phi$ represents an overstress function in terms of a stress intensity $I_D$, threshold stress $\tau_D$ for damage evolution, and a ramp function in Macaulay notation $\langle . \rangle$. The memory length is driven by a time scale $l_t$, which was taken as a fraction of the total time $T$. This model was applied with an Ogden hyperelastic law to patient-specific three-dimensional abdominal aortic aneurysm (AA) for critical zone identification, with obtained fractional order $\alpha = 0.75$.

Additional work on variable-order models in the context of fractional damage, ageing and failure include the following contributions. In \citet{beltempo2018fractional} a variable-order viscoelastic creep model was developed, where the evolution of the fractional order $\alpha(t)$ dictates the process of concrete ageing. The variable-order viscoelastic model developed in \citet{meng2019variable} employed a piecewise constant order followed by two linear decreasing functions for $\alpha(t)$ successfully described the initial viscoelasticity, softening and hardening of amorphous glassy polymers under compression. Finally, variable-order operators also proved to be useful mathematical tools to determine the onset of fracture. \citet{patnaik2021variable} employed a variable fractional-order activation function for damage, where the sharp power-law activation threshold induced by the fractional operator was successfully employed to determine crack propagation and branching of brittle materials. We refer the reader to the recent review works on the use of variable-order \cite{patnaik2020applications} and distributed-order \cite{ding2021distributed} fractional models in viscoelasticity and structural mechanics. In the distributed-order case, fractional derivatives are integrated with respect to a distribution of fractional orders within a certain range of values.

\subsection{Future directions in modeling anomalous materials}\label{sec:material-future}

Although there exists a large spectrum of fractional models in the context of materials science, solid mechanics and rheology, these models are mostly characterized by constant-order fractional operators, for which a significant number of fast time-integration schemes is available. Yet, there is still a need for efficient numerical methods for variable- and distributed-order operators. In fact, although fractional models lead to a compact physical description with reduced number of material parameters, the computational cost is still high when calibrating the models with large experimental data sets. Furthermore, although there exist a increasing number of distributed-order operators in the context of viscoelasticity, structural mechanics, and anomalous diffusion \cite{ding2021distributed}, further validation against experimental data is needed. 

We point out interesting research directions that could involve the use of variable- and distributed-order differential equations in the multi-scale modeling of materials. Recently, nano-scale simulation studies on trapping of nano-particles in hydrogel networks indicated a time-temperature dependency of the MSD in the evolution of anomalous diffusion regimes, where a subdiffusion regime has been found to be of transitional nature at intermediate time scales, with ballistic/normal diffusion dynamics for short/long time scales \cite{wang2020controlled}. This motivates the study of variable-order models in time to compactly describe the macroscopic rheological evolution of such polymer networks. Furthermore, the observation of distributions of power-law scaling parameters in micro-rheology creep experiments on cells \cite{Bonadkar2016,hecht2015imaging} indicate the presence of microstructure-induced randomness in rheological response. In this sense, distributed-order models may arise as interesting approaches to naturally incorporate the stochastic parametric data into the differential operator \cite{kharazmi2017petrov}.

\section{Conclusion}\label{sec:conclusion}
In this work we reviewed fundamental concepts of anomalous transport processes and provided the mathematical and statistical background for understanding them. We then selected three applications where the use of fractional models has experienced dramatic growth and improvement. This set of applications was chosen at our discretion and is, by no means, complete. In fact, several other scientific and engineering fields are currently benefiting from fractional modeling (see, e.g., image processing, finance, machine learning algorithms and many others). However, based on the amount of literature, significance of the applications, and variety of fractional models for their descriptions, we believe that subsurface transport, turbulence, and anomalous materials allowed us to provide insights on the several uses and benefits of fractional modeling. Furthermore, these applications are still the subject of very active fractional research. Finally, given the recent advances in high-performance computing and machine learning, we believe it is now the best time to promote and increase the usability of fractional and nonlocal models for those applications that cannot be adequately described by the classical PDE models.

\section*{Acknowledgements}
Marta D'Elia and Mamikon Gulian were partially supported by the U.S. Department of Energy, Office of Advanced Scientific Computing Research under the Collaboratory on Mathematics and Physics-Informed Learning Machines for Multiscale and Multiphysics Problems (PhILMs) project. Marta D'Elia was also supported by the Sandia National Laboratories Laboratory-directed Research and Development (LDRD) program, project 218318. Mamikon Gulian was also supported by John von Neumann fellowship at Sandia National Laboratores. Mohsen Zayernouri and Jorge L. Suzuki were supported by the ARO Young Investigator Program (YIP) award (W911NF-19-1-0444), the National Science Foundation award (DMS-1923201), and the MURI/ARO grant (W911NF-15-1-0562).

Sandia National Laboratories is a multimission laboratory managed and operated by National Technology and Engineering Solutions of Sandia, LLC., a wholly owned subsidiary of Honeywell International, Inc., for the U.S. Department of Energy's National Nuclear Security Administration contract number DE-NA0003525. This paper, SAND2021-11291 R, describes objective technical results and analysis. Any subjective views or opinions that might be expressed in the paper do not necessarily represent the views of the U.S. Department of Energy or the United States Government.

\section*{Conflict of Interest}
On behalf of all authors, the corresponding author states that there is no conflict of interest. 

\bibliographystyle{spbasic}
\bibliography{reference,SC_references}


\end{document}